\documentclass[a4paper,11pt,reqno]{amsart}

\usepackage{amsmath}
\usepackage{amssymb}
\usepackage{amsfonts}
\usepackage{graphicx}
\usepackage{mathtools}
\usepackage{bbm}
\usepackage[colorlinks]{hyperref}

\renewcommand\eqref[1]{(\ref{#1})} 
\graphicspath{ {images/} }
\title[Hardy--Rellich inequalities for Baouendi--Grushin operators] {Extremizer Stability of Higher-order Hardy--Rellich inequalities for Baouendi--Grushin vector fields}

\author[A. Banerjee]{Avas Banerjee}
	\address{Avas Banerjee
    \endgraf
    Theoretical Statistics and Mathematics Unit
    \endgraf
    Indian Statistical Institute, Delhi Center
    \endgraf
    S.J.S. Sansanwal Marg, New Delhi, Delhi 110016
    \endgraf
    India
	}
	\email{avas24r@isid.ac.in}
    
\author[R. Basak]{Riju Basak}
	\address
    {Riju Basak
    \endgraf
    Department of Mathematics
    \endgraf
    National Taiwan Normal University
    \endgraf
    No. 88, Section 4, Tingzhou Road, Wenshan District, Taipei City 116
    \endgraf
    Taiwan
	}
	\email{rijubasak52@ntnu.edu.tw}
    
\author[P. Roychowdhury]{Prasun Roychowdhury}
\address
	{Prasun Roychowdhury
	\endgraf
	Department of Mathematics
    \endgraf
    Indian Institute of Technology Hyderabad
	\endgraf
	Kandi, Sangareddy, Telangana, 502285
	\endgraf
	India
    }
	\email{prasunrc@math.iith.ac.in}


\date{\today}

\theoremstyle{plain}
\newtheorem{theorem}{Theorem}[section]
\newtheorem{proposition}{Proposition}[section]
\newtheorem{lemma}{Lemma}[section]
\newtheorem{corollary}{Corollary}[section]
\newtheorem{remark}{Remark}[section]

\numberwithin{equation}{section} \allowdisplaybreaks

\usepackage[text={6.8in,10in},centering]{geometry}
\newcommand{\rn}{\mathbb{R}^n}
\newcommand{\re}{\mathbb{R}}
\newcommand{\rno}{\mathbb{R}^{n+1}}
\newcommand{\rko}{\mathbb{R}^{q+1}}
\newcommand{\rmo}{\mathbb{R}^{m+1}}
\newcommand{\m}{\mathbb{R}^{n+1}}
\newcommand{\gradg}{\nabla_G}
\newcommand{\rgradg}{\nabla_{\varrho,G}}
\newcommand{\gradgp}{\nabla_{G'}}
\newcommand{\rgradgp}{\nabla_{\tilde{\varrho},G'}}
\newcommand{\lapg}{\mathcal{L}_G}
\newcommand{\rlapg}{\mathcal{L}_{\varrho,G}}
\newcommand{\dx}{\:{\rm d}x}
\newcommand{\dy}{\:{\rm d}y}
\newcommand{\dsn}{\:{\rm d}\Omega}
\newcommand{\dsnp}{\:{\rm d}\Omega'}
\newcommand{\dt}{\:{\rm d}t}
\newcommand{\ds}{\:{\rm d}s}
\newcommand{\dr}{\:{\rm d}\varrho}
\newcommand{\drt}{\:{\rm d}\tilde{\varrho}}
\newcommand{\dph}{\:{\rm d}\phi}
\newcommand{\dw}{\:{\rm d}w}

\newcommand{\pd}{\partial_\varrho}
\newcommand{\pdp}{\partial_{\tvr}}
\newcommand{\br}{B_R^Q(o)}
\newcommand{\brp}{B_R^{Q'}(o')}
\newcommand{\ba}{B_1^Q(o)}
\newcommand{\sn}{\Omega}
\newcommand{\dv}{\dx\dt}
\newcommand{\vr}{\varrho}
\newcommand{\tvr}{\tilde{\varrho}}
\newcommand{\btvr}{\bar{\varrho}}

\begin{document}
\begin{abstract}
In this paper, our primary aim is to establish the extremizer stability of $L^p$-Rellich and Hardy–Rellich inequalities in the setting of Baouendi–Grushin vector fields. As a consequence, we obtain improved forms of these inequalities and derive an identity relating the subcritical and critical Hardy inequalities, thereby demonstrating their equivalence. These improvements are obtained through a suitably defined distance from extremizers. In the higher-order setting, we derive Hardy–Rellich type inequalities involving all radial operators in the Grushin framework and prove that all resulting constants are sharp. Finally, for the $L^2$-higher-order cases, we compute exact remainder terms by establishing identities rather than inequalities.
\end{abstract}
\maketitle

\section{Introduction}
Hardy inequalities are fundamental tools with wide-ranging applications in analysis, mathematical physics, spectral theory, geometry, and quantum mechanics. In particular, their sharp forms play a crucial role in establishing existence and non-existence results for various classes of partial differential equations. Due to their significance and versatility, these inequalities have been the subject of extensive research. The Hardy inequality was first introduced in the seminal paper \cite{hardy}, marking the beginning of its profound impact on mathematical analysis. The development of the classical Hardy inequality has a rich history; we refer the reader to~\cite{kuf} (see also the preface of~\cite{rs-book}). The radial form of the classical Hardy inequality in Euclidean space $\rn$ (with dimension $n \geq 2$) can be stated as follows: suppose $1 \leq p < n$, then
\begin{equation}\label{Hardy-Leray}
\int_{\rn}\left| \frac{x}{|x|} \cdot \nabla u(x)\right|^{p}\dx\geq \left(\frac{n-p}{p}\right)^p\int_{\rn} \frac{|u(x)|^{p}}{|x|^{p}}\dx
\end{equation}
holds for all $u \in C_{0}^\infty(\rn \setminus \{\vec{0}\})$. This range of $p$ is referred to in the literature as the \emph{subcritical} case of the Hardy inequality. Moreover, it is known that for $1 < p$, the constant $\big(\frac{n-p}{p}\big)^p$ is sharp and never attained, except by the trivial function. In the critical regime $p = n$, the inequality \eqref{Hardy-Leray} fails for any positive constant on the left-hand side. However, in \cite[Theorem~1.1]{ii}, Ioku and Ishiwata established a scale-invariant \emph{critical} version of the Hardy inequality by introducing a logarithmic correction term on the Euclidean ball $B_R:=\{x\in\rn \,: \, |x|<R\}$ with radius $R>0$,  as follows
\begin{equation}\label{Hardy-critical-scale}
 \int_{B_R}\left| \frac{x}{|x|} \cdot \nabla u(x)\right|^{n}\dx\geq \left(\frac{n-1}{n}\right)^n\int_{B_R} \frac{|u(x)|^{n}}{|x|^{n}\left(\ln \frac{R}{|x|}\right)^n}\dx
\end{equation}
for all $u \in C_{0}^\infty(\rn \setminus \{\vec{0}\})$ with $n \geq 2$. Once again, the constant $\big(\frac{n-1}{n}\big)^n$ is sharp and never attained.

Owing to their wide-ranging applications, considerable effort has been devoted to extending Hardy-type inequalities beyond the classical Euclidean setting. In particular, there has been growing interest in functional inequalities and identities of Hardy type associated with the sub-Laplacian, defined as the sum of squares of vector fields on $\mathbb{R}^{n+m}$. In this article, we focus on inequalities such as \eqref{Hardy-Leray}, \eqref{Hardy-critical-scale}, and their higher-order counterparts, specifically, the Rellich inequality for the \emph{Baouendi--Grushin operator}. Our approach involves establishing these inequalities and examining their refinements through stability analysis, including the explicit computation of remainder terms in certain cases. In particular, we consider the operator
\begin{equation*}
 \lapg^{\alpha}  =  \Delta_x  +  |x|^{2\alpha} \Delta_y,
\end{equation*}
where $x \in \rn$, $y \in \mathbb{R}^m$, $n \in \mathbb{N}$, $m \in \mathbb{N}$, $|x| = \sqrt{x_1^2 + \cdots + x_n^2}$, $\alpha \in \mathbb{R}$, and $\Delta_x$ and $\Delta_y$ denote the classical Laplacians with respect to the $x$ and $y$ variables, respectively. The operator $\lapg^{\alpha}$ is commonly referred to as the Grushin (or Baouendi--Grushin) operator. It is immediately evident that this operator is not uniformly elliptic due to the presence of the term $|x|^{2\alpha}$, which leads to degeneracy on the submanifold $\{(\vec{0},y)\in \mathbb{R}^{n+m}\, :\, y\in \mathbb{R}^m\}$. This sub-Laplacian operator was first introduced in a preliminary form by Baouendi \cite{boun} and Grushin \cite{Gr1, Gr2}, and it serves as a bridge between elliptic and non-elliptic theory in the context of partial differential equations.

It is worth noting that when $\alpha \in \mathbb{N}$, the Grushin operator $\lapg^{\alpha}$ can be represented as a finite sum of squares of smooth \emph{Grushin vector fields}. In this setting, H\"ormander’s celebrated theorem \cite{Hor1} applies, as the associated vector fields satisfy the finite rank (or bracket-generating) condition, thereby ensuring that the operator is hypoelliptic. However, in the general case where $\alpha$ is not a natural number, H\"ormander’s condition becomes inapplicable due to the insufficient smoothness of the underlying vector fields. Nonetheless, a substantial body of literature addresses operators generated by non-smooth vector fields, where the classical H\"ormander framework does not apply. Among these contributions, we particularly highlight the influential work of Franchi and Lanconelli \cite{BF1}. The analysis of the Grushin operator $\lapg^{\alpha}$ in this broader context is intricate, and in the present work, we shall primarily focus on the specific case $\alpha = 1$ and $m = 1$. From this point onward, we shall refer to the Baouendi--Grushin vector field simply as the ``Grushin vector field”.

We begin by explicitly describing the vector fields and the Grushin operator associated with this setting. Let $(x, t) \in \mathbb{R}^n \times \mathbb{R}$ with $n \geq  1$. For $i = 1, \ldots, n$, consider the vector fields
\begin{align*}
 X_i = \frac{\partial}{\partial x_i}, \quad Y = |x| \frac{\partial}{\partial t}.
\end{align*}
Equipped with these vector fields, the space $\mathbb{R}^{n+1}$ is referred to as the Grushin space. The Grushin gradient is defined as
\begin{align}\label{grad}
 \gradg = (X_1, \ldots, X_n, Y) = (\nabla_x, |x|\, \partial_t),
\end{align}
where $\nabla_x$ denotes the standard Euclidean gradient on $\mathbb{R}^n$. The corresponding Grushin operator is given by
\begin{align}\label{lap}
 \lapg = \Delta_x + |x|^2 \frac{\partial^2}{\partial t^2},
\end{align}
where $\Delta_x$ is the Laplacian in the $x$-variables. Moreover, for a vector field $X$, if we define the Grushin divergence by $\text{div}_G X = \gradg \cdot X$, then the operator can be expressed in divergence form as $\lapg = \text{div}_G \cdot \gradg$.

We also introduce a natural dilation on $\mathbb{R}^{n+1}$, defined by
\begin{align*}
 \delta_\lambda(x, t) := (\lambda x, \lambda^2 t),
\end{align*}
for $\lambda > 0$ and $(x,t) \in \mathbb{R}^{n+1} \setminus \{o\}$, where $o = (\vec{0}, 0)$. A direct computation shows that $\lapg(f \circ \delta_\lambda) = \lambda^2 (\lapg f) \circ \delta_\lambda$, illustrating the homogeneous nature of the operator. Furthermore, for any Lebesgue measurable set $U \subset \mathbb{R}^{n+1}$, we have $|\delta_\lambda U| = \lambda^Q |U|$, where $Q = n + 2$ is the homogeneous dimension of the Grushin space.

On $\mathbb{R}^{n+1}$, we define the function
\begin{align*}
 \varrho(x,t) := \big(|x|^4 + 4t^2\big)^{1/4},
\end{align*}
that serves as a homogeneous norm compatible with the dilation $\delta_\lambda$. This function is smooth on $\mathbb{R}^{n+1} \setminus \{o\}$, positive definite, and satisfies the homogeneity property $\varrho(\delta_\lambda(x,t)) = \lambda \varrho(x,t)$ for all $\lambda > 0$. Notably, $\varrho$ is closely related to the fundamental solution of $\lapg$ (see \cite[Proposition~2.1]{G1}) and will play a key structural role in our analysis.

Finally, we define the $\varrho$-gauge ball in $G = \mathbb{R}^{n+1}$ of radius $R > 0$ centered at $o$ by
\begin{align*}
 \br=\{(x,t)\in \rno \, : \varrho(x,t)<R\}.
\end{align*}
Before presenting our main results, we note that $\rgradg$ and $\rlapg$ denote, respectively, the radial parts of the Grushin gradient \eqref{grad} and the Grushin operator \eqref{lap} (see \eqref{radial-op} in Section~\ref{prelim} for their full descriptions).

The Hardy inequality for the Grushin operator was first established by Garofalo \cite[Corollary~4.1]{G1} for the case $p=2$, and later extended to the general case $p>1$ by D’Ambrosio \cite[Theorem~3.1]{AL1}. Since then, numerous refinements and generalizations of the Hardy inequality for Grushin-type vector fields have been explored in various directions; we refer, without claiming completeness, to \cite{IK1, AS, sj, KY, ye, RS, cow, bad}.

Subsequently, Laptev, Ruzhansky, and Yessirkegenov \cite{ari-ruz-nur} investigated a version of the inequality involving the radial operator corresponding to the nondegenerate component of the Grushin gradient for $p=2$. For a result that incorporates the full radial part of the Grushin gradient, we refer to \cite[Corollary~3.2]{GJR}. More recently, D’Arca \cite[Corollary~4.3]{darca-1} established a weighted Hardy inequality for the radial Grushin gradient in the case $p \geq 2$. In a very recent work, Yessirkegenov and Zhangirbayev \cite{bad} developed refined general weighted $L^p$-Hardy type inequalities associated with Baouendi–-Grushin vector fields. In particular, \cite[Corollary~4.2]{bad}, with the choices $\gamma=1$, $\beta=p$, and $\alpha=2p+\beta$, yields the weighted Hardy inequality for the radial Grushin gradient in the $L^p$ framework across the entire range of $p$. The result can be stated as follows: Let $G$ be a Grushin space of homogeneous dimension $Q$, with $Q \geq 3$. Let $1 < p < \infty$, and $\beta \in \mathbb{R}$ satisfy $\beta < Q - p$. Then, for all $u \in C_0^\infty(\mathbb{R}^{n+1})$, the following inequality holds
\begin{align}\label{sub-hardy-eq-intro}
 \int_{\mathbb{R}^{n+1}} \frac{|\rgradg u(x,t)|^p}{\varrho^{\beta}(x,t)}\dx\dt \geq \left(\frac{Q - p - \beta}{p}\right)^p \int_{\mathbb{R}^{n+1}} \frac{|u(x,t)|^p |\gradg \varrho(x,t)|^p}{\varrho^{p + \beta}(x,t)}\dx\dt, 
\end{align}
where the constant $\big(\frac{Q - p - \beta}{p}\big)^p$ is sharp. The result in \eqref{sub-hardy-eq-intro} with $\beta = 0$ can also be viewed as a counterpart to \eqref{Hardy-Leray} in the Grushin setting.

To the best of our knowledge, the study of critical Hardy inequalities in the context of Grushin operators was initiated in \cite[Theorem~3.1, (2)]{da-crt}, where the authors employed the divergence theorem. However, the formulation of the critical Hardy inequality involving the radial Grushin gradient appears in the very recent work \cite{bad}. In particular, the result presented in \cite[Corollary~4.12]{bad}, with the choices $\gamma=1$ and $\alpha=-b$, can be stated as follows: Let $G$ be a Grushin space of homogeneous dimension $Q$, with $Q \geq 3$, and let $\br$ denote the $\varrho$-gauge ball of radius $R > 0$. Suppose that $1 < p < \infty$, $b \in \mathbb{R}$, and $b > 1$. Then, for any $u \in C_0^\infty(\br)$, the following inequality holds
\begin{align}\label{lim-tw-wg-hardy-eq-intro}
 \int_{\br} \frac{|\rgradg u(x,t)|^p}{\big(\ln \frac{R}{\varrho(x,t)}\big)^{b-p}} \dx\dt \geq \left(\frac{b-1}{p}\right)^p \int_{\br} \frac{|u(x,t)|^p |\gradg \varrho(x,t)|^p}{\varrho^{p}(x,t)\big(\ln \frac{R}{\varrho(x,t)}\big)^b} \dx\dt,
\end{align}
where the constant $\big(\frac{b-1}{p}\big)^p$ is sharp. Furthermore, the inequality \eqref{lim-tw-wg-hardy-eq-intro} may be compared to the Euclidean critical Hardy inequality presented in \eqref{Hardy-critical-scale}. The inequality \eqref{lim-tw-wg-hardy-eq-intro} was established by selecting an appropriate weight within the framework of general weighted $L^p$-Hardy type inequalities. Alternatively, a more general weighted version of \eqref{lim-tw-wg-hardy-eq-intro} can be obtained through a limiting analysis of certain two-weighted Hardy inequalities (see Theorem~\ref{two-wg-hardy}), combined with a precise implementation of the polar coordinate structure introduced in Section~\ref{prelim} (see Corollary~\ref{lim-tw-wg-hardy} with $a=p$).

Analogously to the way the Hardy inequalities \eqref{Hardy-Leray} and \eqref{Hardy-critical-scale} differ significantly in the Euclidean setting in terms of their form, scaling behavior, and sharp constants, the inequalities \eqref{sub-hardy-eq-intro} and \eqref{lim-tw-wg-hardy-eq-intro} also exhibit fundamental distinctions within the Grushin framework. Notably, Sano and Takahashi \cite{st-cvpde} established that, in the Euclidean setting, the critical Hardy inequality on a ball is equivalent to the subcritical Hardy inequality on the whole space. Inspired by this remarkable correspondence, we investigate the relationship between critical and subcritical Hardy inequalities in the context of Grushin spaces. This analysis constitutes one of the main contributions of the present work. Our first key result is presented below.

\begin{theorem}\label{main-thm-1}
Let $G = \rno$ and $G' = \rmo$ be the Grushin spaces. Consider points $(x,t) \in \rno$ and $(y,s) \in \rmo$ in these spaces, and let $\varrho$ and $\tilde{\varrho}$ denote the respective homogeneous norms. Assume that $Q = n + 2 > Q' = m + 2 \geq 3$, and let $\brp$ denote the $\tilde{\varrho}$-gauge ball in $G'$ with radius $R > 0$ and center $o' = (\vec{0}, 0) \in \rmo$. Suppose that $\beta, b \in \mathbb{R}$ satisfy $\beta < Q - Q'$ and $Q' \leq b$. Then, for any $w \in C_0^\infty(\brp \setminus \{o'\})$, there exists a function $u \in C_0^\infty(\rno \setminus \{o\})$ such that the following identity holds
\begin{align*}
 \int_{\rno}&\frac{|\rgradg u(x,t)|^{Q'}}{\varrho^{\beta}(x,t)}\dx\dt  
 - \left(\frac{Q - Q' - \beta}{Q'}\right)^{Q'}\int_{\rno}\frac{|u(x,t)|^{Q'}|\gradg \varrho(x,t)|^{Q'}}{\varrho^{Q' + \beta}(x,t)}\dx\dt \\
 &= \frac{\omega_n}{\omega_m}
 \left(\frac{Q - Q' - \beta}{b - 1}\right)^{Q' - 1}\bigg[
 \int_{\brp} \frac{|\rgradgp w(y,s)|^{Q'}}{\big(\ln \frac{R}{\tilde{\varrho}(y,s)}\big)^{b - Q'}}\dy\ds \\
 &\quad - \left(\frac{b - 1}{Q'}\right)^{Q'}\int_{\brp} \frac{|w(y,s)|^{Q'}|\gradgp \tilde{\varrho}(y,s)|^{Q'}}{\tilde{\varrho}^{Q'}(y,s)\big(\ln \frac{R}{\tilde{\varrho}(y,s)}\big)^b}\dy\ds
 \bigg],
\end{align*}
where $\omega_n$ and $\omega_m$ denote the surface measures of the unit spheres in $\re^n$ and $\re^m$, respectively.
\end{theorem}

A related result has been investigated in the setting of homogeneous groups in \cite[Theorem~6.1]{rsy}, where the authors adopted the same methodology as in \cite{st-cvpde}, namely the polar coordinate decomposition. However, their analysis was restricted to radial functions, primarily due to the absence of an explicit polar coordinate structure in general homogeneous groups. In contrast, in the case of the Grushin space, the availability of an explicit polar coordinate framework, along with carefully detailed computations, allows us to extend the result to non-radial functions as well. Moreover, by applying Theorem~\ref{main-thm-1} and a suitable weighted Caffarelli--Kohn--Nirenberg inequality, we are able to further improve~\eqref{lim-tw-wg-hardy-eq-intro} for radial functions; this improvement is also discussed in detail.

The Hardy inequality \eqref{sub-hardy-eq-intro} for Grushin vector fields has been improved in various directions; see, for instance, the works \cite{DQN, IK1, KY, SY, YSK, ZHD, nch, al} and the references therein. In this paper, we focus on improving the inequality from the perspective of \emph{stability}. Functional inequalities serve as fundamental tools in the analysis of partial differential equations, particularly in establishing bounds and regularity properties of solutions. Moreover, quantitative stability estimates measure the deficit in the inequality in terms of the distance to the manifold of extremals, providing a finer understanding of the structure of extremizing sequences. Consequently, the stability of functional inequalities has attracted significant attention in recent years. Motivated by these developments, the remaining sections of this paper are devoted to the study of stability aspects of Hardy inequalities associated with the Grushin operator.

From that perspective, let us first define the nonnegative deficit term by
\begin{align*}
 \delta(u):= \int_{\rno}\frac{|\rgradg u(x,t)|^p}{\varrho^{\beta}(x,t)}\dx\dt- \left(\frac{Q-p-\beta}{p}\right)^p\int_{\rno}\frac{|u(x,t)|^p|\gradg \varrho(x,t)|^p}{\varrho^{p+\beta}(x,t)}\dx\dt \quad \forall u\in  C_0^\infty(\rno).
\end{align*}
Now, let us take the collection of smooth extremals
\begin{align*}
    \mathcal{G}:=\{u\in C_0^\infty(\rno)\, :\, \delta(u)=0\}.
\end{align*}
A central question concerning the stability estimate is how to quantify the \emph{Hardy deficit} $\delta(u)$ from below. Specifically, one asks whether the following inequality holds
\begin{align*}
 \delta(u) \geq \Phi(d(u, \mathcal{G})) \quad \text{ for all } u \in C_0^\infty(\mathbb{R}^{n+1}),
\end{align*}
where $\Phi$ is a nonnegative modulus of continuity function and $d(u, \mathcal{G})$ denotes an appropriate distance from $u$ to the set of smooth extremals $\mathcal{G}$. In recent years, the problem of establishing stability estimates has received considerable attention in the analysis of functional and geometric inequalities. Notable results have been obtained for Hardy inequalities \cite{cf-aihp, ban25, ban26}, Caffarelli--Kohn--Nirenberg inequalities \cite{cfll,dfll}, the Hardy--Littlewood--Sobolev inequality \cite{Car17,ex-1, clt26}, log-Sobolev and Moser-Onofri inequalities \cite{lslu, dolcjm, llh25}, and the Heisenberg--Pauli--Weyl uncertainty principle \cite{MV21, DN24, dgll}, among other related results \cite{llls, llr26, dlljfa, chenhei}. More recently, attention has shifted toward obtaining optimal quantitative lower bounds for the stability deficit. Significant progress has been achieved for Sobolev, fractional Sobolev, and Hardy--Littlewood--Sobolev inequalities through the development of refined analytical techniques, including rearrangement flow methods, duality arguments, and orthogonal decomposition techniques; see \cite{dolcjm, llh25, clt26, ex-1}. These advances have substantially strengthened the quantitative theory of geometric inequalities and provide further motivation for investigating analogous stability problems in the Grushin setting. However, a major challenge in this setup is $\mathcal{G}=\{0\}$, i.e., the absence of a nontrivial smooth extremal. Despite this, it is known that the sharpness of the Hardy constant $\big(\frac{Q-p-\beta}{p}\big)^p$ can be demonstrated via approximations using the profile $\varrho^{-(Q-p-\beta)/p}$. By carefully analyzing scaled versions of this function, we investigate the stability properties of the Hardy inequality. This type of extremizer-based stability analysis has been previously studied in the literature; see, for instance, \cite{sano-mia, RSY, st-die, rs-anf,rsy}. However, to our knowledge, this is the first time such an analysis has been carried out in the Grushin operator setting. Motivated by this perspective, we introduce a distance function that measures the deviation of a function from the family of scaled approximate profiles. This choice is consistent with the extremizer-based stability framework developed in the above works and is naturally adapted to the Grushin setting.

For a function  $u\in C_0^\infty(\m\setminus\{o\})$, $1<p<\infty$, $\beta\in\mathbb{R}$ and for a fixed real number $R>0$, we define the following distance function:
\begin{align*} 
 d_H(u,R):= \bigg(\int_{\m}\frac{\big|u(x,t)-R^{\frac{Q-p-\beta}{p}}u\big(R\frac{x}{\varrho(x,t)},R^2\frac{t}{\varrho^2(x,t)}\big)\varrho^{-\frac{Q-p-\beta}{p}}(x,t)\big|^p}{|\ln \frac{R}{\varrho(x,t)}|^p \varrho^{p+\beta}(x)}|\gradg \varrho(x,t)|^p\dv\bigg)^{\frac{1}{p}}.
\end{align*}
Now, we are ready to present our next main result of the paper. This reads as follows.
\begin{theorem}\label{stab-sub-hardy-th}
Let $G$ be a Grushin space with dimension $Q$ with $Q\geq 3$. Let $2\leq p<Q$ and $\beta<Q-p$. Then for $u\in C_0^\infty(\rno\setminus\{o\})$ there holds
\begin{multline*}
 \int_{\rno}\frac{|\rgradg u(x,t)|^p}{\varrho^{\beta}(x,t)}\dx\dt- \left(\frac{Q-p-\beta}{p}\right)^p\int_{\rno}\frac{|u(x,t)|^p|\gradg \varrho(x,t)|^p}{\varrho^{p+\beta}(x,t)}\dx\dt\\\geq c_p\bigg(\frac{p-1}{p}\bigg)^p \sup_{R>0} d_H(u,R)^p,
\end{multline*}
where $c_p$ is the positive constant defined in Lemma \ref{cpl}.
\end{theorem}

It follows immediately from the pointwise comparison $|\gradg u| \geq |\rgradg u|$ that an improvement of \cite[Theorem~3.1]{AL1} can be obtained through stability considerations, which we state in the following remark.
\begin{remark}
 Let $G$ be a Grushin space with dimension $Q$ with $Q\geq 3$. Let $2\leq p<Q$ and $\beta<Q-p$. Then for $u\in C_0^\infty(\rno\setminus\{o\})$ there holds
\begin{align*}
\int_{\rno}\frac{|\gradg u(x,t)|^p}{\varrho^{\beta}(x,t)}{\rm d}x{\rm d}t- \left(\frac{Q-p-\beta}{p}\right)^p\int_{\rno}\frac{|u(x,t)|^p|\gradg \varrho(x,t)|^p}{\varrho^{p+\beta}(x,t)}{\rm d}x{\rm d}t\geq c_p\bigg(\frac{p-1}{p}\bigg)^p \sup_{R>0} d_H(u,R)^p,
\end{align*}
where $c_p$ is the positive constant defined in Lemma~\ref{cpl}.   
\end{remark}

The second-order analogue of Hardy’s inequality is the Rellich inequality, which, in the Euclidean setting, dates back to \cite{rel}. This inequality plays a fundamental role in analytical frameworks requiring control over second derivatives, particularly in establishing regularity and a priori estimates for solutions of elliptic partial differential equations involving singular potentials or complex boundary behavior. A refinement that bridges first-order Hardy-type control with second-order norms is the Hardy--Rellich inequality, introduced in the Euclidean context in \cite{tz}. It is natural to explore such inequalities in the Grushin setting, where the underlying geometry introduces degeneracies absent in the classical case. Over recent decades, higher-order analogues of Hardy-type inequalities associated with the Baouendi--Grushin operator have been extensively investigated, beginning with the pioneering work of Kombe \cite{IK1} and followed by contributions such as \cite{KY, KY1, darca-2, GJR, sj}. In the present work, we extend these inequalities to general values of $p$, with a particular focus on the radial component of the operator. These results, presented in Theorem~\ref{high-hr} for the case $k = 2$ (Rellich inequality) and Theorem~\ref{Hardy--Rellich} for the Hardy–Rellich inequality, serve as key intermediate steps toward the main theorems of this paper.

It is now a well-established technique that by applying induction on the Hardy, Hardy--Rellich, and Rellich inequalities, one can derive inequality results for higher-order derivatives. Several works in this direction, including \cite{Barbatis, Owen, hinz, vhn, prc-20, RSY, ysk, RS2, rsy,ex-2,ex-3,ex-4,ex-5,ex-6, h24, ht25, CF24}, provide a comprehensive account of the development of such inductively defined higher-order operators in various settings. More recently, the theory of Hardy and Hardy--Rellich inequalities has witnessed significant developments through the introduction of Bessel pairs and $p$-Bessel pairs, which provide a unified framework for deriving Hardy-type identities and sharp inequalities in a variety of geometric settings. These techniques have led to several new Hardy, Hardy--Rellich, and Hardy--Sobolev inequalities on Euclidean domains, Cartan--Hadamard manifolds, hyperbolic spaces, and related geometric structures; see, for instance, \cite{gh, du, fl, be, GJR, la, fly}. In our case, we employ the inductive approach to obtain inequalities involving higher-order derivatives (see Theorem~\ref{high-hr}), and we further demonstrate that the corresponding constants are sharp. Building on this framework, we proceed to investigate the stability of these inequalities. In the Euclidean setting, this line of inquiry traces back to the work of Sano \cite{sano-mia}. Extending this approach to the Grushin operator forms the basis of our next main result. Motivated by the approximate profile $\varrho^{-(Q-kp-\beta)/p}$, we introduce a distance function that measures the deviation of a function from the family of its scaled versions. This choice is consistent with the extremizer-based stability framework developed in the Euclidean setting and naturally extends it to the Grushin setting. Thus, for a function $u\in C_{0}^\infty(\mathbb{R}^{n+1}\setminus\{o\})$, $1<p<\infty$, $\beta\in\mathbb{R}$, a fixed real number $R>0$, and an integer $k\ge2$, we define the following distance function
\begin{multline}\label{dist-high}
 d_{R}(u,k,\beta):= \\ \left( \int_{\rno}\frac{\left||u(x,t)|^{\frac{p-2}{2}}u(x,t)-|u_R(x,t)|^{\frac{p-2}{2}}u_R(x,t)\left(\frac{\varrho(x,t)}{R}\right)^{-\frac{Q-kp-\beta}{2}}\right|^2}{\varrho^{kp+\beta}(x,t)\left(\ln\left(\frac{R}{\varrho(x,t)}\right)\right)^2}|\gradg \varrho(x,t)|^p\dv\right)^{\frac{1}{2}},
\end{multline}
where $u_R(x,t)=u\big(R\frac{x}{\varrho(x,t)},R^2\frac{t}{\varrho^2(x,t)}\big)$. Now we are considering the supremum over all possible $R>0$ in the above-defined distance function in \eqref{dist-high}, so that it measures the maximum deviation of the test functions from the profile $\varrho^{-\frac{Q-kp-\beta}{2}}$ while considering all possible scalings, and our main result can be stated as follows.
\begin{theorem}\label{stab-th-higher-hr}
Let $G$ be a Grushin space of dimension $Q$ with $Q\geq 3$. Let $k$ be an integer satisfying $2\leq k<Q$. Let $2\leq p<\frac{Q}{k}$ and $\beta\in\mathbb{R}$. Then for any $u\in C_{0}^\infty(\rno\setminus\{o\})$, there exists a positive constant $C$, depending only on $k$, $p$, and $\beta$, such that the following holds:
\begin{itemize}
 \item[(i)] If $k=2l$, $l\geq 1$ and $-(2l-2)p-Q(p-1)<\beta<Q-2lp$, then we have 
 \begin{align*}
  \int_{\rno}\frac{\left|\rlapg^l u(x,t)\right|^p}{|\gradg \varrho(x,t)|^{(2l-1)p}\varrho^{\beta}(x,t)}\dx\dt-C_{k,p,\beta}\int_{\rno}\frac{|u(x,t)|^p|\gradg \varrho(x,t)|^p}{\varrho^{2lp+\beta}(x,t)}\dx\dt\geq C\sup_{R>0}d_{R}(u,k,\beta)^2;
 \end{align*}        
 \item[(ii)] If $k=2l+1$, $l\geq 1$ and $-(2l-1)p-Q(p-1)<\beta<Q-(2l+1)p$, then we have
 \begin{align*}
  \int_{\rno}\frac{\left|\rgradg \rlapg^l u(x,t)\right|^p}{|\gradg \varrho(x,t)|^{2lp}\varrho^{\beta}(x,t)}\dx\dt-C_{k,p,\beta}\int_{\rno}\frac{|u(x,t)|^p|\gradg \varrho(x,t)|^p}{\varrho^{(2l+1)p+\beta}(x,t)}\dx\dt\geq C\sup_{R>0}d_{R}(u,k,\beta)^2;
 \end{align*}
 where the constants $C_{k,p,\beta}$ are defined in Section~\ref{prelim}; in each case, they are best possible.
\end{itemize}
\end{theorem}

\begin{remark}
It is important to note that the above result for $k=2$ immediately gives the stability result for the weighted Rellich inequality on the Grushin setting. This reads as follows: Let $G$ be a Grushin space of dimension $Q$ with $Q\geq 3$. Assume $2\leq p<\frac{Q}{2}$ and $\beta\in\mathbb{R}$ with $-Q(p-1)<\beta<Q-2p$. Then for any $u\in C_{0}^\infty(\rno\setminus\{o\})$, there exists a positive constant $C$, depending only on $p$, and $\beta$, such that there holds
\begin{multline*}
 \int_{\rno}\frac{\left|\rlapg u(x,t)\right|^p}{|\gradg \varrho(x,t)|^{p}\varrho^{\beta}(x,t)}\dx\dt-\left(\frac{(Q(p-1)+\beta)(Q-\beta-2p)}{p^2}\right)^p\int_{\rno}\frac{|u(x,t)|^p|\gradg \varrho(x,t)|^p}{\varrho^{2p+\beta}(x,t)}\dx\dt\\ \geq C\sup_{R>0}d_{R}(u,2,\beta)^2.   
\end{multline*}
\end{remark}

Also, in a continuation, we want to mention that we can also achieve this type of extremizer stability estimate for critical Hardy and critical Rellich inequalities, and for this, refer to Theorem~\ref{stab-th-crt-hardy} and Theorem~\ref{stab-th-crt-rellich}. 

In the final part of this article, we present the exact remainder terms associated with the higher-order Hardy--Rellich inequalities established in Theorem~\ref{high-hr}. For the cases $k = 1$ (Hardy) and $k = 2$ (Rellich), these remainder terms can be obtained from recent results in \cite[Corollary~3.2]{GJR} and \cite[Corollary~1.2 (b)]{GJR}, where the analysis is carried out explicitly in the Grushin setting. A crucial ingredient in their approach is the delicate use of spherical harmonics developed by Garofalo and Shen \cite{GS} in the Hilbert space $L^2(G)$. More recently, Huang and Ye \cite{he} employed the Hilbert space structure of $L^2(\mathbb{R}^n)$ in the Euclidean setting to compute explicit remainder terms for higher-order Hardy--Rellich inequalities in terms of the radial components of the gradient and Laplacian. Drawing inspiration from their methodology, we adopt a similar approach to compute the explicit remainder terms for our higher-order Hardy--Rellich inequalities. To state our final main result, we first introduce the following operators. For any $u\in C_0^\infty(\rno\setminus\{o\})$ and $\beta\in\mathbb{R}$, we define a weighted first-order differential operator by
\begin{align}\label{tbeta}
 T_{\beta} (u(x,t)):= |\gradg \varrho(x,t)| \left(\frac{\partial u}{\partial \varrho}(x,t)+\frac{(Q-\beta-2)}{2\varrho(x,t)}u(x,t) \right),
\end{align}
and for $k\in \mathbb{N}\cup \{0\}$, we define recursively
\begin{align}\label{remainder-R}
    \mathcal{R}_{\beta, k}(u(x,t)):= T_{\beta}\circ T_{\beta+2}\circ \cdots \circ T_{\beta+2k}(u(x,t)).
\end{align}
We now state the final main result of the paper.
\begin{theorem}\label{main-thm-4}
Let $G$ be a Grushin space with dimension $Q$ with $Q\geq 3$. Assume that $k\geq 1$ is an integer and $\beta \in \mathbb{R}$. Then for any $u\in C_0^\infty(\rno\setminus\{o\})$ following hold: 
\begin{itemize}
 \item[(i)] If $k=2l$ and $l\geq 1$, then we have
 \begin{align}\label{radial-higher-Rellich}
  \int_{\rno}\frac{\left|\rlapg^l u(x,t)\right|^2}{|\gradg \varrho(x,t)|^{2(2l-1)}\varrho^{\beta}(x,t)}\dx\dt=& C_{k,2,\beta}\int_{\rno}\frac{|u(x,t)|^2|\gradg \varrho(x,t)|^2}{\varrho^{\beta+4l}(x,t)}\dx\dt\\
  \nonumber +& \sum_{j=0}^{2l-1} D_{j,l, \beta} \int_{\rno} \frac{|\mathcal{R}_{\beta+2j, 2l-1-j} (u(x,t))|^2}{|\gradg \varrho(x,t)|^{2(2l-1-j)}\varrho^{\beta+2j}(x,t)}  \dx\dt
\end{align}
where the constants $D_{j,l,\beta}$ are given by the iterative formula
\begin{align}\label{coeff-iteration}
 D_{j,m+1,\beta} =  \Lambda_{2m}(Q,2,\beta) D_{j-2,m,\beta} + A_{\beta+4m} D_{j-1,m,\beta} + D_{j,m,\beta},
\end{align}
for all $0 \leq j \leq 2m+2, ~m\geq 0$, with the convention that $D_{0,0,\beta}=1$ for every $\beta$ and $D_{j,m,\beta}=0$ for $j<0$ or $j>2m$ and $\Lambda_{2m}(Q,2,\beta)$, $A_{\beta+4m}$ are defined in \eqref{big-lam}, \eqref{Dm-constant} respectively.
 \item[(ii)] If $k=2l+1$ and $l\geq 0$, then we have
  \begin{align}\label{Higher-order-odd}
  \int_{\rno}\frac{\left|\rgradg \rlapg^l u(x,t)\right|^2}{|\gradg \varrho(x,t)|^{4l}\varrho^{\beta}(x,t)}\dx\dt=& C_{k,2,\beta}\int_{\rno}\frac{|u(x,t)|^2|\gradg \varrho(x,t)|^2}{\varrho^{2(2l+1)+\beta}(x,t)}\dx\dt\\
  \nonumber +& \sum_{j=0}^{2l} \widetilde{D}_{j,l, \beta} \int_{\rno} \frac{|\mathcal{R}_{\beta+2j, 2l-j} (u(x,t))|^2}{|\gradg \varrho(x,t)|^{2(2l-j)}\varrho^{\beta+2j}(x,t)}  \dx\dt
 \end{align}
  where $\widetilde{D}_{j,l, \beta}=B_{j,l,\beta, 0}+\frac{(Q-\beta-2)^2}{4} D_{j-1, l,\beta}$ for all $0\leq j \leq 2l$, where $B_{j,l,\beta, 0}$ and $D_{j, l,\beta}$ are the positive constants defined in \eqref{coeff-R-alpha-ite} and \eqref{coeff-iteration} , respectively.
  \end{itemize}
\end{theorem}

\begin{remark}
Let $k$, $\beta$, $j$, $l$, and $m$ be as in Theorem~\ref{main-thm-4}. We observe that, for all admissible parameters, the constants $C_{k,2,\beta}$ (defined in \eqref{sharp-constant}), $D_{j,l,\beta}$ (defined in \eqref{radial-higher-Rellich}), and $\widetilde{D}_{j,l,\beta}$ (defined in \eqref{Higher-order-odd}) are all nonnegative. Indeed, when $p=2$, the quantities $\Lambda_{\bullet}(\bullet,2,\bullet)$ defined in \eqref{big-lam} and \eqref{big-lam-1} are nonnegative. Moreover, the coefficients $A_{\beta+4m}$ defined in \eqref{Dm-constant} are also nonnegative. Since $D_{j,l,\beta}$ is defined recursively in terms of these nonnegative coefficients, together with nonnegative lower-order terms, an induction on the recursion shows that $D_{j,l,\beta}\ge0$ for all admissible parameters $j$, $l$, and $\beta$. Likewise, the recursively defined coefficients $B_{j,l,\beta,0}$ in \eqref{coeff-R-alpha-ite} are nonnegative, which immediately implies that $\widetilde{D}_{j,l,\beta}\ge 0$.
\end{remark}

The paper is organized as follows. In Section~\ref{prelim}, we recall some necessary preliminaries, including polar coordinates, the description of radial operators, and spherical harmonics. We also establish several essential lemmas and provide a complete description of the relevant constants. In Section~\ref{sh-ch-est}, we derive a two-weighted Hardy inequality (Theorem~\ref{two-wg-hardy}) and, as consequences, we present a critical Hardy inequality and a geometric Hardy inequality. The higher-order analogues of the Hardy inequality, namely the Rellich and Hardy--Rellich inequalities, are established in Section~\ref{r-hr-est}. In Section~\ref{sh-ch-rel}, we prove Theorem~\ref{main-thm-1}, which illustrates the interplay between subcritical and critical Hardy inequalities in the Grushin setting. This section concludes with improvements to the subcritical and critical Hardy inequalities, incorporating positive remainder terms that arise from the Caffarelli–Kohn–Nirenberg inequality. Section~\ref{ex-st-an} is devoted to the proofs of Theorem~\ref{stab-sub-hardy-th} and Theorem~\ref{stab-th-higher-hr}. Finally, in Section~\ref{hr-r-id}, we establish Theorem~\ref{main-thm-4}, which provides an exact expression for the remainder terms in the higher-order Hardy--Rellich inequality. This is achieved by deriving an identity rather than merely proving an inequality.

\section{Preliminary concepts}\label{prelim}
We will collect some notations and background information in this section, which will be utilized throughout the article.

{\bf Polar coordinate.} Let's start by recalling the polar coordinate structure, which was first introduced for Grushin settings in \cite{GS}. Let $G=\rno$ be the Grushin space. For any  $(x,t)\in \rno$, we can write $(x,t)=(\varrho,\phi,\theta_1,\theta_2,\cdots,\theta_{n-1})$ as below:	
\begin{align*}
		&x_1=\varrho \,\sin^{1/2}\phi\, \sin\theta_1\cdots\sin \theta_{n-1},\\
		&x_2=\varrho \,\sin^{1/2}\phi\, \sin\theta_1\cdots\cos \theta_{n-1},\\
		&\cdots\\
		&x_n=\varrho \,\sin^{1/2}\phi\, \cos\theta_1,\\
		&t=\frac{\varrho^2}{2}\cos \phi,
\end{align*}
	where
	$0<\phi<\pi$, $0<\theta_i<\pi$ for $i=1,2,\cdots,n-2$ and $0<\theta_{n-1}<2\pi$. Here, we exclude the endpoints, as they form sets of measure zero. Now, denoting $r = |x| = \sqrt{x_1^2 + \cdots + x_n^2}$, we can observe that $r = \varrho \sin^{1/2} \phi$, and define the \emph{geometric term} $\psi := \sin \phi = \frac{|x|^2}{\varrho^2}$. Hence, this will give 
	\begin{align*}
		|\gradg \varrho(x,t)|^2=|\nabla_x\varrho|^2+|x|^2(\partial_t \varrho)^2=\frac{|x|^2}{\varrho^2(x,t)}=\psi(x,t)=\sin\phi = \frac{r^2}{\varrho^2},
	\end{align*}	
	and 
	\begin{align*}
		\lapg\varrho (x,t) = (Q-1)\frac{|x|^2}{\varrho^3(x,t)}=(Q-1)\frac{\psi(x,t)}{\varrho(x,t)}.
	\end{align*}

{\bf Volume elements and polar integration.} Using polar coordinates on Grushin space $G=\rno$, we obtain $\dx=r^{n-1}{\rm d}r\,{\rm d}w$, and 
\begin{align}\label{volume}
        \nonumber\dx\dt&=\frac{1}{2}\varrho^{n+1}(\sin\phi)^{\frac{n-2}{2}}\dr\dph\dw\\&=\frac{1}{2}\varrho^{n+1}(\sin\phi)^{\frac{n-2}{2}}(\sin \theta_1)^{n-2}(\sin \theta_2)^{n-3}\ldots(\sin \theta_{n-2})\dr\dph\,{\rm d}\theta_1\,{\rm d}\theta_2\ldots\,{\rm d}\theta_{n-2}\,{\rm d}\theta_{n-1}
\end{align}
	where $\dw$ is the spherical measure on the standard round sphere $\mathbb{S}^{n-1}$ in the Euclidean space $\mathbb{R}^n$. Now consider the unit $\varrho$-sphere on the Grushin space as follows:
	\begin{align*}
		\Omega=\biggl\{(x,t)\in \rno : \varrho(x,t)=1\biggr\}.
	\end{align*}
	Moreover, using the standard argument, there is the following polar integration. Let $u$ be some integrable function, and then we have
	\begin{align*}
		\int_{\rno}u(x,t) \dx \dt=\int_{\Omega}\int_{0}^{\infty}u(\varrho,\phi,w)\frac{1}{2}\varrho^{n+1}(\sin\phi)^{\frac{n-2}{2}}\dr\dph\dw=\int_{\Omega}\int_{0}^{\infty}u(\varrho,\sigma)\frac{\varrho^{n+1}}{2\psi}\dr\dsn,
	\end{align*}
	where $${\rm d}\Omega=(\sin\phi)^{\frac{n}{2}}\dph\dw, \text{ and } \sigma=(\phi,w), \text{ with }w\in \mathbb{S}^{n-1}.$$
	
{\bf Operators and their radial part.} Now, we will write all the operators in terms of the polar coordinate structure on $G=\rno$. The Grushin gradient in polar coordinates may be defined as follows,
\begin{align*}
\gradg=\psi^{1/2}\bigg(\partial_\varrho, \frac{2}{\varrho}\partial_\phi,\frac{1}{\varrho\sin\phi}\nabla_w\bigg),
\end{align*}
	where $(x,t)=(\varrho,\phi,w)$, $w\in \mathbb{S}^{n-1}$, $\nabla_w$ is the gradient on $\mathbb{S}^{n-1}$. Also, it can be written as follows:
\begin{align*}
\gradg=\psi^{1/2}\bigg(\partial_\varrho, \frac{2}{\varrho}\nabla_\sigma\bigg),
	\end{align*}
	where $\nabla_\sigma$ is the gradient operator on the Grushin $\varrho$-sphere. We refer to \cite[Page. 11]{flynn} for more details.	
	In polar coordinates, the Grushin operator takes the form 
	\begin{align*}
		\lapg=\psi\biggl\{\frac{\partial^2}{\partial \varrho^2}+\frac{Q-1}{\varrho}\frac{\partial}{\partial \varrho}+\frac{4}{\varrho^2}\mathcal{L}_\sigma\biggr\},
	\end{align*}
where $\sigma=(\phi,w)$ and $\mathcal{L}_\sigma$ is the Laplacian on Grushin $\varrho$-sphere.
	
To this end, we define the radial contribution of the operators. The definitions of the \emph{radial Grushin gradient} and the \emph{radial Grushin operator} are
	\begin{align}\label{radial-op}
		\nabla_{\varrho,G}=\psi^{1/2}\big(\partial_\varrho,0,0\big), \quad \text{ and } \quad \mathcal{L}_{\varrho,G}=\psi\biggl\{\frac{\partial^2}{\partial \varrho^2}+\frac{Q-1}{\varrho}\frac{\partial}{\partial \varrho}\biggr\}.
	\end{align}

{\bf Spherical harmonics.} Assume that the Grushin space is $G=\rno$ and that its unit $\varrho$-sphere is $\Omega$. For $k=0,1,\cdots,$ we form the function $u_k(\varrho, \sigma):=\varrho^k g(\sigma)$. Then, $u_k(\varrho, \sigma)$ is a solution of $\mathcal{L}_{G}u_k = 0$ if and only if 
	\begin{align*}
		\mathcal{L}_\sigma g=-\frac{k(n+k)}{4}g.
	\end{align*}
    From \cite[Lemma~2.1]{GS}, we have $u(x, t)=u(\varrho,\sigma)\in C_{0}^\infty(\rno\setminus\{o\})$, with $n\geq2$, $\varrho\in({0},\infty)$ and $\sigma\in \Omega$, thus we can write
	\begin{equation*}
		u(\varrho,\sigma)=\sum_{k=0}^{\infty}d_{k}(\varrho)\Phi_k(\sigma)
	\end{equation*}
	in $L^2(\rno)$, where $\{ \Phi_k \}$ is an orthonormal system of spherical harmonics in $L^2(\Omega)$ and 
	\begin{equation*}
		d_{k}(\varrho)=\int_{\Omega}u(\varrho,\sigma)\Phi_k(\sigma) \ {\rm d}\Omega\,.
	\end{equation*} 
	A spherical harmonic $\{ \Phi_k \}$ of order $k$ satisfies $$-\mathcal{L}_{\sigma}\Phi_k=\lambda_k\Phi_k,$$
	for all $k\in\mathbb{N}\cup\{0\}$, where $\lambda_k=\frac{k(k+n)}{4}$. 

We now state a basic lemma that provides a bound for the integral of the geometric function $\psi$ in terms of the $\varrho$-sphere component. 
\begin{lemma}\label{i-alpha}
  Let $\alpha\geq 0$ be any real number and the integral $\mathcal{I}_\alpha$ is defined by
  \begin{align*}
     \mathcal{I}_\alpha:=\frac{1}{2}\int_{\mathbb{S}^{n-1}}\int_{0}^\pi(\sin\phi)^\alpha\,{\rm d}\phi\,{\rm d}w. 
  \end{align*}
  Then there holds
  \begin{align*}
     \frac{\pi \, \omega_n}{3\cdot2^{\alpha}}\leq  \mathcal{I}_\alpha \leq \frac{\pi\, \omega_n}{2},
  \end{align*}
  where $\omega_n$ is the spherical measure of the round Euclidean unit sphere $\mathbb{S}^{n-1}$ in $\re^n$.
\end{lemma}
\begin{proof}
    Using the pointwise estimate $\frac{1}{2}\leq \sin\phi \leq 1$ for any $\phi\in [\frac{\pi}{6},\frac{5\pi}{6}]$, the estimate follows.
\end{proof}

We now recall an important lemma (see \cite[Formula (2.13)]{FS08}) that will be useful for our further analysis.
\begin{lemma}[\cite{FS08}]\label{cpl}
	Let $p \geq 2$. Then there exists $c_p > 0$ such that
	\begin{align*}
		|a-b|^p \geq |a|^p - p|a|^{p-2}{\rm Re\;} \overline{a}\cdot b + c_p|b|^p
	\end{align*}
	holds for all vectors $a, b \in \mathbb{C}^N$, where 
	\begin{align*}
		c_p=\min_{0<t< 1/2}[(1-t)^p-t^p+pt^{p-1}]
	\end{align*}
	is sharp in this inequality.
\end{lemma}
 
{\bf Description of $\Lambda_{k}(Q,p,\beta)$ and $C_{k,p,\beta}$.} We define for $\beta$, $Q$, and $p$ the following notations
    \begin{align}\label{big-lam-1}
        \Lambda_{1}(Q,p,\beta)=\left(\frac{Q-p-\beta}{p}\right)^p,
    \end{align}
    and
\begin{align}\label{big-lam}
    \Lambda_{k}(Q,p,\beta)= \bigg(\frac{(Q(p-1)+\beta+(k-2)p)(Q-\beta-kp)}{p^2}\bigg)^p \text{ for integer } k\geq 2.
\end{align}
Also, define the following constants 
\begin{align*}
    C_{0,p,\beta}:=1,
\end{align*}
and for $k\geq 1$ integer
\begin{align}\label{sharp-constant}
		C_{k,p,\beta} :=\left\{
\begin{array}{ll}
			\prod_{i=1}^{l} \Lambda_{2i}(Q,p,\beta) , &\quad k=2l,\\
			\prod_{i=0}^{l} \Lambda_{2i+1}(Q,p,\beta) , &\quad k=2l+1.\\
		\end{array} 
		\right.
\end{align}

We now state two lemmas that will serve as key ingredients in establishing inequalities for the higher-order Grushin operator.
\begin{lemma}\label{ind-lemm}
Let $G = \rno$ denote the Grushin space of homogeneous dimension $Q = n+2$. Let $k\in\mathbb{N}$, and let $p$ and $\beta$ be real numbers. Then there holds
\begin{align*}
		\left\{
		\begin{array}{ll}
			\psi^{-l}\left|\rlapg^l \varrho^{-\frac{Q-\beta-kp}{p}}\right|=C_{k,p,\beta}^{\frac{1}{p}}\varrho^{-\frac{Q-\beta}{p}}, & \text{ for } \quad k=2l,\\
			\psi^{-\frac{2l+1}{2}}\left|\rgradg \rlapg^l \varrho^{-\frac{Q-\beta-kp}{p}}\right|=C_{k,p,\beta}^{\frac{1}{p}}\varrho^{-\frac{Q-\beta}{p}}, & \text{ for } \quad k=2l+1.\\
		\end{array} 
		\right.
	\end{align*}
\end{lemma}
\begin{proof}
We will proceed by induction. Before that, notice that for any real number $\nu$, there holds
\begin{align}\label{nu-ind}
    \left(\pd^2  + (Q-1)\varrho^{-1}\pd\right)\varrho^{\nu}=\nu(Q+\nu-2)\varrho^{\nu-2}. 
\end{align}
Let us first consider the case $k=1$, which gives $l=0$. Then we have
\begin{align*}
    \psi^{-\frac{1}{2}}|\rgradg \varrho^{-\frac{Q-\beta-p}{p}}|=\left(\frac{Q-\beta-p}{p}\right)\varrho^{-\frac{Q-\beta}{p}}=C_{1,p,\beta}^{\frac{1}{p}}\varrho^{-\frac{Q-\beta}{p}}.
\end{align*}
Next, consider the case $k=2$ and so $l=1$. This gives using \eqref{nu-ind} with $\nu=-\frac{Q-\beta-2p}{p}$, that
\begin{align*}
    \psi^{-1}|\rlapg \varrho^{-\frac{Q-\beta-2p}{p}}|=\left(\frac{Q-\beta-2p}{p}\right)\left(\frac{Q(p-1)+\beta}{p}\right)\varrho^{-\frac{Q-\beta}{p}}=C_{2,p,\beta}^{\frac{1}{p}}\varrho^{-\frac{Q-\beta}{p}}.
\end{align*}
Next, consider the case $k=3$ and so $l=1$. Then we have
\begin{align*}
    \psi^{-\frac{3}{2}}|\rgradg \rlapg \varrho^{-\frac{Q-\beta-3p}{p}}|=\Lambda_3^{
    \frac{1}{p}}\left|\pd \varrho^{-\frac{Q-\beta-p}{p}}\right|=C_{3,p,\beta}^{\frac{1}{p}}\varrho^{-\frac{Q-\beta}{p}}.
\end{align*}
The rest of the part will follow from induction and repeated use of \eqref{nu-ind}.
\end{proof}
\begin{lemma}\label{lib-ind}
   Let $G$ be a Grushin space of homogeneous dimension $Q$. Let $f\in C^\infty(\rno)$ be a radial function, and let $\nu$ be any real number. Suppose $k\geq 1$ is an integer. Then, for each $k$, there exist $(k+1)$ constants depending only on $Q$ and $\nu$ such that, for any fixed $\tau>0$, the identity
   \begin{align*}
		\left\{
		\begin{array}{ll}
			\psi^{-l}\rlapg^l \left(f(\varrho/\tau)\varrho^{\nu}\right)= \sum_{i=0}^{2l} \tau^{-i}f^{(i)}(\varrho/\tau)\varrho^{\nu-2l+i}h_{k,i}, & \text{for} \quad k=2l,\\
			\psi^{-l}\pd \rlapg^l\left(f(\varrho/\tau)\varrho^{\nu}\right)=\sum_{i=0}^{2l+1} \tau^{-i}f^{(i)}(\varrho/\tau)\varrho^{\nu-2l-1+i}h_{k,i}, & \text{for}  \quad k=2l+1,\\
		\end{array} 
		\right.
	\end{align*}
    holds, where ${h_{k,i}}$ are independent of $\tau$ for each $k$ and $i$, and $f^{(i)}$ denotes $\partial^i f$ for each $i$.
\end{lemma}
\begin{proof}
Again, the result follows by induction and the Leibniz rule. For $k=1$, we have $l=0$, and thus \begin{align*} \partial\left(f(\varrho/\tau)\varrho^{\nu}\right) = \nu f(\varrho/\tau)\varrho^{\nu-1} + \tau^{-1}f^{(1)}(\varrho/\tau)\varrho^{\nu}, \end{align*} and the case is completed by choosing $h_{1,0} = \nu$ and $h_{1,1} = 1$.
Next, consider the case $k=2$, so that $l=1$. Then we compute
\begin{align*}
\psi^{-1}\rlapg\left(f(\varrho/\tau)\varrho^{\nu}\right) &= \nu(\nu-1) f(\varrho/\tau)\varrho^{\nu-2} + \nu\tau^{-1} f^{(1)}(\varrho/\tau)\varrho^{\nu-1} \\ &\quad + \tau^{-2} f^{(2)}(\varrho/\tau)\varrho^{\nu} + \tau^{-1}\nu f^{(1)}(\varrho/\tau)\varrho^{\nu-1} \\ &\quad + (Q-1)\nu f(\varrho/\tau)\varrho^{\nu-2} + (Q-1)\tau^{-1}f^{(1)}(\varrho/\tau)\varrho^{\nu-1}.
\end{align*}
Setting $h_{2,0} = \nu^2 + (Q-2)\nu$, $h_{2,1} = 2\nu + (Q-1)$, and $h_{2,2} = 1$, the case $k=2$ is completed. The general case follows by induction.
\end{proof}

{\bf Basic facts and notations.} In this article, we will use the following notation and terminology:
\begin{itemize}
    \item {\bf Grushin Spaces and Dimensions}
    \begin{itemize}
        \item Primary space: $G=\rno(n\geq 1)$ with homogeneous dimension $Q=n+2$.
        \item Auxiliary spaces: $G'=\rmo$, $G''=\re^{q+1}$ with $m\geq 1$, $n\geq 1$ and $Q'=m+2$, $Q''=q+2$.
        \item Homogeneous norms: $\varrho$ on $G$, $\tilde{\varrho}$ on $G'$, and $\bar{\varrho}$ on $G''$.
        \item Space variables: $(x,t)$ on $G$, $(y,s)$ on $G'$, and $(z,r)$ on $G''$.
        \item Origins: $o$ on $G$, $o'$ on $G'$, and $o''$ on $G''$.
        \item Grushin gradient: $\nabla_{G}$ on $G$, $\nabla_{G'}$ on $G'$, and $\nabla_{G''}$ on $G''$. 
        \item Radial Grushin gradient:  $\rgradg$ on G, $\nabla_{\tilde{\varrho},G'}$ on $G'$, and $\nabla_{\bar{\varrho},G''}$ on $G'$. 
    \end{itemize}
    \item {\bf Function Spaces}
    
    For any open $\mathcal{J}\subset G$
    \begin{itemize}
        \item $C_0^\infty(\mathcal{J})$: Compactly supported smooth real valued functions on $\mathcal{J}$.
        \item $C_{0,rad}^\infty(\mathcal{J})$: Compactly supported smooth radial real valued functions on $\mathcal{J}$.
        \item $L^p(\mathcal{J})$: $p$-integrable real valued functions $(1\leq p<\infty)$.
        \item $W^{k,p}_0(\mathcal{J})$: $k$-th order Sobolev space with vanishing trace.
    \end{itemize}
    \item {\bf Geometric Structures}
    \begin{itemize}
        \item Gauge unit spheres: $\Omega$ (for $G$), $\Omega'$ (for $G'$), and $\Omega''$ (for $G''$).
        \item Geometric terms: $\psi$ (for $G$), $\tilde{\psi}$ (for $G'$), and $\bar{\psi}$ (for $G''$).
        \item Ball of radius $R>0$: $\br$ (for $G$), $\brp$ (for $G'$), and $B_R^{Q''}(o'')$ (for $G''$).
        \end{itemize}
        \item {\bf Integral Constants}

        For each Grushin space:
        \begin{itemize}
            \item On $G$: For $k_1 \geq -n$, $C^{\Omega,k_1}
    := \int_{\Omega} \psi^{\frac{k_1}{2}} \dsn$.
    \item On $G'$: For $k_2 \geq -m$, $C^{\Omega,k_2}
    := \int_{\Omega} \psi^{\frac{k_2}{2}} \dsnp$.
    \item On $G''$: For $k_3 \geq -q$, $C^{\Omega,k_3}
    := \int_{\Omega} \psi^{\frac{k_3}{2}} \dsn''$.
        \end{itemize}
        \item {\bf Key Properties}
        \begin{itemize}
            \item Density: $C_0^\infty(\br)$ is dense in $W^{k,p}_0(\br)$, as the boundary $\partial \br$ is smooth and non-tangential to the degenerate set $\{0\} \times \mathbb{R}$ (see \cite[Section~5]{Laura}).
            \item Boundedness: By applying Lemma~\ref{i-alpha} and using the volume description of the gauge unit sphere, all integral constants are bounded from above and below.
        \end{itemize}        
\end{itemize}

\section{Hardy-type inequalities}\label{sh-ch-est}
This section is devoted to the establishment of some Hardy-type inequalities. First, we consider a Hardy inequality for the difference between a function and its scaled version. This formulation serves as the core for establishing stability-type results. In the Euclidean setting, such studies can be traced back to~\cite{mow, mow-2}. Here, we present an analogous result in the Grushin setting.
\begin{theorem}\label{diff-crit-hardy}
     Let $G$ be a Grushin space with dimension $Q$ with $Q\geq 3$. Assume $1<\gamma<\infty$ and $\max\{1,\gamma-1\}<p<\infty$. Then, for all $u\in C_0^\infty(\rno)$ and all, $R>0$ we have the inequality
     \begin{multline}\label{diff-crit-hardy-eq}
         \bigg|\bigg|\varrho^{\frac{p-Q}{p}}(x,t)\bigg(\ln \frac{R}{\varrho(x,t)}\bigg)^{\frac{p-\gamma}{p}}|\rgradg u(x,t)|\bigg|\bigg|_{L^p(\rno)}\\ \geq \frac{(\gamma-1)}{p}\bigg|\bigg|\frac{u(x,t)-u_R(x,t)}{\varrho^{\frac{Q}{p}}(x,t)\big(\ln\frac{R}{\varrho(x,t)}\big)^{\frac{\gamma}{p}}}|\gradg \varrho(x,t)|\bigg|\bigg|_{L^p(\rno)},
     \end{multline}
    where $u_R(x,t):=u\big(R\frac{x}{\varrho(x,t)},R^2\frac{t}{\varrho^2(x,t)}\big)$ and the constant $\frac{(\gamma-1)}{p}$ is sharp.
\end{theorem}
\begin{proof}
     Let us start the computation on the $\varrho$-gauge ball $\br$ with radius $R>0$. Using the polar coordinate, we deduce the following:
     \begin{align*}
        \int_{\br}&\frac{|u-u_R|^p}{|\varrho^{\frac{Q}{p}}\big(\ln\frac{R}{\varrho}\big)^{\frac{\gamma}{p}}|^p}|\gradg \varrho|^p\dx\dt=\int_{0}^R\int_{\sn}\frac{|u(\varrho,\sigma)-u(R,\sigma)|^p}{\big(\ln\frac{R}{\varrho}\big)^{\gamma}}\frac{\psi^{\frac{p}{2}-1}\varrho^{-1}}{2}\dsn\dr\\&=\int_{0}^R \pd\bigg(\int_{\sn}\frac{|u(\varrho,\sigma)-u(R,\sigma)|^p}{(\gamma-1)\big(\ln\frac{R}{\varrho}\big)^{\gamma-1}}\frac{\psi^{\frac{p}{2}-1}}{2}\dsn\bigg)\dr\\&-\frac{p}{\gamma-1}\int_{0}^R \int_{\sn}\frac{|u(\varrho,\sigma)-u(R,\sigma)|^{p-2}(u(\varrho,\sigma)-u(R,\sigma))}{\big(\ln\frac{R}{\varrho}\big)^{\gamma-1}}\pd(u(\varrho,\sigma))\frac{\psi^{\frac{p}{2}-1}}{2}\dsn\dr.
    \end{align*}
    The first integral vanishes due to $\gamma-1>0$ and the $u(\varrho,\sigma)-u(R,\sigma)=0$ when $\varrho=R$. So, writing back in terms of the original coordinate and using H\"older inequality, we deduce
    \begin{align*}
        \int_{\br}\frac{|u-u_R|^p}{|\varrho^{\frac{Q}{p}}\big(\ln\frac{R}{\varrho}\big)^{\frac{\gamma}{p}}|^p}|\gradg \varrho|^p\dx\dt& \leq\frac{p}{(\gamma-1)}\int_{\br}\frac{|u-u_R|^{p-1}|\gradg\varrho|^{p-1}}{\varrho^{Q-1}\big(\ln\frac{R}{\varrho}\big)^{\gamma-1}}|\rgradg u|\dx\dt\\&\leq \frac{p}{(\gamma-1)}\bigg(\int_{\br}\frac{|u-u_R|^{p}}{\big|\varrho^{\frac{Q}{p}}\big(\ln\frac{R}{\varrho}\big)^{\frac{\gamma}{p}}\big|^p}|\gradg \varrho |^p\dx\dt\bigg)^{(p-1)/p}\\& \quad \quad \times \bigg(\int_{\br}\varrho^{(p-Q)}\bigg(\ln\frac{R}{\varrho}\bigg)^{(p-\gamma)}|\rgradg u|^p\dx\dt\bigg)^{1/p}.
    \end{align*}
    Thus, we obtain 
    \begin{align*}
        \bigg(\int_{\br}\varrho^{(p-Q)}\bigg|\ln\frac{R}{\varrho}\bigg|^{(p-\gamma)}|\rgradg u|^p\dx\dt\bigg)^{\frac{1}{p}}\geq   \frac{(\gamma-1)}{p} \bigg(\int_{\br}\frac{|u-u_R|^p}{|\varrho^{\frac{Q}{p}}\big(\ln\frac{R}{\varrho}\big)^{\frac{\gamma}{p}}|^p}|\gradg \varrho|^p\dx\dt\bigg)^{\frac{1}{p}}.
    \end{align*}
       Similarly, one can observe that
       \begin{align*}
           \bigg(\int_{\br^c}\varrho^{(p-Q)}\bigg|\ln\frac{R}{\varrho}\bigg|^{(p-\gamma)}|\rgradg u|^p\dx\dt\bigg)^{\frac{1}{p}} \geq \frac{(\gamma-1)}{p}\bigg(\int_{\br^c}\frac{|u-u_R|^p}{|\varrho^{\frac{Q}{p}}\big(\ln\frac{R}{\varrho}\big)^{\frac{\gamma}{p}}|^p}|\gradg \varrho|^p\dx\dt\bigg)^{\frac{1}{p}}.
       \end{align*}
    Combining both, we deduce the required inequality. 
    
Now we verify the optimality of the constant $\frac{(\gamma-1)}{p}$. First note that from \eqref{diff-crit-hardy-eq} it implies
  \begin{align*}
       \bigg(\int_{\br}\varrho^{(p-Q)}\bigg(\ln\frac{R}{\varrho}\bigg)^{(p-\gamma)}|\rgradg f|^p\dx\dt\bigg)^{\frac{1}{p}}\geq \frac{(\gamma-1)}{p} \bigg(\int_{\br}\frac{|f|^p}{\varrho^{Q}\big(\ln\frac{R}{\varrho}\big)^{\gamma}}|\gradg \varrho|^p\dx\dt\bigg)^{\frac{1}{p}}
  \end{align*}
    for all $f\in C_0^\infty(\br)$.
    Therefore, it is enough to consider optimality for the above inequality. Suppose $\delta$ is a small positive number less than $R$, and then consider the sequence of functions
\begin{equation}\label{seq-fun}
 f_\delta(x,t)=f_\delta(\varrho,\sigma):=
\begin{dcases}
\bigg(\ln\frac{R}{\delta}\bigg)^{\frac{\gamma-1}{p}} &  \text{ when } \varrho\leq \delta; \\
\bigg(\ln\frac{R}{\varrho}\bigg)^{\frac{\gamma-1}{p}} & \text{ when } \delta\leq \varrho\leq  \frac{R}{2}; \\
\left(\ln 2 \right)^{\frac{\gamma-1}{p}}\frac{2}{R}(R-\varrho) &  \text{ when }\frac{R}{2}\leq \varrho\leq  R. \\
\end{dcases}
\end{equation}
This gives that 
\begin{equation}\label{seq-fun-der}
\pd f_\delta=
\begin{dcases}
0 &  \text{ when } \varrho<\delta; \\
-\left(\frac{\gamma-1}{p}\right)\bigg(\ln\frac{R}{\varrho}\bigg)^{\frac{\gamma-1}{p}-1}\varrho^{-1} & \text{ when } \delta < \varrho <  \frac{R}{2}; \\
-\left(\ln 2 \right)^{\frac{\gamma-1}{p}}\frac{2}{R} & \text{ when }\frac{R}{2} < \varrho <  R. \\
\end{dcases}
\end{equation}
It is easy to check that $f_\delta \in W^{1,p}_0(\br)$ and, using density arguments, it will serve as a test function. Now, from polar coordinate decomposition and Lemma~\ref{i-alpha} for $\alpha=\frac{n+p}{2}-1\geq 0$, it follows that
 \begin{align*}
       \int_{\br}\frac{|f_\delta|^p}{\varrho^{Q}\big(\ln\frac{R}{\varrho}\big)^{\gamma}}|\gradg \varrho|^p\dx\dt&\geq\int_{B^Q_{\frac{R}{2}}(o)\setminus B^Q_\delta(o)}\frac{|f_\delta|^p}{\varrho^{Q}\big(\ln\frac{R}{\varrho}\big)^{\gamma}}|\gradg \varrho|^p\dx\dt \\ &= \int_{\sn}\int_\delta^{\frac{R}{2}} \bigg(\ln\frac{R}{\varrho}\bigg)^{-1}\frac{\psi^{\frac{p}{2}-1}}{2\varrho } \dr\dsn \\&= \mathcal{I}_{\frac{n+p}{2}-1}\int_\delta^\frac{R}{2} \bigg(\ln\frac{R}{\varrho}\bigg)^{-1}\frac{1}{\varrho} \dr\\&=\mathcal{I}_{\frac{n+p}{2}-1}\left(\ln\left(\ln\frac{R}{\delta}\right)-\ln(\ln2))\right).
    \end{align*}
    So, we have 
    \begin{align*}
        \lim_{\delta\rightarrow 0^+} \int_{\br}\frac{|f_\delta|^p}{\varrho^{Q}\big(\ln\frac{R}{\varrho}\big)^{\gamma}} 
        |\gradg \varrho|^p\dv = +\infty.
    \end{align*}    
    Also, we deduce
    \begin{align*}
        &\int_{\br}\varrho^{(p-Q)}\big(\ln\frac{R}{\varrho}\big)^{(p-\gamma)}|\rgradg f_\delta|^p\dx\dt\\&\leq \bigg(\frac{\gamma-1}{p}\bigg)^p\int_{B^Q_{\frac{R}{2}}(o)\setminus B^Q_\delta(o)}\frac{|f_\delta|^p}{\varrho^{Q}\big(\ln\frac{R}{\varrho}\big)^{\gamma}}|\gradg \varrho|^p\dx\dt\\&+(\ln 2)^{\gamma-1}2^pR^{-p}\int_{B^Q_R(o)\setminus B^Q_{\frac{R}{2}}(o)}\varrho^{(p-Q)}\left(\ln\frac{R}{\varrho}\right)^{(p-\gamma)}|\gradg \varrho|^p\dv.
    \end{align*}
    Again using Lemma~\ref{i-alpha} for $\alpha=\frac{n+p}{2}-1\geq 0$, it follows that
    \begin{align*}
        \int_{B^Q_R(o)\setminus B^Q_{\frac{R}{2}}(o)}\varrho^{(p-Q)}\left(\ln\frac{R}{\varrho}\right)^{(p-\gamma)}|\gradg \varrho|^p\dv&\leq R^p \mathcal{I}_{\frac{n+p}{2}-1}\int_{\frac{R}{2}}^R\left(\ln\frac{R}{\varrho}\right)^{(p-\gamma)}\frac{1}{\varrho}\dr\\&=\frac{R^p \mathcal{I}_{\frac{n+p}{2}-1}}{(p-\gamma+1)}\left(\ln2\right)^{(p-\gamma+1)}<+\infty,
    \end{align*}
    where we used the fact that $p-\gamma+1>0$. This gives, in combining that
    \begin{align*}
       \lim_{\delta\rightarrow 0^+}\frac{\int_{\br}\varrho^{(p-Q)}\big(\ln\frac{R}{\varrho}\big)^{(p-\gamma)}|\rgradg f_\delta|^p\dx\dt}{\int_{\br}\frac{|f_\delta |^p}{\varrho^{Q}\big(\ln\frac{R}{\varrho}\big)^{\gamma}}|\gradg \varrho|^p\dx\dt} =\bigg(\frac{\gamma-1}{p}\bigg)^p.
    \end{align*}
\end{proof}

We now establish the two-weighted form of the Hardy inequality in the radial Grushin setting, which is of independent interest in the literature. Although the underlying method for deriving such inequalities is classical in the Euclidean setting and can be traced back to~\cite{sano-mjm}, the present work extends it to the Grushin framework.
\begin{theorem}\label{two-wg-hardy}
   Let $G$ be a Grushin space with dimension $Q$ with $Q\geq 3$. Assume $\br$ be the $\varrho-$gauge ball with radius $R>0$. Let $1< p<\infty$, $a< Q$, $p\leq b<\infty$, $c>0$, and $a\leq Q- (b-1)c$. Then for any $u\in C_0^\infty (\br)$ there holds
   \begin{multline}\label{two-wg-hardy-eq}
        \int_{\br} \frac{|\rgradg u(x,t)|^p}{\varrho^{a-p}(x,t)\left(1-\left(\frac{\varrho(x,t)}{R}\right)^c\right)^{b-p}} \dx\dt \geq  \bigg(\frac{(b-1)}{p}c\bigg)^p\int_{\br} \frac{|u(x,t)|^{p}|\gradg \varrho(x,t)|^p}{\varrho^{a}(x,t)\left(1-\left(\frac{\varrho(x,t)}{R}\right)^c\right)^{b}} \dx\dt,
   \end{multline}
where the constant $\big(\frac{(b-1)}{p}c\big)^p$ is sharp.
\end{theorem}
\begin{proof}
Let us start with
    \begin{align*}
       &\frac{(b-1)c}{R} \int_{\br} \frac{|u|^p|\gradg \varrho|^p}{(\varrho/R)^a(1-(\varrho/R)^c)^b}\dx\dt \\& =\frac{R^{Q-1}}{2}\int_{\sn}\int_0^R\frac{|u|^p\psi^{\frac{p}{2}-1}}{(\varrho/R)^{a+(b-1)c-Q}}\pd\bigg(\frac{1}{((\varrho/R)^{-c}-1)^{b-1}}\bigg)\dr\dsn \\&=-\frac{p}{2}R^{Q-1}\int_{\sn}\int_{0}^{R} \frac{|u|^{p-2}u\pd(u)\psi^{\frac{p}{2}-1}}{(\varrho/R)^{a+(b-1)c-Q}}\frac{1}{((\varrho/R)^{-c}-1)^{b-1}}\dr\dsn\\&-\frac{(Q-a-(b-1)c)}{2}R^{Q-2}\int_{\sn}\int_{0}^{R} \frac{|u|^p\psi^{\frac{p}{2}-1}}{(\varrho/R)^{a+(b-1)c-Q+1}}\frac{1}{((\varrho/R)^{-c}-1)^{b-1}}\dr\dsn\\&=-p\int_{\br} \frac{|u|^{p-2}u\pd(u)\psi^{\frac{p}{2}}}{(\varrho/R)^{a+(b-1)c-1}}\frac{1}{((\varrho/R)^{-c}-1)^{b-1}}\dx\dt\\&-\frac{(Q-a-(b-1)c)}{R}\int_{\br} \frac{|u|^p\psi^{\frac{p}{2}}}{(\varrho/R)^{a+(b-1)c}}\frac{1}{((\varrho/R)^{-c}-1)^{b-1}}\dx\dt\\&=-p\int_{\br} \frac{|u|^{p-2}u\pd(u)\psi^{\frac{p}{2}}}{(\varrho/R)^{a-1}(1-(\varrho/R)^{c})^{b-1}}\dx\dt\\&-\frac{(Q-a-(b-1)c)}{R}\int_{\br} \frac{|u|^p\psi^{\frac{p}{2}}}{(\varrho/R)^{a}(1-(\varrho/R)^{c})^{b-1}}\dx\dt\\&\leq p\bigg(\int_{\br} \frac{|\rgradg u|^p}{(\varrho/R)^{a-p}(1-(\varrho/R)^{c})^{b-p}} \dx\dt\bigg)^{\frac{1}{p}}\bigg(\int_{\br} \frac{|u|^{p}|\gradg\varrho|^p}{(\varrho/R)^{a}(1-(\varrho/R)^{c})^{b}} \dx\dt\bigg)^{1-\frac{1}{p}}\\&-\frac{(Q-a-(b-1)c)}{R}\int_{\br} \frac{|u|^p|\gradg \varrho|^p}{(\varrho/R)^{a}(1-(\varrho/R)^{c})^{b-1}}\dx\dt. 
    \end{align*}
    Therefore, dividing by $\big(\int_{\br} \frac{|u|^{p}|\gradg\varrho|^p}{(\varrho/R)^{a}(1-(\varrho/R)^{c})^{b}} \dx \, \dt\big)^{1-\frac{1}{p}}$ for the non-identically zero function, neglecting the negative term, and adjusting the radius $R$ we deduce the result \eqref{two-wg-hardy-eq}. Also, due to the compact support and condition $a< Q$, we have that the boundary term vanishes during integration by parts.
 
    Now, we verify the optimality of the constant. Take $\mu>\frac{b-1}{p}$, and small positive $\delta$ such that $2\delta\in (0,R)$. Let us define the following radially symmetric function
    \begin{align*}
        f_\mu(x,t)=\chi(\varrho(x,t)) \left(1-\left(\frac{\varrho(x,t)}{R}\right)^c\right)^\mu,
    \end{align*}
    where the cutoff function $\chi:[0,\infty)\rightarrow \mathbb{R}$ satisfies the following properties:
\begin{itemize}
    \item[1.] $\chi(r)\in [0,1]$ for all $r\in [0,\infty)$ and $\chi$ is smooth and compactly supported;
    \item[2.] 
\begin{equation*}
\chi(r)=
\begin{dcases}
1, & R-\delta\leq r\leq R; \\
0, & 0\leq r< R-2\delta. \\
\end{dcases}
\end{equation*}
\end{itemize}
Using the condition $p\leq b$ and support of cutoff functions, we can verify $f_\mu\in W^{1,p}_0(\br)$, for each $\mu$. Then, using polar coordinates and Lemma~\ref{i-alpha} for $\alpha=\frac{n+p}{2}-1>0$, we deduce
\begin{align*}
    \int_{\br} &\frac{|f_\mu|^{p}|\gradg \varrho|^p}{\varrho^{a}(1-\frac{\varrho^{c}}{R^c})^{b}} \dx\dt  = \frac{1}{2}\int_{\sn}\psi^{\frac{p}{2}-1}\int_{R-2\delta}^{R} \frac{|\chi(\varrho) (1-\frac{\varrho^{c}}{R^c})^\mu|^{p}}{\varrho^{1+a-Q}(1-\frac{\varrho^{c}}{R^c})^{b}} \dr\dsn \\& \geq \frac{1}{2}\int_{\sn}\psi^{\frac{p}{2}-1}\int_{R-\delta}^{R} \varrho^{Q-1-a}(1-\frac{\varrho^{c}}{R^c})^{\mu p-b}\dr\dsn\\&=\mathcal{I}_{\frac{n+p}{2}-1}\bigg(-\frac{R^c}{c}\bigg)\int_{R-\delta}^{R} \varrho^{Q-c-a}(1-\frac{\varrho^{c}}{R^c})^{\mu p-b}\pd(1-\frac{\varrho^{c}}{R^c})\\&\geq \mathcal{I}_{\frac{n+p}{2}-1}\frac{R^c}{c(\mu p-b+1)}\max\{R^{Q-c-a},(R-\delta)^{Q-c-a}\}\bigg[(1-\frac{\varrho^{c}}{R^c})^{\mu p-b+1}\bigg]^{R-\delta}_R\\&=\mathcal{I}_{\frac{n+p}{2}-1}\frac{R^c}{c(\mu p-b+1)}\max\{R^{Q-c-a},(R-\delta)^{Q-c-a}\}\bigg[1-\frac{(R-\delta)^{c}}{R^c}\bigg]^{\mu p-b+1},
    \end{align*}
which tends to $+\infty$ as $\mu \rightarrow \left(\frac{b-1}{p}\right)^+$. Therefore,
\begin{align*}
  \lim_{\mu\rightarrow \left(\frac{b-1}{p}\right)^+}  \int_{\br} \frac{|f_\mu|^{p}|\gradg \varrho|^p}{\varrho^{a}(1-(\varrho/R)^{c})^{b}} \dx\dt =+\infty.
\end{align*}
On the other hand, we have
\begin{align*}
  \pd( f_\mu)= \chi^\prime(\varrho)(1-(\varrho/R)^c)^\mu-\mu c\chi(\varrho)(1-(\varrho/R)^c)^{\mu-1}\varrho^{c-1}R^{-c} .
\end{align*}
Therefore, we establish 
\begin{multline*}
     \int_{\br} \frac{|\chi^\prime(\varrho)(1-(\varrho/R)^c)^\mu|^p}{\varrho^{a-p}(1-(\varrho/R)^{c})^{b-p}}\psi^{\frac{p}{2}}\dx\dt=\mathcal{I}_{\frac{n+p}{2}-1}\int_{R-2\delta}^{R-\delta}|\chi^\prime(\varrho)|^p\varrho^{Q-a+p-1}(1-(\varrho/R)^c)^{\mu p-b+p}\dr=O(1),
\end{multline*}
and 
\begin{multline*}
    (\mu c)^p\int_{\br} \frac{|\chi(\varrho)(1-(\varrho/R)^c)^{\mu-1}\varrho^{c-1}R^{-c}|^p}{\varrho^{a-p}(1-(\varrho/R)^{c})^{b-p}}\psi^{\frac{p}{2}}\dx\dt\\=(\mu c)^p\int_{\br} (\varrho/R)^{pc}\frac{|\chi(\varrho)(1-(\varrho/R)^c)^{\mu}|^p}{\varrho^{a}(1-(\varrho/R)^{c})^{b}}\psi^{\frac{p}{2}}\dx\dt\leq (\mu c)^p \int_{\br} \frac{|f_\mu|^{p}|\gradg \varrho|^p}{\varrho^{a}(1-(\varrho/R)^{c})^{b}} \dx\dt
\end{multline*}
for $\mu \rightarrow \left(\frac{b-1}{p}\right)^+$. That is, 
\begin{align*}
    \lim_{\mu \rightarrow \left(\frac{b-1}{p}\right)^+}\frac{\int_{\br} \frac{|\rgradg f_\mu(x,t)|^p}{\varrho^{a-p}(x,t)\left(1-\left(\frac{\varrho(x,t)}{R}\right)^c\right)^{b-p}} \dx\dt}{\int_{\br} \frac{|f_\mu(x,t)|^{p}|\gradg \varrho(x,t)|^p}{\varrho^{a}(x,t)\left(1-\left(\frac{\varrho(x,t)}{R}\right)^c\right)^{b}} \dx\dt}=\lim_{\mu \rightarrow \left(\frac{b-1}{p}\right)^+} (\mu c)^p= \bigg(\frac{(b-1)}{p}c\bigg)^p.
    \end{align*}
\end{proof}
\begin{remark}
It is worth noting that our result, Theorem~\ref{two-wg-hardy}, with $a = b = p$ and $c = 1$, can be compared to~\cite[Corollary~1.1]{bad} with $\gamma = 1$.
\end{remark}

We now present the weighted critical version of the Hardy inequality, arising as a limiting case, and establish the sharpness of the constant.
 \begin{corollary}\label{lim-tw-wg-hardy}
  Let $G$ be a Grushin space with dimension $Q$ with $Q\geq 3$. Assume $\br$ be the $\varrho-$gauge ball with radius $R>0$.    Let $1< p<\infty$, $p\leq b<\infty$, and $a\leq Q$. Then for any $u\in C_0^\infty (\br)$ there holds 
    \begin{align}\label{lim-tw-wg-hardy-eq}
       \int_{\br} \frac{|\rgradg u(x,t)|^p}{\varrho^{a-p}(x,t)(\ln \frac{R}{\varrho(x,t)})^{b-p}} \dx\dt\geq  \bigg(\frac{b-1}{p}\bigg)^p\int_{\br} \frac{|u(x,t)|^{p}|\gradg \varrho(x,t)|^p}{\varrho^{a}(x,t)(\ln \frac{R}{\varrho(x,t)})^{b}} \dx\dt,
    \end{align}   
    where the constant $\big(\frac{b-1}{p}\big)^p$ is sharp. 
 \end{corollary}
 \begin{proof}
     We know that $1-t^\alpha=\alpha\ln \frac{1}{t}+ o(\alpha)$ as $\alpha\rightarrow 0$. Using this in Theorem~\ref{two-wg-hardy}, we deduce the required inequality. Now, we will establish that the constant is optimal using a minimizing sequence. Take $\mu>\frac{b-1}{p}$, and small positive $\delta$ such that $2\delta\in (0,R)$. Let us define the following radially symmetric function
    \begin{align*}
        f_\mu(x,t)=\chi(\varrho(x,t)) \left(\ln \frac{R}{\varrho(x,t)}\right)^\mu,
    \end{align*}
    where the cutoff function $\chi:[0,\infty)\rightarrow \mathbb{R}$ satisfies the following properties:
\begin{itemize}
    \item[1.] $\chi(r)\in [0,1]$ for all $r\in [0,\infty)$ and $\chi$ is smooth and compactly supported;
    \item[2.] 
\begin{equation*}
\chi(r)=
\begin{dcases}
1, & R-\delta\leq r\leq R; \\
0, & 0\leq r< R-2\delta. \\
\end{dcases}
\end{equation*}
\end{itemize}
It is easy to verify that $f_\mu\in W^{1,p}_0(\br)$. Then, using Lemma~\ref{i-alpha} for $\alpha=\frac{n+p}{2}-1>0$ and following a similar program like showing the optimality of the constant in Theorem~\ref{two-wg-hardy}, it turns out 
\begin{align*}
    \lim_{\mu\rightarrow \left(\frac{b-1}{p}\right)^+}\frac{ \int_{\br} \frac{|\rgradg f_{\mu}(x,t)|^p}{\varrho^{a-p}(x,t)(\ln \frac{R}{\varrho(x,t)})^{b-p}} \dx\dt}{\int_{\br} \frac{|f_\mu(x,t)|^{p}|\gradg \varrho(x,t)|^p}{\varrho^{a}(x,t)(\ln \frac{R}{\varrho(x,t)})^{b}} \dx\dt}=\lim_{\mu\rightarrow\left(\frac{b-1}{p}\right)^+} \mu^p= \bigg(\frac{b-1}{p}\bigg)^p.
\end{align*}
 \end{proof}
\begin{remark}
The above Corollary~\ref{lim-tw-wg-hardy} remains valid for a more general weight function. If we compare~\eqref{lim-tw-wg-hardy-eq} with $a = p$ and $b=-\alpha$, it recovers the version of \cite[Corollary~4.12]{bad} with $\gamma=1$. However, in our case, we allow a more general weight function on the left-hand side of the inequality.
\end{remark}

The geometric Hardy inequality is another important functional inequality in the literature. For the Grushin operator, this appears to be the first instance demonstrating a geometric Hardy inequality. For any $(x,t)\in \br$, we denote the distance function by $\text{dist}(x,\partial \br) = R - \varrho(x,t)$. By substituting $a = p$ and $c = 1$ in Theorem~\ref{two-wg-hardy}, we obtain the following remark.
\begin{remark}
     Let $G$ be a Grushin space with dimension $Q$ with $Q\geq 3$. Assume $\br$ be the $\varrho-$gauge ball with radius $R>0$. Let $1< p<Q$, $p\leq b<\infty$, and $p\leq Q-b+1$. Then for any $u\in C_0^\infty (\br)$ there holds
\begin{align*}
       \int_{\br} \frac{|u(x,t)|^{p}|\gradg\varrho(x,t)|^p}{\text{dist}(x,\partial \br)^b} \dx\dt &\leq  \int_{\br} \frac{|u(x,t)|^{p}|\gradg \varrho(x,t)|^p}{\left(\frac{\varrho(x,t)}{R}\right)^p(R-\varrho(x,t))^{b}} \dx\dt\\&\leq \bigg(\frac{p}{b-1}\bigg)^p \int_{\br} \frac{|\rgradg u(x,t)|^p}{\text{dist}(x,\partial \br)^{b-p}} \dx\dt\\&\leq \bigg(\frac{p}{b-1}\bigg)^p \int_{\br} \frac{|\gradg u(x,t)|^p}{\text{dist}(x,\partial \br)^{b-p}} \dx\dt,
    \end{align*}  
where the constant  $\big(\frac{p}{b-1}\big)^p$ is sharp. In the end, we used pointwise estimates $|\rgradg u(x,t)|\leq |\gradg u(x,t)|$ for every $(x,t)\in \rno$. 
\end{remark}

In Corollary~\ref{lim-tw-wg-hardy}, we established a critical Hardy inequality for functions supported in the ball $\br$. Also, notice that in Corollary~\ref{lim-tw-wg-hardy}, we cannot push the parameter $b$ up to $0$. However, one can also prove a different form of a critical Hardy inequality for functions supported on the entire Grushin space $G$, where logarithmic terms only appear on the left-hand side. For homogeneous groups, this phenomenon was addressed in~\cite[Theorem~4.1]{rs-fenn}. We conclude this section by demonstrating the same phenomenon in the Grushin setting.
 \begin{theorem}
Let $G$ be a Grushin space with dimension $Q$ with $Q\geq 3$. Assume $1<p<\infty$. Then for any $u\in C_0^\infty (\rno \setminus\{o\} )$ there holds
\begin{equation}\label{anoth-crit-hardy-eq}
   \int_{\rno} \frac{|\ln \varrho(x,t)|^p\big|\rgradg u(x,t)\big|^p}{\varrho^{Q-p}(x,t)}\dx\dt\geq \frac{1}{p^p} \int_{\rno}\frac{|u(x,t)|^p |\gradg\varrho(x,t)|^p}{\varrho^Q(x,t)}\dx\dt,
\end{equation}
where the constant $\frac{1}{p^p}$ is sharp.
\end{theorem}
\begin{proof}
Let us start with a function $u\in C_0^\infty(\rno \setminus\{o\})$. Then there exists a sufficiently large $R>0$ such that $\text{supp}(u) \subset \br$. Using polar coordinates and integration by parts, we obtain
\begin{align*}
    &\int_{\rno}\frac{|u(x,t)|^p |\gradg\varrho(x,t)|^p}{\varrho^Q(x,t)}\dx\dt\\&=\frac{1}{2}\int_{\Omega} \psi^{\frac{p}{2}-1}\int_0^R |u|^p \pd(\ln \varrho)\dr\dsn\\&=-\frac{p}{2}\int_\Omega \psi^{\frac{p}{2}-1}\int_0^R |u|^{(p-2)}u\pd(u)(\ln\varrho)\dr\dsn\\&\leq p\int_{\br}\frac{|\rgradg u||\ln \varrho|}{\varrho^{\frac{Q-p}{p}}}\frac{|u|^{(p-1)}\psi^{\frac{p-1}{2}}}{\varrho^{\frac{Q(p-1)}{p}}}\dx\dt\\
       &\leq p\Bigg(\int_{\br} \frac{|\ln \varrho(x,t)|^p\big|\rgradg u(x,t)\big|^p}{\varrho^{Q-p}(x,t)}\dx\dt\Bigg)^{\frac{1}{p}}
       \Bigg( \int_{\br}\frac{|u(x,t)|^p |\gradg\varrho(x,t)|^p}{\varrho^Q(x,t)}\dx\dt \Bigg)^{\frac{p-1}{p}},
\end{align*}
and by arranging both the terms, we obtain \eqref{anoth-crit-hardy-eq}. 

Next, we prove the constant $\frac{1}{p^p}$ is optimal. For each small $\delta>0$ such that $\delta<\min\{\frac{1}{4},\frac{1}{p}\}$, we define the radially symmetric function  $\{f_\delta\}$ as follows
\begin{equation*}
f_{\delta}(x,t)= f_\delta(\varrho, \sigma):=
\begin{dcases}
4\left(\ln 4 \right)^{-\frac{1}{p}+\delta}\frac{(\varrho-\delta)}{(1-4\delta)} & \text{ when } \delta\leq \varrho\leq  \frac{1}{4}; \\
\left(-\ln \varrho\right)^{-\frac{1}{p}+\delta} & \text{ when }  \frac{1}{4} \leq \varrho\leq  \frac{1}{1+\delta};\\
\frac{(1+\delta)}{\delta}\left(\ln (1+\delta) \right)^{-\frac{1}{p}+\delta}(1-\varrho) &  \text{ when }\frac{1}{1+\delta}\leq \varrho\leq  1; \\
0 &  \text{ when } 0\leq \varrho \leq \delta \text{ or } 1\leq \varrho <  \infty.
\end{dcases}
\end{equation*}
It is easy to verify
\begin{equation*}
 \partial_{\varrho} f_\delta(\varrho, \sigma)=
\begin{dcases}
\frac{4}{\left(1-4\delta\right)}\left(\ln 4 \right)^{-\frac{1}{p}+\delta} & \text{ when } \delta< \varrho<  \frac{1}{4}; \\
-\left(-\frac{1}{p}+\delta\right)\left(-\ln \varrho\right)^{-\frac{1}{p}+\delta-1}\frac{1}{\varrho} & \text{ when }  \frac{1}{4} < \varrho<  \frac{1}{1+\delta};\\
-\frac{(1+\delta)}{\delta}\left(\ln (1+\delta) \right)^{-\frac{1}{p}+\delta} &  \text{ when }\frac{1}{1+\delta}< \varrho<  1; \\
0 &  \text{ when } 0< \varrho < \delta \text{ or } 1< \varrho <  \infty.
\end{dcases} 
\end{equation*}
and therefore for each fixed $\delta>0$, $f_\delta\in W_0^{1,p}(\rno\setminus \{o\})$. Now, using polar coordinates and Lemma~\ref{i-alpha} for $\alpha=\frac{n+p}{2}-1>0$, we get
\begin{align*}
    \int_{\rno} \frac{|f_\delta(x,t)|^p|\gradg\varrho(x,t)|^p}{\varrho^{Q}(x,t)}\dx \dt \geq& \mathcal{I}_{\frac{n+p}{2}-1} \int_{1/4}^{1/(1+\delta)}(-\ln \varrho)^{-1+\delta p}\varrho^{-1}\dr\\
    =&\mathcal{I}_{\frac{n+p}{2}-1}\left[\frac{(\ln 4)^{\delta p}-(\ln (1+\delta))^{\delta p}}{\delta p}\right].
\end{align*}
Since
\begin{align*}
    \lim_{\delta \rightarrow 0^+}\frac{(\ln 4)^{\delta p}-(\ln (1+\delta))^{\delta p}}{\delta p}=+\infty, 
\end{align*}
we get
\begin{align*}
    \int_{\rno} \frac{|f_\delta(x,t)|^p|\gradg\varrho(x,t)|^p}{\varrho^{Q}(x,t)}\dx \dt \rightarrow +\infty \text{ as } \delta\rightarrow 0^+.
\end{align*}
On the other hand, we have 
\begin{align*}
   & \int_{\rno} \frac{|\ln \varrho|^p\big|\rgradg f_{\delta}\big|^p}{\varrho^{Q-p}}\dx\dt\\& \leq \left(\frac{1}{p}-\delta\right)^{p} \int_{ B_{\frac{1}{1+\delta}}^{Q}(o)\setminus B_{\frac{1}{4}}^{Q}(o)}\frac{|\ln \varrho|^{-1+\delta p}}{\varrho^{Q}} |\gradg \varrho|^p \dx\dt\\& +\frac{4^p}{(1-4\delta)^p}\left(\ln 4 \right)^{-1+\delta p} \int_{B_{\frac{1}{4}}^{Q}(o)\setminus B_{\delta}^{Q}(o)}\frac{|\ln \varrho|^{p}}{\varrho^{Q-p}}|\gradg \varrho|^p\dx\dt
   \\&+\frac{(1+\delta)^p}{\delta^p}\left(\ln (1+\delta) \right)^{-1+\delta p} \int_{B_{1}^{Q}(o)\setminus B_{\frac{1}{1+\delta}}^{Q}(o)}\frac{|\ln \varrho|^{p}}{\varrho^{Q-p}}|\gradg \varrho|^p\dx\dt \\
   &=\left(\frac{1}{p}-\delta\right)^{p} \int_{ B_{\frac{1}{1+\delta}}^{Q}(o)\setminus B_{\frac{1}{4}}^{Q}(o)}\frac{|f_\delta|^{p}}{\varrho^{Q}}|\gradg \varrho|^p\dx\dt+F_1+F_2,
\end{align*}
where 
\begin{align*}
    F_1&:=\frac{4^p}{(1-4\delta)^p}\left(\ln 4 \right)^{-1+\delta p} \int_{B_{\frac{1}{4}}^{Q}(o)\setminus B_{\delta}^{Q}(o)}\frac{|\ln \varrho|^{p}}{\varrho^{Q-p}}|\gradg \varrho|^p\dx\dt\\&= \frac{4^p}{(1-4\delta)^p}\left(\ln 4 \right)^{-1+\delta p} \mathcal{I}_{\frac{n+p}{2}-1} \int_{\delta}^{1/4}\varrho^{p-1}(-\ln \varrho)^p\dr = O(1),
\end{align*}
and 
\begin{align*}
    F_2&:=\frac{(1+\delta)^p}{\delta^p}\left(\ln (1+\delta) \right)^{-1+\delta p} \int_{B_{1}^{Q}(o)\setminus B_{\frac{1}{1+\delta}}^{Q}(o)}\frac{|\ln \varrho|^{p}}{\varrho^{Q-p}}|\gradg \varrho|^p\dx\dt\\&= \frac{(1+\delta)^p}{\delta^p}\left(\ln (1+\delta) \right)^{-1+\delta p} \mathcal{I}_{\frac{n+p}{2}-1}\int_{\frac{1}{1+\delta}}^{1}\varrho^{p-1}(-\ln \varrho)^p\dr \\& \leq \mathcal{I}_{\frac{n+p}{2}-1}\frac{(1+\delta)^p}{\delta^p(p+1)}\left(\ln (1+\delta) \right)^{p+\delta p}=O(1).
\end{align*}
For the term $F_1$, we used the fact $\lim_{\varrho\rightarrow 0^+}\varrho^{p-1}(-\ln \varrho)^p=0$, where $p>1$. For $F_2$, we used the fact $ \lim_{\delta \rightarrow 0^+}\frac{(1+\delta)^p}{\delta^p}\left(\ln (1+\delta) \right)^{p+\delta p}=1$. Hence, combining all the above estimates, we conclude
\begin{align*}
    \lim_{\delta\rightarrow 0^+}\frac{\int_{\rno} \frac{|\ln \varrho(x,t)|^p\big|\rgradg f_{\delta}(x,t)\big|^p}{\varrho^{Q-p}(x,t)}\dx\dt}{\int_{\rno}\frac{|f_{\delta}(x,t)|^p |\gradg\varrho(x,t)|^p}{\varrho^Q(x,t)}\dx\dt}=  \frac{1}{p^p},
\end{align*}
which gives the desired optimality.
\end{proof}

\section{Hardy--Rellich  and Rellich inequalities}\label{r-hr-est}
The establishment of Hardy--Rellich and Rellich inequalities for radial Grushin operators is the main goal of this section. Our initial objective is to address a particular type of inequality involving only the first-order radial partial derivative $\pd$. The work \cite[Lemma~4.1]{vhn} appears to be the first to develop such an approach, which we adapt here to the Grushin setting.
 \begin{theorem}\label{rellich-type}
        Let $G$ be a Grushin space with dimension $Q$ with $Q\geq 3$. Let $1<p<Q$ and $-Q(p-1)<\beta<Q-p$. Then for $u\in C_0^\infty(\rno)$ there holds
        \begin{align}\label{rellich-type-eq}
            \int_{\rno}\frac{\left|\pd u(x,t)+\frac{(Q-1)}{\varrho(x,t)}u(x,t)\right|^p}{\varrho^{\beta}(x,t)}\dx\dt\geq \left(\frac{Q(p-1)+\beta}{p}\right)^p\int_{\rno}\frac{|u(x,t)|^p}{\varrho^{p+\beta}(x,t)}\dx\dt,
        \end{align}
        where the constant $\big(\frac{Q(p-1)+\beta}{p}\big)^p$ is best possible.
    \end{theorem}
    \begin{proof}
        The polar coordinate decomposition, the integration by parts formula, and the condition $Q - p - \beta > 0$ yield
        \begin{align*}
           \int_{\rno}&\frac{|u|^p}{\varrho^{p+\beta}}\dx\dt= \int_{\Omega}\int_{0}^{\infty}\frac{|u|^p}{\varrho^{p+\beta}}\frac{\varrho^{n+1}}{2\psi}\dr\dsn\\&=-\frac{p}{2(Q-p-\beta)}\int_{\Omega}\int_{0}^{\infty}|u|^{p-2}u\,\pd(u)\psi^{-1}\varrho^{-p-\beta+n+2}\dr\dsn\\&=-\frac{p}{(Q-p-\beta)}\int_{\rno}\frac{\left[\pd (u)+\frac{(Q-1)}{\varrho}u\right]}{\varrho^{\frac{\beta}{p}}}\frac{|u|^{p-2}u}{\varrho^{\frac{(\beta+p)(p-1)}{p}}}\dx\dt+\frac{p(Q-1)}{(Q-p-\beta)}\int_{\rno}\frac{|u|^p}{\varrho^{\beta+p}}\dx\dt.
        \end{align*}
        After arranging terms, this implies  
        \begin{align*}
            \left(\frac{Q(p-1)+\beta}{p}\right)&\int_{\rno}\frac{|u|^p}{\varrho^{p+\beta}}\dx\dt= \int_{\rno}\frac{\left[\pd (u)+\frac{(Q-1)}{\varrho}u\right]}{\varrho^{\frac{\beta}{p}}}\frac{|u|^{p-2}u}{\varrho^{\frac{(\beta+p)(p-1)}{p}}}\dx\dt,
        \end{align*}
        and applying H\"older inequality result \eqref{rellich-type-eq} follows.

        Now, it remains to show the optimality of the constant. Choose $0<\epsilon<\frac{1}{2}$. Take a cutoff function $\chi:\rno\rightarrow [0,1]$, which is radial, and it is one on $B_{\frac{1}{2}}(o)$ and zero outside $B_{1}(o)$. Now define the minimizing sequence by $H_{\epsilon}(\varrho)=\left(1-\chi(\varrho/\epsilon)\right)\chi(\varrho)\varrho^{-\frac{Q-\beta-p}{p}}$ and this clearly belongs to $C_0^\infty(\rno\setminus\{o\})\subset C_0^\infty(\rno)$. Now using Lemma~\ref{i-alpha} for $\alpha=\frac{n}{2}-1\geq 0$, we obtain
        \begin{align*}
            \int_{\rno}\frac{|H_{\epsilon}(\varrho)|^p}{\varrho^{p+\beta}}\dx\dt\geq \int_{B_{\frac{1}{2}}^Q(o)\setminus B_\epsilon^Q(o)}\frac{|H_{\epsilon}(\varrho)|^p}{\varrho^{p+\beta}}\dx\dt=\mathcal{I}_{\frac{n}{2}-1}\int_{\epsilon}^{\frac{1}{2}} \varrho^{-1}\dr=\mathcal{I}_{\frac{n}{2}-1}\ln\left(\frac{1}{2\epsilon}\right)\rightarrow +\infty
        \end{align*}
        for $\epsilon\rightarrow 0^+$. On the other hand, we have 
        \begin{multline*}
            \pd H_\epsilon +\frac{(Q-1)}{\varrho}H_\epsilon=-\epsilon^{-1}\chi'(\varrho/\epsilon)\chi(\varrho)\varrho^{-\frac{Q-\beta-p}{p}}\\+\left(1-\chi(\varrho/\epsilon)\right)\chi'(\varrho)\varrho^{-\frac{Q-\beta-p}{p}}+\frac{(Q(p-1)+\beta)}{p}\left(1-\chi(\varrho/\epsilon)\right)\chi(\varrho)\varrho^{-\frac{Q-\beta-p}{p}-1}.
        \end{multline*}
        We can check now 
        \begin{align*}
        \int_{\rno}\frac{\left|\epsilon^{-1}\chi'(\varrho/\epsilon)\chi(\varrho)\varrho^{-\frac{Q-\beta-p}{p}}\right|^p}{\varrho^{\beta}(x,t)}\dx\dt&\leq \mathcal{I}_{\frac{n}{2}-1}\int_{0}^{\infty}\epsilon^{-p}\left|\chi'(\varrho/\epsilon)\right|^p\varrho^{p-1}\dr\\&=\mathcal{I}_{\frac{n}{2}-1}\int_{0}^{\infty}\left|\chi'(s)\right|^p s^{p-1}\ds=O(1),
        \end{align*}
        \begin{align*}
        \int_{\rno}\frac{\left|\left(1-\chi(\varrho/\epsilon)\right)\chi'(\varrho)\varrho^{-\frac{Q-\beta-p}{p}}\right|^p}{\varrho^{\beta}(x,t)}\dx\dt=\int_{\ba\setminus B_{\frac{1}{2}}^Q(o)}\frac{\left|\left(1-\chi(\varrho/\epsilon)\right)\chi'(\varrho)\varrho^{-\frac{Q-\beta-p}{p}}\right|^p}{\varrho^{\beta}(x,t)}\dx\dt=O(1),
        \end{align*}
        for $\epsilon\rightarrow 0^+$ and 
        \begin{align*}
        \int_{\rno}\frac{\left|\left(1-\chi(\varrho/\epsilon)\right)\chi(\varrho)\varrho^{-\frac{Q-\beta-p}{p}-1}\right|^p}{\varrho^{\beta}(x,t)}\dx\dt=\int_{\rno}\frac{|H_{\epsilon}(\varrho)|^p}{\varrho^{p+\beta}}\dx\dt.
        \end{align*}
        Therefore, we deduce
        \begin{align*}
            \lim_{\epsilon\rightarrow 0^+}\frac{ \int_{\rno}\frac{\left|\pd H_{\epsilon}(\varrho)+\frac{(Q-1)}{\varrho}H_{\epsilon}(\varrho)\right|^p}{\varrho^{\beta}}\dx\dt}{\int_{\rno}\frac{|H_{\epsilon}(\varrho)|^p}{\varrho^{p+\beta}}\dx\dt}= \left(\frac{Q(p-1)+\beta}{p}\right)^p.
        \end{align*}
    \end{proof}

Now, we are ready to present the Hardy--Rellich inequality in the Grushin setting with the full radial operator.
     \begin{theorem}\label{Hardy--Rellich}
        Let $G$ be a Grushin space with dimension $Q$ with $Q\geq 3$. Let $1<p<Q$ and $-Q(p-1)<\beta<Q-p$. Then for $u\in C_0^\infty(\rno)$ there holds
        \begin{align}\label{Hardy--Rellich-eq}
            \int_{\rno}\frac{\left|\rlapg u(x,t)\right|^p}{|\gradg \varrho(x,t)|^p\varrho^{\beta}(x,t)}\dx\dt\geq \left(\frac{Q(p-1)+\beta}{p}\right)^p\int_{\rno}\frac{|\rgradg u(x,t)|^p}{\varrho^{p+\beta}(x,t)}\dx\dt,
        \end{align}
        where the constant $\big(\frac{Q(p-1)+\beta}{p}\big)^p$ is best possible.
    \end{theorem}
    \begin{proof}
        Replacing $u$ by $\pd \big(\psi^{\frac{1}{2}} u\big)$ in Theorem~\ref{rellich-type}, and using the independence of $\psi$ w.r.t. $\varrho$, we immediately obtain \eqref{Hardy--Rellich-eq}. The constant is sharp; otherwise, suppose \eqref{Hardy--Rellich-eq} holds for $u$ with some bigger constant, and then we will get a contradiction of the optimality of \eqref{rellich-type-eq} with the smooth compactly supported function $\psi^{\frac{1}{2}}u_\varrho$.
    \end{proof}

    It is well known that Rellich-type inequalities can be obtained by iterating the Hardy-type inequalities. Also, we can obtain the critical higher-order version in the inductive steps using the critical version \eqref{lim-tw-wg-hardy-eq} with $b=p$, in the inductive steps. Here, for brevity, we are only mentioning the subcritical results for higher-order derivatives in the Grushin setting. We also refer to recent work \cite[Corollary~4.4]{darca-1} in this context. 

    \begin{theorem}\label{high-hr}
       Let $G$ be a Grushin space with dimension $Q$ with $Q\geq 3$. Assume an integer $k\geq 1$ that is less than $Q$. Let $1<p<\frac{Q}{k}$ and $\beta\in\mathbb{R}$. Then for any $u\in C_0^\infty(\rno)$ the following holds:
       \begin{itemize}
           \item[(i)] If $k=2l$, $l\geq 1$ and $-(2l-2)p-Q(p-1)<\beta<Q-2lp$, then we have 
           \begin{align*}
            \int_{\rno}\frac{\left|\rlapg^l u(x,t)\right|^p}{|\gradg \varrho(x,t)|^{(2l-1)p}\varrho^{\beta}(x,t)}\dx\dt\geq C_{k,p,\beta}\int_{\rno}\frac{|u(x,t)|^p|\gradg \varrho(x,t)|^p}{\varrho^{2lp+\beta}(x,t)}\dx\dt.
        \end{align*}
        \item[(ii)] If $k=2l+1$, $l\geq 0$ and $-(2l-1)p-Q(p-1)<\beta<Q-(2l+1)p$, then we have 
           \begin{align*}
            \int_{\rno}\frac{\left|\rgradg \rlapg^l u(x,t)\right|^p}{|\gradg \varrho(x,t)|^{2lp}\varrho^{\beta}(x,t)}\dx\dt\geq C_{k,p,\beta}\int_{\rno}\frac{|u(x,t)|^p|\gradg \varrho(x,t)|^p}{\varrho^{(2l+1)p+\beta}(x,t)}\dx\dt.
            \end{align*}
       \end{itemize}
       Moreover, all the constants are sharp.
    \end{theorem}
    \begin{proof}
       First, the case $k=1$ is the subcritical Hardy inequality \eqref{sub-hardy-eq-intro}. Next, we will establish the case $k=2$, which is the classical Rellich inequality for the Grushin setting. Start with $u\in C_0^\infty(\rno)$. Then from \eqref{sub-hardy-eq-intro} with replaced weight $\beta$ by $\beta+p$, we have
     \begin{align*}
            \int_{\rno}\frac{|\rgradg u(x,t)|^p}{\varrho^{\beta+p}(x,t)}\dx\dt\geq \left(\frac{Q-\beta-2p}{p}\right)^p\int_{\rno}\frac{|u(x,t)|^p|\gradg \varrho(x,t)|^p}{\varrho^{2p+\beta}(x,t)}\dx\dt.
        \end{align*}
        Using this in \eqref{Hardy--Rellich}, we deduce
        \begin{multline}\label{rellich-eq}
             \int_{\rno}\frac{\left|\rlapg u(x,t)\right|^p}{|\gradg \varrho(x,t)|^p\varrho^{\beta}(x,t)}\dx\dt \geq \left(\frac{(Q(p-1)+\beta)(Q-\beta-2p)}{p^2}\right)^p\int_{\rno}\frac{|u(x,t)|^p|\gradg \varrho(x,t)|^p}{\varrho^{2p+\beta}(x,t)}\dx\dt.
        \end{multline}
        
       Now we will consider the case $k=3$, i.e., $l=1$. Using \eqref{sub-hardy-eq-intro} for the function $\frac{\rlapg u(x,t)}{|\gradg \varrho (x,t)|^2}$, we deduce
       \begin{align*}
           \int_{\rno}\frac{\left|\rgradg \rlapg u(x,t)\right|^p}{|\gradg \varrho(x,t)|^{2p}\varrho^{\beta}(x,t)}\dx\dt\geq \Lambda_{1}\int_{\rno}\frac{|\rlapg u(x,t)|^p}{|\gradg \varrho(x,t)|^{p}\varrho^{p+\beta}(x,t)}\dx\dt
       \end{align*}
       and applying \eqref{rellich-eq} for the weight $p+\beta$ instead of $\beta$, we deduce
       \begin{align*}
           \int_{\rno}\frac{\left|\rgradg \rlapg u(x,t)\right|^p}{|\gradg \varrho(x,t)|^{2p}\varrho^{\beta}(x,t)}\dx\dt\geq \Lambda_{1}\Lambda_3\int_{\rno}\frac{|u(x,t)|^p|\gradg \varrho(x,t)|^p}{\varrho^{3p+\beta}(x,t)}\dx\dt,
       \end{align*}
       which completes the case.

  Next, we will show the case $k=4$, where $l=2$. Take the function $\frac{\rlapg u(x,t)}{|\gradg \varrho (x,t)|^2}$ and apply it in \eqref{rellich-eq}, we deduce
       \begin{align*}
            \int_{\rno}\frac{\left|\rlapg^2 u(x,t)\right|^p}{|\gradg \varrho(x,t)|^{3p}\varrho^{\beta}(x,t)}\dx\dt\geq  \Lambda_{2}\int_{\rno}\frac{|\rlapg u(x,t)|^p}{|\gradg \varrho(x,t)|^p\varrho^{2p+\beta}(x,t)}\dx\dt,
        \end{align*}
        and again applying \eqref{rellich-eq} with weight $2p+\beta$, we have 
        \begin{align*}
            \int_{\rno}\frac{\left|\rlapg^2 u(x,t)\right|^p}{|\gradg \varrho(x,t)|^{3p}\varrho^{\beta}(x,t)}\dx\dt\geq  \Lambda_{2}\Lambda_{4}\int_{\rno}\frac{| u(x,t)|^p|\gradg \varrho(x,t)|^p}{\varrho^{4p+\beta}(x,t)}\dx\dt,
        \end{align*}
        and this completes $k=4$ case. Also, in the above, we used the fact that $|\gradg \varrho(x,t)|=\psi(x,t)$ is independent of the homogeneous distance function $\varrho$. The rest of the proof follows from induction.

 Now, it remains to show the optimality of the constant. Let $\delta>0$ and choose $0<\epsilon<\frac{\delta}{2}$. Take a cutoff function $\varphi:\rno\rightarrow [0,1]$, which is radial and it is one on $B_{\frac{1}{2}}(o)$ and zero outside $B_{1}(o)$. First, consider the case $k=2l$ for $l\geq 1$. Now define the minimizing sequence by $H_{\delta,\epsilon}(\varrho)=\left(1-\varphi(\varrho/\epsilon)\right)\varphi(\varrho/\delta)\varrho^{-\frac{Q-\beta-2lp}{p}}$ and this clearly belongs to $C_0^\infty(B_\delta^Q(o)\setminus B_{\frac{\epsilon}{2}}^Q(o))\subset C_0^\infty(\rno)$. Now using Lemma~\ref{i-alpha} for $\alpha=\frac{n+p}{2}-1>0$, we obtain
        \begin{align*}
            \int_{\rno}\frac{|H_{\delta,\epsilon}(\varrho)|^p|\gradg \varrho|^p}{\varrho^{2lp+\beta}}\dx\dt&\geq \int_{B_{\frac{\delta}{2}}^Q(o)\setminus B_\epsilon^Q(o)}\frac{|H_{\delta,\epsilon}(\varrho)|^p|\gradg \varrho|^p}{\varrho^{2lp+\beta}}\dx\dt\\&=\mathcal{I}_{\frac{n+p}{2}-1}\int_{\epsilon}^{\frac{\delta}{2}} \varrho^{-1}\dr=\mathcal{I}_{\frac{n+p}{2}-1}\ln\left(\frac{\delta}{2\epsilon}\right)\rightarrow +\infty
        \end{align*}
        for $\epsilon\rightarrow 0^+$ and $\delta$ fixed. On the other hand, we have 
        \begin{align*}
            \int_{\rno}\frac{\left|\rlapg^l H_{\delta,\epsilon}(\varrho)\right|^p}{|\gradg \varrho|^{(2l-1)p}\varrho^{\beta}}\dx\dt& = \int_{B_{\frac{\delta}{2}}^Q(o)\setminus B_\epsilon^Q(o)}\frac{\left|\rlapg^l \varrho^{-\frac{Q-\beta-2lp}{p}}\right|^p}{|\gradg \varrho|^{(2l-1)p}\varrho^{\beta}}\dx\dt\\&+\int_{B_{\delta}^Q(o)\setminus B_\frac{\delta}{2}^Q(o)}\frac{\left|\rlapg^l \left(\varphi(\varrho/\delta)\varrho^{-\frac{Q-\beta-2lp}{p}}\right)\right|^p}{|\gradg \varrho|^{(2l-1)p}\varrho^{\beta}}\dx\dt\\&+\int_{B_{\epsilon}^Q(o)\setminus B_\frac{\epsilon}{2}^Q(o)}\frac{\left|\rlapg^l\left(\left(1-\varphi(\varrho/\epsilon)\right) \varrho^{-\frac{Q-\beta-2lp}{p}}\right)\right|^p}{|\gradg \varrho|^{(2l-1)p}\varrho^{\beta}}\dx\dt\\&:=\mathcal{M}_1 +\mathcal{M}_2 +\mathcal{M}_3 
        \end{align*}
        for a fixed $\delta$. We will evaluate each term separately. First applying Lemma~\ref{ind-lemm}, we deduce
        \begin{align*}
            \mathcal{M}_1\leq C_{k,p,\beta} \int_{B_{\frac{\delta}{2}}^Q(o)\setminus B_\epsilon^Q(o)}|\gradg \varrho|^p\varrho^{-Q}\dx\dt=C_{k,p,\beta}\int_{B_{\frac{\delta}{2}}^Q(o)\setminus B_\epsilon^Q(o)}\frac{|H_{\delta,\epsilon}(\varrho)|^p|\gradg \varrho|^p}{\varrho^{2lp+\beta}}\dx\dt.
        \end{align*}
        Next, applying Lemma~\ref{lib-ind}, Minkowski's inequality, we have 
        \begin{align*}
            \mathcal{M}_2^{\frac{1}{p}}&\leq C\sum_{i=0}^{2l}\left(\delta^{-pi}|h_{k,i}|^p \int_{B_{\delta}^Q(o)\setminus B_\frac{\delta}{2}^Q(o)}|\gradg\varrho|^p\varrho^{-Q+ip}\dx\dt \right)^{\frac{1}{p}}\\&\leq C\mathcal{I}_{\frac{n+p}{2}-1}^{\frac{1}{p}}\sum_{i=0}^{2l}\left(\delta^{-pi}|h_{k,i}|^p \int_{\frac{\delta}{2}}^\delta \varrho^{ip-1}\dr \right)^{\frac{1}{p}}\\&=C\mathcal{I}_{\frac{n+p}{2}-1}^{\frac{1}{p}}\left[|h_{k,0}|(\ln 2)^{\frac{1}{p}}+\sum_{i=1}^{2l}|h_{k,i}| (1-2^{-ip})^{\frac{1}{p}}\right]
        \end{align*}
        where $C$ is dependent on the bound of derivatives of $\varphi$ and this gives $\mathcal{M}_2<+\infty$ for any $\delta$ and $\epsilon$. Finally, we have $\mathcal{M}_3<+\infty$ by similar analysis, the bound is independent of all $\epsilon$ and $\delta$. Combining all these, we deduce
        \begin{align*}
           \lim_{\epsilon\rightarrow 0^+}\frac{\int_{\rno}\frac{\left|\rlapg^l H_{\delta,\epsilon}(\varrho)\right|^p}{|\gradg \varrho|^{(2l-1)p}\varrho^{\beta}}\dx\dt}{\int_{\rno}\frac{|H_{\delta,\epsilon}(\varrho)|^p|\gradg \varrho|^p}{\varrho^{2lp+\beta}}\dx\dt} \leq C_{k,p,\beta}
        \end{align*}
        which gives the optimality, and similarly, we can get the odd case.
        \end{proof}

We will conclude this section by establishing the critical Rellich inequality in the Grushin setting.
        \begin{theorem}\label{rellich-crit}
        Let $G$ be a Grushin space with dimension $Q$ with $Q\geq 3$. Assume $\br$ be the $\varrho-$gauge ball with radius $R>0$. Let $1< p<Q$. Then for any $u\in C_0^\infty (\br)$ there holds
        \begin{align*}
            \int_{\br}\frac{\left|\rlapg u(x,t)\right|^p}{|\gradg \varrho(x,t)|^p\varrho^{Q-2p}(x,t)}\dx\dt\geq \bigg(\frac{(Q-2)(p-1)}{p}\bigg)^p\int_{\br} \frac{|u(x,t)|^{p}|\gradg \varrho(x,t)|^p}{\varrho^{Q}(x,t)(\ln \frac{R}{\varrho(x,t)})^{p}} \dx\dt,
        \end{align*}
        where the constant $\big(\frac{(Q-2)(p-1)}{p}\big)^p$ is best possible.
    \end{theorem}
    \begin{proof}
    Using Theorem~\ref{Hardy--Rellich} for the weight $\beta$ replaced by $Q-2p$, we have
    \begin{align*}
        \int_{\br}\frac{\left|\rlapg u(x,t)\right|^p}{|\gradg \varrho|^p\varrho^{Q-2p}(x,t)}\dx\dt\geq \left(Q-2\right)^p\int_{\br}\frac{|\rgradg u(x,t)|^p}{\varrho^{Q-p}(x,t)}\dx\dt,
    \end{align*}
    and using above in \eqref{lim-tw-wg-hardy-eq} for $b=p$, and $a=Q$, we deduce the result. Suppose $\delta$ is a small positive number less than $R$, and consider the sequence of functions $\{f_\delta\}$ defined in \eqref{seq-fun} for $\gamma=p$. Then from \eqref{seq-fun-der} we can estimate the first derivatives, and also we calculate
\begin{equation*}
\psi^{-1}\rlapg f_\delta=
\begin{dcases}
0 &  \text{ when } \varrho<\delta; \\
-\frac{(p-1)(Q-2)}{p}\bigg(\ln\frac{R}{\varrho}\bigg)^{-\frac{1}{p}}\varrho^{-2}-\left(\frac{p^2-1}{p^2}\right)\bigg(\ln\frac{R}{\varrho}\bigg)^{-\frac{1+p}{p}}\varrho^{-2} & \text{ when } \delta < \varrho <  \frac{R}{2}; \\
-\left(\ln 2 \right)^{\frac{p-1}{p}}\frac{2}{R}(Q-1)\varrho^{-1} & \text{ when }\frac{R}{2} < \varrho <  R. \\
\end{dcases}
\end{equation*}
It is easy to check that $f_\delta \in W^{2,p}_0(\br)$ and, using density arguments, it will serve as a test function. Now, following same estimate as Theorem~\ref{diff-crit-hardy}, we have 
    \begin{align*}
        \lim_{\delta\rightarrow 0^+} \int_{\br} \frac{|f_\delta|^{p}|\gradg \varrho|^p}{\varrho^{Q}(\ln \frac{R}{\varrho})^{p}} \dx\dt = +\infty.
    \end{align*}    
    Also, we deduce
    \begin{align*}
        &\int_{\br}\frac{\left|\rlapg f_\delta\right|^p}{|\gradg \varrho|^p\varrho^{Q-2p}}\dx\dt\\&\leq \bigg(\frac{(p-1)(Q-2)}{p}\bigg)^p\int_{B^Q_{\frac{R}{2}}(o)\setminus B^Q_\delta(o)}\frac{|f_\delta|^p}{\varrho^{Q}\big(\ln\frac{R}{\varrho}\big)^{p}}|\gradg \varrho|^p\dx\dt\\&+\left(\frac{p^2-1}{p^2}\right)^p\int_{B^Q_{\frac{R}{2}}(o)\setminus B^Q_\delta(o)}\bigg(\ln\frac{R}{\varrho}\bigg)^{-1-p}\varrho^{-Q}|\gradg \varrho|^p\dv\\&+(\ln 2)^{p-1}2^pR^{-p}(Q-1)^p\int_{B^Q_R(o)\setminus B^Q_{\frac{R}{2}}(o)}\varrho^{-Q+p}|\gradg \varrho|^p\dv.
    \end{align*}
    Immediately, we received that the last term of the above estimate is bounded, and also, we have
    \begin{align*}
        \int_{B^Q_{\frac{R}{2}}(o)\setminus B^Q_\delta(o)}\bigg(\ln\frac{R}{\varrho}\bigg)^{-1-p}\varrho^{-Q}|\gradg \varrho|^p\dv&\leq \mathcal{I}_{\frac{n+p}{2}-1}\int_{\delta}^{\frac{R}{2}}\left(\ln\frac{R}{\varrho}\right)^{-1-p}\frac{1}{\varrho}\dr\\&=-\frac{\mathcal{I}_{\frac{n+p}{2}-1}}{p}\left[\left(\ln\delta\right)^{-p}-(\ln 2)^{-p}\right]=O(1),
    \end{align*}
    for $\delta\rightarrow 0^+$. This gives, in combining that
    \begin{align*}
       \lim_{\delta\rightarrow 0^+}\frac{\int_{\br}\frac{\left|\rlapg f_\delta\right|^p}{|\gradg \varrho|^p\varrho^{Q-2p}}\dx\dt}{\int_{\br}\frac{|f_\delta |^p}{\varrho^{Q}\big(\ln\frac{R}{\varrho}\big)^{p}}|\gradg \varrho|^p\dx\dt} =\bigg(\frac{(Q-2)(p-1)}{p}\bigg)^p.
    \end{align*}
    \end{proof}

\section{Interplay of Critical and Subcritical Hardy inequalities}\label{sh-ch-rel}
This section is devoted to studying the equivalence between the subcritical and critical weighted Hardy inequalities for the radial Grushin operator. In the latter part of the section, we improve both the critical and subcritical Hardy inequalities by adding positive remainder terms, which arise from a careful application of weighted Caffarelli--Kohn--Nirenberg (CKN) inequalities in the Grushin setting. We begin with the proof of Theorem~\ref{main-thm-1}, which illustrates the connection between the critical and subcritical weighted Hardy inequalities.
\begin{proof}[Proof of Theorem~\ref{main-thm-1}]
	Let us recall the polar coordinate representation. For any, $(x,t)\in \rno$ we can write it as 
	\begin{align*}
		(x,t)= (\varrho, \phi, \theta_1 ,...,\theta_{n-1}) \text{ where} \, \phi, \theta_i\in(0,\pi) \, \text{for}\, i=1,\cdots, (n-2), \,\text{and}  \,\theta_{n-1}\in (0,2\pi).
	\end{align*}
	Also for any $(y,s)\in\brp\subset\rmo$, we write
	\begin{align*}
		(y,s)= (\tvr, \tilde{\phi}, \tilde{\theta_1} ,...,\tilde{\theta}_{m-1}) \text{ where} \, \tilde{\phi}, \tilde{\theta_i}\in(0,\pi) \, \text{for}\, i=1,\cdots, (m-2), \,\text{and}  \,\tilde{\theta}_{m-1}\in (0,2\pi).
	\end{align*}
	
	Now, for a given $w\in C_{0}^\infty(\brp \setminus \{o'\})$, we define a new function $u\in C^\infty_{0}(\rno\setminus\{o\})$ as follows 
	\begin{align*}
		u(\varrho, \phi, \theta_1,...,\theta_{n-1})=u(\varrho, \phi, \theta_{n-m+1},...,\theta_{n-1})
		:=w(\tvr, \tilde{\phi}, \tilde{\theta_1} ,...,\tilde{\theta}_{m-1})
	\end{align*}
	where, 
	\begin{align*}
		\varrho=\left(\ln\frac{R}{\tvr}\right)^{-\frac{1}{\alpha}}, \text{ with } \alpha=\left(\frac{Q-Q'-\beta}{b-1}\right) \text{ and } \tilde{\phi}=\phi, \, \tilde{\theta_i}=\theta_{n-m+i} \text{ for } 1\leq i\leq m-1.
	\end{align*}
	Next, by polar coordinate decomposition, the above change of variable, independence of $u$ for the variables $\theta_1,\ldots,\theta_{n-m}$, and volume element representation \eqref{volume}, we deduce
	\begin{align*}
		&\int_{\rno}\frac{|\rgradg u(x,t)|^{Q'}}{\varrho^{\beta}(x,t)}\dx\dt\\
		&=\frac{1}{2}\int_{\Omega} \int_0^{\infty}(\sin\phi)^{\frac{m}{2}}|\pd u|^{Q'}\varrho^{Q-\beta-1}\dr\dsn\\
		&=\left(\prod_{i=1}^{n-m}\int_{\theta_i=0}^{\pi}(\sin\theta_{i})^{n-1-i}{\rm d}\theta_{i}\right)\times\left(\frac{1}{2} \int_{\Omega '}\int_0^R (\sin\phi ')^{\frac{m}{2}} |\pdp w|^{Q'}\bigg|\frac{\drt}{\dr}\bigg|^{Q'-1}\left(\ln\frac{R}{\tvr}\right)^{-\frac{Q-\beta-1}{\alpha}} \drt \dsnp\right)\\
		&=\frac{\omega_n}{\omega_m}\frac{1}{2}\int_{\Omega '}\int_0^R(\sin\phi ')^{\frac{m}{2}}|\pdp w|^{Q'}\left(\alpha \tvr\right)^{Q'-1}\left(\ln\frac{R}{\tvr}\right)^{-\frac{Q-\beta-1-(\alpha+1)(Q'-1)}{\alpha}}\ds\dsnp\\
		&=\frac{\omega_n}{\omega_m}\alpha^{Q'-1} \int_{\brp} |\rgradgp w(y,s)|^{Q'}\left(\ln\frac{R}{\tvr(y,s)}\right)^{Q'-b}\dy\ds,
	\end{align*}
		where we used $-(Q-\beta-1)+(\alpha +1)(Q'-1)=(Q'-b)\alpha$. In the same way as earlier, we have
		\begin{align*}
		&	\int_{\rno}\frac{|u(x,t)|^{Q'}|\gradg \varrho(x,t)|^{Q'}}{\varrho^{Q' + \beta}(x,t)}\dx\dt\\&
			=\frac{1}{2}\int_{\Omega} \int_0^{\infty}(\sin\phi)^{\frac{m}{2}}| u|^{Q'}\varrho^{Q-Q'-\beta-1}\dr\dsn\\&=\left(\prod_{i=1}^{n-m}\int_{\theta_i=0}^{\pi}(\sin\theta_{i})^{n-1-i}{\rm d}\theta_{i}\right)\times\left(\frac{1}{2} \int_{\Omega '}\int_0^R (\sin\phi ')^{\frac{m}{2}} | w|^{Q'}\left(\ln\frac{R}{\tvr}\right)^{-\frac{Q-Q'-\beta+\alpha}{\alpha}} (\alpha\tilde{\varrho})^{-1}\drt \dsnp\right)\\&
		=\frac{\omega_n}{\omega_m}\alpha^{-1}\int_{\brp}\frac{|w(y,s)|^{Q'}|\gradgp \tvr (y,s)|^{Q'}}{\tvr (y,s)^{Q'}\left( \ln\frac{R}{\tvr (y,s)}\right)^{b}}\dy \ds.
	\end{align*}
	Hence, by combining these two, we will get the result.
\end{proof}

In the non-radial function setting, we have observed that the critical Hardy inequality in lower dimensions is equivalent to the subcritical Hardy inequality in higher dimensions. However, the converse does not hold: for non-radial functions, it is not possible to embed the subcritical Hardy inequality (in higher dimensions) into the critical one (in lower dimensions), as certain information associated with the spherical variables may be lost. In contrast, for radial functions, the converse also holds. Moreover, in the radial case, we observe a distinct behavior of the constants in both directions, a phenomenon that highlights a key difference between the Grushin and Euclidean settings. We present this result next.

 \begin{theorem}\label{w-r-t}
	Suppose $G=\rno$ and $G'=\rmo$ are the Grushin spaces. Assume $Q=n+2>Q'=m+2\geq 3$ and radius $R>0$. Suppose $\beta,\, b\in \mathbb{R}$ with $\beta<Q-Q'$ and $Q'\leq b$. Then for any $w\in C_{0,rad}^\infty(\brp \setminus \{o'\})$ (resp. $u\in C^\infty_{0,rad}(\rno \setminus \{o\})$), there exists $u\in C^\infty_{0,rad}(\rno\setminus\{o\})$ (resp. $w\in C^\infty_{0,rad}(\brp \setminus \{o'\})$), and there hold 
	\begin{align*}
		\int_{\rno}&\frac{|\rgradg u(x,t)|^{Q'}}{\varrho ^{\beta}(x,t)}\dx\dt - \bigg(\frac{Q-Q'-\beta}{Q'}\bigg)^{Q'}\int_{\rno}\frac{|u(x,t)|^{Q'}|\gradg\varrho(x,t)|^{Q'}}{\varrho^{Q'+\beta}(x,t)}\dx\dt \\&
		=\frac{C^{\Omega,m}}{C^{\Omega',m}}\left(\frac{Q - Q' - \beta}{b - 1}\right)^{Q' - 1}\bigg[
		\int_{\brp} \frac{|\rgradgp w(y,s)|^{Q'}}{\left(\ln \frac{R}{\tilde{\varrho}(y,s)}\right)^{b - Q'}}\dy\ds \\
		&\quad - \left(\frac{b - 1}{Q'}\right)^{Q'}\int_{\brp} \frac{|w(y,s)|^{Q'}|\gradgp \tilde{\varrho}(y,s)|^{Q'}}{\tilde{\varrho}^{Q'}(y,s)\left(\ln \frac{R}{\tilde{\varrho}(y,s)}\right)^b}\dy\ds
		\bigg],
	\end{align*}
    where the integral constants $C^{\Omega,m}$ and $C^{\Omega',m}$ are defined in notation subsection.
\end{theorem}
\begin{proof}
	We start with the homogeneous norm $\varrho$ and $\tilde{\varrho}$ on $\rno$ and $\rmo$ respectively. For a given $w\in C_{0,rad}^\infty(\brp \setminus \{o'\})$ (resp. $u\in C^\infty_{0,rad}(\rno\setminus \{o\})$), we define $u\in C^\infty_{0,rad}(\rno\setminus\{o\})$ (resp. $w\in C_{0,rad}^\infty(\brp \setminus \{o'\})$) by 
	\begin{equation*}
		u(\vr)= w(\tvr),
	\end{equation*}
	where $\vr$ and $\tvr$ satisfies the following transformation
	\begin{align}\label{w-r}
		\varrho=\left(\ln\frac{R}{\tvr}\right)^{-\frac{1}{\alpha}}, \quad \text{ where } \alpha=\left(\frac{Q-Q'-\beta}{b-1}\right).
	\end{align}
	Then, using polar coordinate decomposition and a change of variable \eqref{w-r}, we have    
	\begin{align*}
		\int_{\mathbb{R}^{n+1}}\frac{|\rgradg u(x,t)|^{Q'}}{\varrho^{\beta}(x,t)}\dx\dt
		&=\frac{C^{\Omega,m}}{2}\int_0^{\infty}|\pd u(\varrho)|^{Q'}\varrho^{Q-\beta-1}\dr\\
		=&\frac{C^{\Omega,m}}{2}\alpha^{Q'-1}\int_0^{R}|\pdp w(\tvr)|^{Q'}\tvr^{Q'-1}\left(\ln\frac{R}{\tvr}\right)^{-\frac{(Q-\beta-1)-(\alpha +1)(Q' -1)}{\alpha}}\drt\\
		=&\frac{C^{\Omega,m}}{C^{\Omega',m}}\alpha^{Q'-1}\int_{\brp}|\rgradgp w(y,s)|^{Q'}\left(\ln\frac{R}{\tvr(y,s)}\right)^{Q'-b} \dy \ds,
	\end{align*}
	where we used $-(Q-\beta-1)+(\alpha +1)(Q'-1)=(Q'-b)\alpha$. Similarly, for the other term, we deduce
	\begin{align*}
		\int_{\rno}\frac{|u(x,t)|^{Q'}|\gradg\varrho(x,t)|^{Q'}}{\varrho^{Q'+\beta}(x,t)}\dx\dt=\frac{C^{\Omega,m}}{C^{\Omega',m}}\alpha^{-1}\int_{\brp}\frac{|w(y,s)|^{Q'}|\gradgp \tvr (y,s)|^{Q'}}{\tvr (y,s)^{Q'}\left( \ln\frac{R}{\tvr (y,s)}\right)^{b}}\dy \ds
	\end{align*}    
	Therefore, by combining both terms, we get the desired result.
\end{proof}

In the remainder of this section, we aim to improve \eqref{lim-tw-wg-hardy-eq-intro} by adding positive remainder terms. This improvement arises from a delicate application of Theorem~\ref{w-r-t} together with the weighted CKN inequalities on the Grushin space established by Song and Li \cite{song}.

In special parameters, CKN inequalities \cite[Theorem~1.5]{song} reads as follows: Let $G''=\re^{q+1}$ with $q\geq 1$, with $Q''=q+2$. Suppose $a\in (0,1)$, $p\gamma<p-Q''$, and $2\leq p\leq Q''$. Then there exist a positive constant $C_{CKN}=C_{CKN}(Q'',p,\gamma,a)$ such that for all $v\in C_0^\infty(\rko\setminus\{o''\})$ there holds
 \begin{equation}\label{CKN-ineq}
            \|\btvr^{\frac{p-Q''}{p}}\nabla_{G''} v\|_{{L^p}(\rko)}^a \,\|\btvr^{\frac{p-Q''}{p}} v\|_{{L^p}(\rko)}^{1-a} \geq  C_{CKN} \|\bar{\psi}^{\frac{a-\gamma}{2}}\btvr^{\gamma} v\|_{L^{\frac{Q''}{1-a-\gamma}}(\rko)}.
\end{equation}

First, using this, we establish an improvement of the subcritical Hardy inequality for radial functions, presented as a preparatory result below.
\begin{proposition}\label{i-r}
    Let $G=\rno$ be a Grushin space with homogeneous dimension $Q=n+2\geq 3$. Assume that $q\in \mathbb{N}$, $a\in (0,1)$, and $2\leq p<\min\{Q, Q''\}$, where $Q''$ is the homogeneous dimension of $\mathbb{R}^{q+1}$. Suppose the following conditions hold: $p\gamma < p - Q''$ and $a - \gamma - 2 \geq -q$. Further, let $\beta < Q - p$ and define
    \begin{align*}
       \delta = Q'' - Q + \frac{Q'}{1 - a - \gamma} \left( \gamma + \frac{Q - p - \beta}{p} \right). 
    \end{align*}
   Then, for every $u\in C_{0,\mathrm{rad}}^\infty(\rno\setminus\{o\}) \setminus \{0\}$, the following inequality holds
    \begin{multline*}
        \int_{\rno} \frac{|\rgradg u(x,t)|^p}{\varrho^{\beta}(x,t)} \dx\dt - \left( \frac{Q - p - \beta}{p} \right)^p \int_{\rno} \frac{|u(x,t)|^p |\gradg \varrho(x,t)|^p}{\varrho^{p+\beta}(x,t)} \dx\dt\\
        \geq C \frac{\left( \int_{\rno} |u(x,t)|^{\frac{Q''}{1 - a - \gamma}} \varrho^{\delta}(x,t) \dx\dt \right)^{\frac{p(1 - a - \gamma)}{Q'' a}}}{\left( \int_{\rno} |u(x,t)|^p \varrho^{-\beta}(x,t) \dx\dt \right)^{\frac{1 - a}{a}}},
    \end{multline*}
    where $C$ is a constant depending on $Q''$, $Q$, $p$, $\gamma$, and $a$.
\end{proposition}
\begin{proof}
Take $u\in C_{0,rad}^\infty(\rno\setminus\{o\})\setminus\{0\}$ and define $v(\varrho(x,t)):= \varrho^{\frac{Q-p-\beta}{p}}(x,t)u(\varrho(x,t))$ for $(x,t)\in\rno$.  By using the polar coordinate decomposition, we obtain 
\begin{equation}\label{hardy-remain}
	\begin{split}
		J(u):=&\int_{\m}\frac{|\rgradg u|^p}{\varrho^\beta}\dv - \bigg(\frac{Q-p-\beta}{p}\bigg)^p\int_{\m}\frac{|u|^p|\gradg \varrho|^p}{\varrho^{p+\beta}}\dv\\&=\frac{1}{2}\int_{\sn}\int_0^\infty\bigg|\frac{\partial v}{\partial \varrho} \varrho^{-\frac{Q-p-\beta}{p}}-\bigg(\frac{Q-p-\beta}{p}\bigg)v \varrho^{-\frac{Q-\beta}{p}}\bigg|^p\varrho^{Q-1-\beta}\psi^{\frac{p}{2}-1}\dr\dsn\\&-\frac{1}{2}\bigg(\frac{Q-p-\beta}{p}\bigg)^p\int_{\sn}\int_0^\infty|u|^p \varrho^{Q-1-p-\beta}\psi^{\frac{p}{2}-1}\dr\dsn.
	\end{split}
\end{equation}

Now applying Lemma \ref{cpl}, for real valued case and the expression of $u$, we deduce 
\begin{align*}
	&\bigg|\frac{\partial v}{\partial \varrho} \varrho^{-\frac{Q-p-\beta}{p}}-\bigg(\frac{Q-p-\beta}{p}\bigg)v \varrho^{-\frac{Q-\beta}{p}}\bigg|^p\\&\geq \bigg(\frac{Q-p-\beta}{p}\bigg)^p|u|^p\varrho^{-p}-p\bigg(\frac{Q-p-\beta}{p}\bigg)^{p-1}|v|^{p-2}v\frac{\partial v}{\partial \varrho} \varrho^{-Q+\beta+1}+c_p\bigg|\frac{\partial v}{\partial \varrho}\bigg|^p\varrho^{-Q+p+\beta}.
\end{align*}

Integration by parts and boundary condition of the compactly supported smooth function $v$  gives 
\begin{align*}
	\frac{1}{2}\int_{\sn}\int_0^\infty|v|^{p-2}v\frac{\partial v}{\partial \varrho} \varrho^{-Q+\beta+1}\varrho^{Q-1-\beta} \psi^{\frac{p}{2}-1}\dr\dsn=\frac{1}{2}\int_{\sn}\int_0^\infty\frac{\partial |v|^p}{\partial \varrho} \psi^{\frac{p}{2}-1}\dr\dsn=0.
\end{align*}

Therefore, we deduce by summing all these terms as follows
\begin{align}\label{hardy-remain-last}
	J(u)\geq  \frac{c_p}{2}\int_{\sn}\int_0^
	\infty \bigg|\frac{\partial v}{\partial \varrho}\bigg|^p\varrho^{p-1}\psi^{\frac{p}{2}-1}\dr\dsn= c_p\int_{\m}\big|\rgradg v\big|^p\varrho^{p-Q}\dv.
\end{align}

 Now, again using the polar coordinate decomposition, we deduce
        \begin{align*}
    J(u)&\geq \frac{c_p}{2}\int_{\sn}\int_0^
    \infty \bigg|\frac{\partial v}{\partial \varrho}\bigg|^p\varrho^{p-1}\psi^{\frac{p}{2}-1}\dr\dsn\\&=c_p\frac{C^{\Omega,p-2}}{C^{\Omega'',p-2}}\int_{\mathbb{R}^{q+1}}|\nabla_{G''} v(\btvr(z,r))|^p|\btvr(z,r)|^{p-Q'}\,{\rm d}z \,{\rm d}r
        \end{align*}
where we apply the change of variable $\varrho=\btvr$ which is homogeneous distance on $\mathbb{R}^{q+1}$ on the variable $(z,r)$. Now, applying CKN inequality \eqref{CKN-ineq} for the function $v$, and further change of variable and simplification gives
\begin{align*}
     &J(u)\geq c_p\frac{C^{\Omega,p-2}}{C^{\Omega'',p-2}}       C_{CKN}\frac{\bigg(\int_{\mathbb{R}^{q+1}}\bar{\psi}^{\frac{a-\gamma}{2}}\btvr^{\frac{\gamma Q''}{1-a-\gamma}}|v|^{\frac{Q''}{1-a-\gamma}}\,{\rm d}z\,{\rm d}r\bigg)^{\frac{p(1-a-\gamma)}{Q''a}}}{\bigg(\int_{\mathbb{R}^{q+1}}\btvr^{p-Q''}|v|^p\,{\rm d}z\,{\rm d}r\bigg)^{\frac{1-a}{a}}}\\&\geq c_p\frac{C^{\Omega,p-2}}{C^{\Omega'',p-2}} C_{CKN}\frac{\left(C^{\Omega'',a-\gamma-2}\int^{\infty}_0 |u(\varrho)|^{\frac{Q''}{1-a-\gamma}} \varrho^{\frac{\gamma Q''}{1-a-\gamma}+\frac{(Q-p-\beta)Q''}{p(1-a-\gamma)}+Q''-Q}\varrho^{Q-1}\dr\right)^{\frac{p(1-a-\gamma)}{Q''a}}}{\bigg(C^{\Omega'',-1}\int_0^{\infty} |u(\varrho)|^p \varrho^{-\beta} \varrho^{Q-1}\dr\bigg)^{\frac{1-a}{a}}}\\&=C\frac{\big(\int_{\rno}|u(x,t)|^{\frac{Q''}{1-a-\gamma}}\varrho^{\delta}(x,t)\dx\dt\big)^{\frac{p(1-a-\gamma)}{Q''a}}}{\big(\int_{\rno}|u(x,t)|^p\varrho^{-\beta}(x,t)\dx\dt\big)^{\frac{1-a}{a}}},
\end{align*}
where $C=c_p\frac{C^{\Omega,p-2}}{C^{\Omega'',p-2}}C_{CKN}\frac{\left(C^{\Omega'',a-\gamma-2}\right)^{\frac{p(1-a-\gamma)}{Q''a}}}{\left(C^{\Omega'',-1}\right)^{\frac{1-a}{a}}}\frac{\left(C^{\Omega,-1}\right)^{-\frac{p(1-a-\gamma)}{Q''a}}}{\left(C^{\Omega,-1}\right)^{\frac{a-1}{a}}}$ and this completes the result.
\end{proof}

Now, as a straightforward consequence of Proposition~\ref{i-r} with $p=Q'$, and by subsequently applying this result in Theorem~\ref{w-r-t}, we establish the existence of a nontrivial remainder term in the critical Hardy inequality. For brevity, we omit the details.
\begin{theorem}
  Suppose $G=\rno$ and $G'=\rmo$ are the Grushin spaces. Assume that $Q=n+2>Q'=m+2\geq 3$ and let the radius $R>0$. Let $q\in\mathbb{N}$ be such that $m<q$. Suppose $\beta, \, b\in \re$ satisfy $0<\beta<Q'$, $Q'\leq b$, $\beta<Q-Q'$, and $\gamma Q'<Q'-(q+2)$. Define
  \begin{align*}
      \delta = (q+2-Q)+\frac{(q+2)(\gamma Q'+Q-Q'-\beta)}{(Q'-\beta-\gamma Q')}, \quad  M=1+\frac{(b-1)(Q+\delta)}{(Q-Q'-\beta)}, \quad  
        N=\frac{(b\beta-bQ+Q')}{(Q-Q'-\beta)}.
  \end{align*}
  Then, for every $u\in C_{0,rad}^\infty(\brp\setminus\{o'\})\setminus\{0\}$, there exists a constant $C=C(Q,Q',q,\beta,\gamma)>0$ such that the following inequality holds
  \begin{multline*}
        \int_{\brp} \frac{|\rgradgp u(y,s)|^{Q'}}{\big(\ln\frac{R}{\tvr(y,s)}\big)^{b-Q'}}\dy\ds - \left(\frac{b-1}{Q'}\right)^{Q'}\int_{\brp} \frac{|u(y,s)|^{Q'}|\gradg \tvr(y,s)|^{Q'}}{\tvr^{Q'}(y,s)\big(\ln\frac{R}{\tvr(y,s)}\big)^{b}}\dy\ds\\
       \geq  C\frac{\left(\int_{\brp}|u(y,s)|^{\frac{Q'(q+2)}{(Q'-\beta-\gamma Q')}}\tvr(y,s)^{-Q'}\left(\ln\frac{R}{\tvr (y,s)}\right)^{-M}\dy\ds\right)}{\left(\int_{\brp}|u(y,s)|^{Q'}\tvr(y,s)^{-Q'}\left(\ln \frac{R}{\tvr(y,s)}\right)^{N}\dy\ds\right)}.
  \end{multline*}
\end{theorem}

\begin{remark}
It is worth noting that the author in \cite{st-cvpde} employs the radialization method to obtain analogous results for non-radial functions in the Euclidean setting. For further details on this method, we refer the reader to \cite[Lemma~4.1]{pems} and \cite[Equation~18, p.~69]{st-cvpde} for general $p$. However, this method is not available in the Grushin setting due to the presence of the sub-Riemannian geometric term $\psi$ in the polar coordinate decomposition. As a consequence, extending the earlier results to non-radial functions remains an open problem, as we are currently able to establish them only for radial functions.
\end{remark}

\section{Extremizer Stability Analysis}\label{ex-st-an}
The primary objective of this section is to demonstrate stability. First, we will prove Theorem~\ref{stab-sub-hardy-th}, which establishes the stability of the weighted Grushin Hardy inequality.
\begin{proof}[Proof of Theorem~\ref{stab-sub-hardy-th}]
Let us begin with any  function $u\in C_0^\infty(\m\setminus\{o\})$, and define the following transformation, $v(x,t):=\varrho^{\frac{Q-p-\beta}{p}}u(x,t)$ which is also a member of $C_0^\infty(\m\setminus\{o\})$. Recall the definition of $J(u)$ in \eqref{hardy-remain} and then by \eqref{hardy-remain-last} we deduce
\begin{align*}
    J(u)\geq  c_p\int_{\m}\big|\rgradg v\big|^p\varrho^{p-Q}\dv.
\end{align*}
Now, in the above using \eqref{diff-crit-hardy-eq} for $\gamma=p$ and for the function $v$, we get 
\begin{align*}
    J(u)&\geq c_p \frac{(p-1)^p}{p^p}\int_{\m}\frac{|v(x,t)-v_R(x,t)|^p}{\varrho^{Q}(x,t)\big|\ln\frac{R}{\varrho(x,t)}\big|^{p}}|\gradg \varrho(x,t)|^p\dv\\&=c_p \frac{(p-1)^p}{p^p}\int_{\m}\frac{|\varrho^{\frac{Q-p-\beta}{p}}u(x,t)-R^{\frac{Q-p-\beta}{p}}u\big(R\frac{x}{\varrho(x,t)},R^2\frac{t}{\varrho^2(x,t)}\big)|^p}{\varrho^{Q}(x,t)\big|\ln\frac{R}{\varrho(x,t)}\big|^{p}}|\gradg \varrho(x,t)|^p\dv,
\end{align*}
for any $R>0$. After simplifying it 
\begin{align*}
    J(u)\geq c_p \frac{(p-1)^p}{p^p}d_H(u,R)^p
\end{align*}
and taking the supremum over $R>0$, we deduce the desired result.
\end{proof}

In the same spirit, we can discuss the stability of the critical Hardy inequality (\eqref{lim-tw-wg-hardy-eq} with $b=p$) for the Grushin operator, and for brevity, we are concentrating on the case of unit radius. To do this, let us define the following distance function. For a function  $u\in C_0^\infty(\ba\setminus\{o\})$ and for a fixed real number $R>0$, we define the following distance function
\begin{multline*}
    d_C(u,R):= \bigg(\int_{\ba}\frac{|u(x,t)-\big(\ln \frac{1}{\varrho(x,t)}\big)^{\frac{p-1}{p}}R^{\frac{p-1}{p}}u\big(e^{-\frac{1}{R}}\frac{x}{\varrho(x,t)},e^{-\frac{2}{R}}\frac{t}{\varrho^2(x,t)}\big)|^p}{\varrho^Q(x,t)\big(\ln \frac{1}{\varrho(x,t)}\big)^{p}\big|\ln \big(R\ln \frac{1}{\varrho(x,t)}\big)\big|^{p}}|\gradg \varrho(x,t)|^p\dv\bigg)^{\frac{1}{p}}.
\end{multline*}
Now we have the following stability result in the critical setup:
\begin{theorem}\label{stab-th-crt-hardy}
Assume $2\leq p<\infty$. Let $G$ be a Grushin space with dimension $Q$ with $Q\geq 3$. Then for $u\in C_0^\infty(\ba\setminus\{o\})$ there holds 
\begin{multline*}
    \int_{\ba} \frac{|\rgradg u(x,t)|^p}{\varrho^{Q-p}(x,t)} \dx\dt-  \bigg(\frac{p-1}{p}\bigg)^p\int_{\ba} \frac{|u(x,t)|^{p}|\gradg \varrho(x,t)|^p}{\varrho^{Q}(x,t)(\ln \frac{1}{\varrho(x,t)})^{p}} \dx\dt\\ \geq c_p\bigg(\frac{p-1}{p}\bigg)^p \sup_{R>0} d_C(u,R)^p,
\end{multline*}
where $c_p$ is defined in Lemma~\ref{cpl}.   
\end{theorem}
\begin{proof}
    Let us start with any functions $u\in C_0^\infty(\ba\setminus\{o\})$, we define the new function in terms of polar coordinates by $g(s,\sigma):=s^{\frac{p-1}{p}}u(e^{-s^{-1}},\, \sigma)$, where $s\in (0,\infty)$ and $\sigma\in \Omega$. We notice $g(o)=\lim_{s\rightarrow 0} g(s,\sigma)=u(o)\lim_{s\rightarrow 0}s^{\frac{p-1}{p}}=0$ and as $u$ has compact support, so $g\in C_0^\infty(\m\setminus\{o\})$. The rest of the proof will be exactly similar to the case of stability of the subcritical Hardy inequality. Also, refer to the proof in \cite[Theorem~3.4]{RSY} for details. We omit the details for brevity.
\end{proof}

In the same manner, we can have the extremizer stability of the Critical Rellich Inequality Theorem~\ref{rellich-crit} in the case of unit radius. This can be stated as follows:
 \begin{theorem}\label{stab-th-crt-rellich}
     Let $G$ be a Grushin space with dimension $Q$ with $Q\geq 3$. Assume $\ba$ be the $\varrho-$gauge ball with radius $1$. Let $2\leq p<Q$. Then for any $u\in C_0^\infty (\ba\setminus\{o\})$ there holds
        \begin{multline*}
            \int_{\ba}\frac{\left|\rlapg u(x,t)\right|^p}{|\gradg \varrho(x,t)|^p\varrho^{Q-2p}(x,t)}\dx\dt - \bigg(\frac{(Q-2)(p-1)}{p}\bigg)^p\int_{\ba} \frac{|u(x,t)|^{p}|\gradg \varrho(x,t)|^p}{\varrho^{Q}(x,t)(\ln \frac{1}{\varrho(x,t)})^{p}} \dx\dt\\ \geq c_p\bigg(\frac{(Q-2)(p-1)}{p}\bigg)^p \sup_{R>0} d_C(u,R)^p ,
        \end{multline*}
        where $c_p$ is defined in Lemma~\ref{cpl}
 \end{theorem}
 \begin{proof}
Using Theorem~\ref{Hardy--Rellich} for the weight $\beta$ replaced by $Q-2p$, and then applying Theorem~\ref{stab-th-crt-hardy}, we derive the desired result.
\end{proof}

We now proceed to prove Theorem~\ref{stab-th-higher-hr}, following the approach from~\cite{sano-mia}. We incorporate a second-order operator with a special parameter and, by carefully applying integration by parts and the polar coordinate decomposition, we achieve the desired result. The details are presented below.

\begin{proof}[Proof of Theorem~\ref{stab-th-higher-hr}]
We start with $u\in C_0^\infty(\m\setminus\{o\})$, and then we define the following transformation, $v(x,t):=\varrho^{\frac{Q-kp-\beta}{p}}u(x,t)$ which is also a member of $C_0^\infty(\m\setminus\{o\})$. First, consider the even case $k=2l$ for $l\geq 1$. By Theorem~\ref{high-hr}, part (i), we can write
\begin{align}\label{h0}
            \int_{\rno}\frac{\left|\rlapg^l u(x,t)\right|^p}{|\gradg \varrho(x,t)|^{(2l-1)p}\varrho^{\beta}(x,t)}\dx\dt&= \int_{\rno}\frac{\left|\rlapg^{l-1}\left(\frac{\rlapg u(x,t)}{|\gradg \varrho(x,t)|^2}\right)\right|^p}{|\gradg \varrho(x,t)|^{(2l-3)p}\varrho^{\beta}(x,t)}\dx\dt \nonumber \\&\geq C_{k-2,p,\beta}\int_{\rno}\frac{|\rlapg u(x,t)|^p}{\varrho^{(k-2)p+\beta}(x,t)|\gradg \varrho(x,t)|^p}\dx\dt.
\end{align}

After a simple computation, we can write 
\begin{equation*}
-\rlapg u = \psi \varrho^{-\frac{Q+(2-k)p-\beta}{p}}\left(\Lambda_{k}^{\frac{1}{p}}(Q,p,\beta)\, v -\varrho ^2 \Delta_{\theta_k}v\right),
\end{equation*}
where $\Lambda_{k}(Q,p,\beta)$ is defined in \eqref{big-lam}, and
\begin{align*}
   \Delta_{\theta_k} v= \frac{\partial^2 v}{\partial \varrho^2}+ \frac{\theta_k -1}{\varrho} \frac{\partial v}{\partial \varrho}\quad \text{ with } \quad \theta_k=2k+\frac{Q(p-2)+2\beta}{p}.
\end{align*}

Now, let us define the following deficit term by 
\begin{align*}
   H(u):= \int_{\rno}\frac{|\rlapg u(x,t)|^p}{\varrho^{(k-2)p+\beta}(x,t)|\gradg \varrho(x,t)|^p}\dv-\Lambda_{k}(Q,p,\beta) \int_{\rno}\frac{|u(x,t)|^p|\gradg \varrho(x,t)|^p}{\varrho^{kp+\beta}(x,t)}\dv. 
\end{align*}
Therefore, using the transformation, polar coordinates, Lemma~\ref{cpl}, and neglecting the non-negative term, we have
\begin{align}\label{h1}
    \nonumber H(u)&=\frac{1}{2}\int_{\Omega}\int_{0}^{\infty}\left[ \left|\Lambda_{k}^{\frac{1}{p}}(Q,p,\beta)\, v -\varrho ^2 \Delta_{\theta_k}v\right|^p-\Lambda_k(Q,p,\beta)|v|^p\right]\psi^{\frac{p}{2}-1} \varrho^{-1}\dr\dsn \nonumber\\&\geq -\Lambda_{k}^{\frac{p-1}{p}}(Q,p,\beta)\frac{p}{2}\int_{\Omega}\int_{0}^{\infty}|v|^{p-2}v(\Delta_{\theta_k}v)\psi^{\frac{p}{2}-1} \varrho\dr\dsn \nonumber \\
    &=-\Lambda_{k}^{\frac{p-1}{p}}(Q,p,\beta)\frac{p}{2}\int_{\Omega}\int_{0}^{\infty}|v|^{p-2}v\, \frac{\partial^2 v}{\partial \varrho^2}\, \psi^{\frac{p}{2}-1} \varrho\dr\dsn,
\end{align}
where in the end line we used the expression of $\Delta_{\theta_k}v$ and  $\int_{\Omega}\int_{0}^{\infty}|v|^{p-2}v(\pd v)\psi^{\frac{p}{2}-1}\dr\dsn=0$.

Furthermore, applying integration by parts w.r.t. radial component, we have
\begin{align*}
    -\frac{1}{2}\int_{\Omega}\int_{0}^{\infty}|v|^{p-2}v\, \frac{\partial^2 v}{\partial \varrho^2}\, \psi^{\frac{p}{2}-1} \varrho\dr\dsn&=\frac{2(p-1)}{p^2}\int_{\Omega}\int^{\infty}_0 |\pd(|v|^{\frac{p-2}{2}}v)|^2\varrho \, \psi^{\frac{p}{2}-1}\dr\dsn\\&   =\frac{4(p-1)}{p^2}\|\rgradg(|v|^\frac{p-2}{2}v\psi^{\frac{p-2}{4}})\varrho^{\frac{2-Q}{2}}\|^2_{L^2(\mathbb{R}^{n+1})}.
\end{align*}

Now applying Theorem~\ref{diff-crit-hardy} for $p=\gamma=2$ and for the function $|v|^\frac{p-2}{2}v\psi^{\frac{p-2}{4}}$ in above, then using it \eqref{h1} and writing back in terms of $u$, and using $\psi(Rx/\varrho^2, R^2t/\varrho^2)=\psi(x,t)$, we deduce 
\begin{align*}
    H(u)\geq \Lambda_{k}^{\frac{p-1}{p}}(Q,p,\beta)\frac{4(p-1)}{p}\|\rgradg(|v|^\frac{p-2}{2}v\psi^{\frac{p-2}{4}})\varrho^{\frac{2-Q}{2}}\|^2_{L^2(\mathbb{R}^{n+1})} \geq \Lambda_{k}^{\frac{p-1}{p}}(Q,p,\beta)\frac{(p-1)}{p} d_{R}(u,k,\beta)^2,
\end{align*}
where distance function $d_{R}(u,k,\beta)$ is defined in \eqref{dist-high}. Now we finish this part by using the above, \eqref{h0}, $2l=k$, and \eqref{sharp-constant}  that
\begin{align*}
      &\int_{\rno}\frac{\left|\rlapg^l u(x,t)\right|^p}{|\gradg \varrho(x,t)|^{(2l-1)p}\varrho^{\beta}(x,t)}\dx\dt-C_{k,p,\beta}\int_{\rno}\frac{|u(x,t)|^p|\gradg \varrho(x,t)|^p}{\varrho^{2lp+\beta}(x,t)}\dx\dt\\&\geq  C_{k-2,p,\beta}\int_{\rno}\frac{|\rlapg u(x,t)|^p}{\varrho^{(k-2)p+\beta}(x,t)|\gradg \varrho(x,t)|^p}\dx\dt-C_{k,p,\beta}\int_{\rno}\frac{|u(x,t)|^p|\gradg \varrho(x,t)|^p}{\varrho^{kp+\beta}(x,t)}{\rm d}x{\rm d}t\\&\geq H(u)\,C_{k-2,p,\beta}+\left(\Lambda_{k}(Q,p,\beta)C_{k-2,p,\beta}-C_{k,p,\beta}\right) \int_{\rno}\frac{|u(x,t)|^p|\gradg \varrho(x,t)|^p}{\varrho^{kp+\beta}(x,t)}{\rm d}x{\rm d}t\\& \geq C_{k-2,p,\beta}\Lambda_{k}^{\frac{p-1}{p}}(Q,p,\beta)\frac{(p-1)}{p} d_{R}(u,k,\beta)^2,
\end{align*}
and by taking $C=C_{k-2,p,\beta}\Lambda_{k}^{\frac{p-1}{p}}(Q,p,\beta)\frac{(p-1)}{p} $ and supremum over $R>0$, part (i) follows.

Next, we will deal with the odd case $k=2l+1$ with $l\geq 1$. First, from \eqref{sub-hardy-eq-intro}, we notice
\begin{align*}
     \int_{\rno}\frac{\left|\rgradg \rlapg^l u(x,t)\right|^p}{|\gradg \varrho(x,t)|^{2lp}\varrho^{\beta}(x,t)}\dx\dt\geq \Lambda_1(Q,p,\beta)\int_{\rno}\frac{\left|\rlapg^l u(x,t)\right|^p}{|\gradg \varrho(x,t)|^{(2l-1)p}\varrho^{p+\beta}(x,t)}\dx\dt.
\end{align*}
Next in the r.h.s. of the above using part (i) with weight $p+\beta$ instead of $\beta$, we deduce
 \begin{multline*}
                \int_{\rno}\frac{\left|\rlapg^l u(x,t)\right|^p}{|\gradg \varrho(x,t)|^{(2l-1)p}\varrho^{p+\beta}(x,t)}\dx\dt\\ \geq C_{2l,p,p+\beta}\int_{\rno}\frac{|u(x,t)|^p|\gradg \varrho(x,t)|^p}{\varrho^{2lp+\beta+p}(x,t)}\dx\dt+ C'\sup_{R>0}d_{R}(u,2l,p+\beta)^2,
\end{multline*}  
where $C'=C_{2l-2,p,p+\beta}\Lambda_{2l}^{\frac{p-1}{p}}(Q,p,p+\beta)\frac{(p-1)}{p}$. Therefore, for $k=2l+1$, combining the last two estimates, we deduce
\begin{align*}
     \int_{\rno}\frac{\left|\rgradg \rlapg^l u(x,t)\right|^p}{|\gradg \varrho(x,t)|^{2lp}\varrho^{\beta}(x,t)}\dx\dt &\geq \Lambda_1(Q,p,\beta) C_{2l,p,p+\beta}\int_{\rno}\frac{|u(x,t)|^p|\gradg \varrho(x,t)|^p}{\varrho^{2lp+\beta+p}(x,t)}\dx\dt\\&+\Lambda_1(Q,p,\beta)C'\sup_{R>0}d_{R}(u,2l,p+\beta)^2.
\end{align*}
Finally, we notice from the definition of distance functions and the defined constants in \eqref{big-lam-1}, \eqref{big-lam}, and \eqref{sharp-constant} that
\begin{align*}
    d_{R}(u,2l,p+\beta)=d_{R}(u,2l+1,\beta) \quad \text{ and } \quad C_{2l+1,p,\beta}=\Lambda_1(Q,p,\beta) C_{2l,p,p+\beta}.
\end{align*}
Therefore, part (ii) is completed by considering the generic constant $C=\Lambda_1(Q,p,\beta)C'$.
\end{proof}

\section{Higher order Hardy--Rellich identity}\label{hr-r-id}
Up to this point, we have established various Hardy--Rellich-type results for higher-order iterative Grushin radial operators in the $L^p(G)$ setting. Unfortunately, similar to the pointwise comparison between the gradient and the radial gradient, a pointwise comparison between these higher-order iterative radial operators and their original (non-radial) counterparts is not available in the literature. In the Grushin setting, a comparison can be made in terms of the $L^2$ norm, as established in \cite[Remark~5.2]{GJR}. This is stated as follows: let $Q \geq 4$ and $u \in C_0^\infty(\rno \setminus {o})$, then \begin{align}\label{rad_lap_comp}
 \int_{\rno}\frac{|\lapg u|^2}{|\gradg \varrho|^2}\dx\dt \geq \int_{\rno}\frac{|\rlapg u|^2}{|\gradg \varrho|^2}\dx\dt.
 \end{align}
Such comparison results with more general weights have been well studied in the context of Riemannian model manifolds whose sectional curvature is bounded above by minus one (see \cite[Lemma~6.1]{jmaa}). Following this line of reasoning, we now aim to improve \eqref{rad_lap_comp} by establishing a weighted version.
\begin{lemma}\label{rad_lap}
	Let $G$ be a Grushin space with dimension $Q$ with $Q\geq 4$. Assume $\alpha\in \mathbb{R}$ and $-2\leq \alpha<Q-4$. Then for all $u\in C_0^\infty(\rno\setminus\{o\})$, there holds 
	\begin{equation}\label{eq_rad_lap}
		\int_{\rno}\frac{|\lapg u|^2}{|\gradg \varrho|^2\varrho^\alpha}\dx\dt \geq \int_{\rno}\frac{|\rlapg u|^2}{|\gradg \varrho|^2\varrho^\alpha}\dx\dt.
	\end{equation}
	Moreover, equality is held when $u$ is a radial function.		
\end{lemma}	
\begin{proof}		
	Consider any $u\in C_0^\infty(\rno\setminus\{o\})$, then we know	
	\begin{align*}
		(\lapg u)^2  = \psi^2 \left( \partial^2_{\varrho}u + \frac{(Q-1)}{\varrho} \partial_{\varrho}u \right)^2  + \psi^2
		\frac{16}{\varrho^4} (\mathcal{L}_{\sigma} u)^2  + 2\psi^2 \left( \partial^2_{\varrho}u + \frac{(Q-1)}{\varrho} \partial_{\varrho}u \right) \frac{4}{\varrho^2} (\mathcal{L}_{\sigma} u) u,
	\end{align*}
	and
	\begin{align*}
		(\rlapg u)^2 & = \psi^2 \left( \partial^2_{\varrho}u + \frac{(Q-1)}{\varrho} \partial_{\varrho}u \right)^2\,.
	\end{align*}
	
	Now, by Grushin's spherical harmonics decomposition, we write 
	\begin{align*}
		u(x,t):= u(\varrho, \sigma) =\sum_{k=0}^{\infty}d_{k}(\varrho)\Phi_k(\sigma).
	\end{align*}
	
	This decomposition of $u$ gives, after using polar coordinate and orthonormal properties of $\{\Phi_k\}$, the following
	\begin{align*}
		&\int_{\rno}\frac{|\lapg u|^2}{|\gradg \varrho|^2\varrho^\alpha}\dx\dt - \int_{\rno}\frac{|\rlapg u|^2}{|\gradg \varrho|^2\varrho^\alpha}\dx\dt\\&=\frac{1}{2}  \sum_{k = 0}^{\infty} \bigg[ 
		16 \lambda_{k}^2\int_0^\infty d_{k}^2 \varrho^{n-3-\alpha} \dr - 8\lambda_k \int_0^\infty  \left( d_{k}^{\prime \prime} + \frac{(Q-1)}{\varrho}d_{k}^{\prime} \right) d_k \varrho^{n-1-\alpha}\dr\bigg],
	\end{align*}
	
	Hence, using this to prove \eqref{eq_rad_lap}, it is enough to show that 
	\begin{align*}
		R_k:=2\lambda_k\int_{\rno}\frac{d_k^2|\gradg\varrho|^2}{\varrho^{\alpha+4}}\dx\dt -  \int_{\rno}
		\frac{d_k(\rlapg d_k)}{\varrho^{\alpha+2}}\dx\dt\geq 0\quad \text{  for all } k\geq 1.
	\end{align*}
	
	Let us fix $k$. Now notice that we have $d_k(\rlapg d_k)=\frac{1}{2}\rlapg(d_k^2)-|\rgradg d_k|^2$, and using \eqref{sub-hardy-eq-intro} for $p=2$, $u=d_k$, $\beta=\alpha+2$ we deduce 
	\begin{align*}
		R_k\geq \left[2\lambda_k+\frac{(Q-\alpha-4)^2}{4}\right]\int_{\rno}\frac{d_k^2|\gradg\varrho|^2}{\varrho^{\alpha+4}}\dx\dt -  \frac{1}{2}\int_{\rno}
		\frac{\rlapg(d_k^2)}{\varrho^{\alpha+2}}\dx\dt.
	\end{align*}
	
	Now, the integration by parts formula yields
	\begin{align*}
		\frac{1}{2}\int_{\rno}
		\frac{\rlapg(d_k^2)}{\varrho^{\alpha+2}}\dx\dt&= \frac{1}{2}\int_{\Omega}\left(\int_{0}^\infty \left[ \left(\partial^2_{\varrho}d_k^2\right)\varrho^{n-\alpha-1} + (Q-1) \left(\partial_{\varrho}d_k^2\right) \varrho^{n-\alpha-2}\right]\dr \right)\frac{1}{2\psi}\dsn\\&=-\frac{(\alpha+2)(n-\alpha-2)}{2}\int_{\rno}\frac{d_k^2}{\varrho^{\alpha+4}}\dx\dt.
	\end{align*}
	By exploiting this in the remainder term and using the condition $-2\leq \alpha \leq Q-4$ and the fact $|\gradg\varrho|\leq 1$, it simplifies that for fixed $k\geq 1$, there holds
	\begin{align*}
		R_k\geq \left[2\lambda_k+\frac{(Q-\alpha-4)^2}{4}+\frac{(\alpha+2)(n-\alpha-2)}{2}\right]\int_{\rno}\frac{d_k^2|\gradg\varrho|^2}{\varrho^{\alpha+4}}\dx\dt.
	\end{align*}
	This essentially gives it is enough to check that for $k\geq 1$, there holds
	\begin{align*}
		2\lambda_k+\frac{(Q-\alpha-4)^2}{4}+\frac{(\alpha+2)(n-\alpha-2)}{2}\geq 0.
	\end{align*}
	Therefore, by using $\lambda_k=\frac{k(k+n)}{4}$ and $Q=n+2$ we deduce
	\begin{align*}
		2\lambda_k+\frac{(Q-\alpha-4)^2}{4}+\frac{(\alpha+2)(n-\alpha-2)}{2}=\frac{k(k+Q-2)}{2}+\frac{(Q-\alpha-4)(Q+\alpha)}{4}.
	\end{align*}
	Now, using the conditions $ -Q\leq -2\leq \alpha<Q-4$, we immediately comment that the above algebraic term is non-negative, and this concludes the proof.		
\end{proof}	

All the $L^p(G)$ Hardy and $L^2(G)$ Rellich inequalities, as well as their improvements discussed in this article for the radial Grushin operator, can immediately be viewed as improved versions of the existing results for the full Grushin operator. This follows from the pointwise comparison between the radial Grushin gradient and the full Grushin gradient, and from Lemma~\ref{rad_lap}. However, such comparisons are not available for higher-order operators. Even in the Euclidean setting, where the Grushin structure is absent, analogous pointwise comparisons between $\lapg^m$ and $\rlapg^m$ for integers $m \geq 2$ are not available in the literature.

From this perspective, it remains an interesting question to find an identity or remainder term for the Hardy inequality \eqref{sub-hardy-eq-intro}, or for higher-order cases in the general $L^p(G)$ setting when $p \neq 2$, even in the Euclidean case. However, by using spherical harmonic decomposition and exploiting the Hilbert space structure of $L^2(G)$, an explicit identity is available for $p = 2$ (see \cite{GJR}) for both the Hardy and Rellich inequalities in the Grushin setting. A crucial ingredient in this development was the recent introduction of spherical derivatives in the study of inequalities for the Grushin setting, as established in \cite{GJR}, which are not available in the context of homogeneous groups.

In the case $p = 2$, Huang and Ye \cite{he} recently established higher-order Hardy--Rellich identities for polyharmonic radial operators in the Euclidean setting. In a similar spirit, we investigate analogous phenomena in the Grushin setting. It is worth noting that these results essentially provide explicit remainder terms for Theorem~\ref{high-hr} in the case $p = 2$ only. This restriction arises because our approach fundamentally relies on the Hilbert space structure of $L^2$, which is absent for $p \neq 2$. These technical limitations currently prevent the extension of such functional identities to the general case $p > 1$.

For $\beta\in \mathbb{R}$, let us define a weighted $L^2$-norm as follows:
\begin{align*}
    \|u\|_{\beta}^2:= \int_{\mathbb{R}^{n+1}} \frac{|u|^2}{\varrho^{\beta}} \dx \dt,
\end{align*}
and recall the weighted first-order differential operator defined in \eqref{tbeta} in the following form:
\begin{align*}
 T_{\beta} (u):= \psi^{1/2} \left(\partial_{\varrho}u+\frac{Q-\beta-2}{2\varrho}u \right) .
\end{align*}

\begin{lemma}\label{lemma-T-alpha}
 Let $\beta\in \mathbb{R}$. Then for $u\in C_{0}^{\infty}(\mathbb{R}^{n+1})$ there holds
 \begin{align*}
     \|T_{\beta}u\|_{\beta}^2= \|\psi^{1/2} \partial_{\varrho}u\|_{\beta}^2- \frac{(Q-\beta-2)^2}{4} \|\psi^{1/2} u\|_{\beta+2}^2
 \end{align*}
\end{lemma}
\begin{proof}
First, we write
\begin{align}\label{T-alpha}
 \|T_{\beta}u\|_{\beta}^2
 =  \int_{\mathbb{R}^{n+1}} \frac{\psi \left| \partial_{\varrho}u\right|^2} {\varrho^{\beta}}\dv  + (Q-\beta-2)\int_{\mathbb{R}^{n+1}} \frac{\psi ~ \partial_{\varrho}u \cdot u }{\varrho^{\beta+1}}\dv
   + \frac{(Q-\beta-2)^2}{4} \int_{\mathbb{R}^{n+1}} \frac{\psi u}{\varrho^{\beta+2}}\dx \dt   
\end{align}
Now, writing the integral into polar coordinates and using the integration by parts in the $\varrho$ variable, we get
\begin{align}\label{int-T-alpha}
(Q-\beta-2)\int_{\mathbb{R}^{n+1}} \frac{\psi ~ \partial_{\varrho}u \cdot u }{\varrho^{\beta+1}}\dv
   = -\frac{(Q-\beta-2)^2}{2} \int_{\mathbb{R}^{n+1}} \frac{\psi ~  |u|^2  }{\varrho^{\beta+2}}\dv.
\end{align}
Finally, putting the estimate \eqref{int-T-alpha} in \eqref{T-alpha}, we get the desired result. 
\end{proof}

Let $\langle \cdot \rangle_{\beta}$ be an inner product defined by
\begin{align*}
  \langle u,v \rangle_{\beta} =\int_{\mathbb{R}^{n+1}} \frac{u(x,t)v(x,t)}{\varrho^{\beta}(x,t)} \dx \dt.
\end{align*} 
Then we have the following equality 
\begin{align}\label{parial-deri-innprdt}
    \langle \partial_{\varrho}u, v\rangle_{\beta}= - \langle u, \partial_{\varrho} v \rangle_{\beta} - (Q-\beta-1) \langle u,v \rangle_{\beta +1},
\end{align}
for any $u, v \in  C_0^\infty(\rno\setminus\{o\})$.
\begin{lemma}
    Let $\beta, \alpha, \gamma \in \mathbb{R}$. Then, for $u \in C_0^\infty(\rno\setminus\{o\})$, we have
    \begin{align}
        &\|T_{\alpha}(u)\|_{\beta}^{2}= \|T_{\beta}(u)\|_{\beta}^{2}+ \frac{(\alpha-\beta)^2}{4} \|\psi^{1/2}u\|_{\beta+2}^2, \label{T-alpha-beta}\\
        & T_{\beta}(\varrho^{\gamma}u)= \varrho^{\gamma} T_{\beta-2\gamma}(u), \quad \|T_{\beta}(\varrho^{\gamma}u)\|_{\beta}^2=\|T_{\beta-2\gamma}(u)\|_{\beta-2\gamma}^2, \label{T-alpha-gamma}
    \end{align}
    and 
    \begin{align}
         \langle T_{\beta}(\partial_{\varrho}u), T_{\beta+2}(u) \rangle_{\beta+1}=\frac{(4+\beta-Q)}{2} \|T_{\beta+2}(u)\|_{\beta+2}^2. \label{T-alpha-inner-product}
    \end{align}
\end{lemma}
\begin{proof}
The heart of the proof of this lemma is the following equality: for any $\eta \in C_{0}^{\infty}(0,\infty)$ and $\alpha \in \mathbb{R}$, there holds
\begin{align}\label{functional-equa-Huang}
    \int_{0}^{\infty} |\eta^{'}(s)+\frac{M_1}{s}\eta|^2 ~ s^{\alpha} {\rm d}s - \int_{0}^{\infty} |\eta^{'}(s)+\frac{M_2}{s}\eta|^2 ~ s^{\alpha} {\rm d}s
    =  (M_2-M_1)(M_2+M_1-\alpha+1) \int_{0}^{\infty} |\eta|^2 ~ s^{\alpha-2} {\rm d}s.
\end{align}
One can prove the above equality \eqref{functional-equa-Huang} by just using integration by parts appropriately for each term of the left-hand side of \eqref{functional-equa-Huang}. See also equation no. (2.5) of \cite{he}.

Note that for $u\in C_{0}^{\infty}(\mathbb{R}^{n+1}\setminus \{0\})$, we write
\begin{align*}
    \|T_{\sigma}(u)\|_{\beta}^2= \int_{\Omega} \int_{0}^{\infty} \left| \partial_{\varrho}u+ \frac{Q-2-\sigma}{2\varrho} u\right|^2 \frac{\varrho^{n+1-\beta}}{2} \dr \dsn.
\end{align*}
Therefore, using \eqref{functional-equa-Huang} with $M_1= \frac{Q-2-\alpha}{2}$, $M_2= \frac{Q-2-\beta}{2}$ and $\alpha= n+1-\beta$, we get
\begin{align*}
        \|T_{\alpha}(u)\|_{\beta}^2- \|T_{\beta}(u)\|_{\beta}^2
        =& \frac{(\beta-\alpha)}{2} \left(\beta-\frac{\alpha+\beta}{2}\right) \int_{\Omega} \int_{0}^{\infty} \left| u\right|^2 \frac{\varrho^{n-1-\beta}}{2} {\rm d}\varrho {\rm d}\Omega
        =\frac{(\alpha-\beta)^2}{4} \int_{\mathbb{R}^{n+1}} \psi ~\frac{|u|^2}{\varrho^{\beta+2}} {\rm d}x{\rm d}t,
\end{align*}
which completes the proof of \eqref{T-alpha-beta}. The proof of \eqref{T-alpha-gamma} is straightforward and follows from the definition of $T_{\beta}$ and $\|\cdot\|_{\beta}$. To prove \eqref{T-alpha-inner-product}, we first observe that
\begin{align}\label{T-alpha-deri}
    \partial_{\varrho}(\langle T_{\beta+2}(u) \rangle)=T_{\beta}(\partial_{\varrho}u)- \frac{1}{\varrho} T_{\beta+2}(u).
\end{align}
Using the above observation and \eqref{parial-deri-innprdt}, we get
\begin{align*}
 \langle T_{\beta}(\partial_{\varrho}u), T_{\beta+2}(u) \rangle_{\beta+1}
 =& \langle \partial_{\varrho} T_{\beta+2}(u)+ \frac{1}{\varrho}T_{\beta+2}(u), T_{\beta+2}(u) \rangle_{\beta+1}\\
 =& \langle \partial_{\varrho} T_{\beta+2}(u), T_{\beta+2}(u) \rangle_{\beta+1} + \|T_{\beta+2}(u)\|_{\beta+2}^2\\
 = & -\frac{(Q-\beta-2)}{2} \|T_{\beta+2}(u)\|_{\beta+2}^2 + \|T_{\beta+2}(u)\|_{\beta+2}^2\\
 =&\frac{(4+\beta-Q)}{2} \|T_{\beta+2}(u)\|_{\beta+2}^2,   
\end{align*}
which proves \eqref{T-alpha-inner-product}.
\end{proof}

An important consequence of \eqref{T-alpha-beta} is the following equality
\begin{align}\label{consequence}
    \|\mathcal{L}_{\varrho, G}u\|_{\beta}^2= \|T_{-Q}(\psi^{1/2}\partial_{\varrho}u)\|_{\beta}^2
    = \|T_{\beta}(\psi^{1/2}\partial_{\varrho}u)\|_{\beta}^2 + \frac{(\beta+Q)^2}{4} \|(\psi\partial_{\varrho}u)\|_{\beta+2}^2.
\end{align}

Next, we define higher-order operators which will appear in the remainder terms for higher-order Hardy--Rellich inequalities in Theorem \ref{main-thm-4}. For $\beta \in \mathbb{R}$, $k\in \mathbb{N}\cup \{0\}$, the following properties of $\mathcal{R}_{\beta,k}$ (defined in \eqref{remainder-R}) are useful in proving our main results.
\begin{lemma}
For any $\beta, \gamma \in \mathbb{R}$, $k\in \mathbb{N}$ and $u \in C_0^\infty(\rno\setminus\{o\})$, there holds
\begin{align}
    &\mathcal{R}_{\beta,k}(\varrho^{\gamma}u)= \varrho^{\gamma} \mathcal{R}_{\beta-2\gamma,k}(u),\label{R-gamma}\\
    & \|\mathcal{R}_{\beta,k}(\psi^{1/2}\partial_{\rho}u)\|_{\beta}^2= \frac{(Q-2k-\beta-4)^2}{4} \|\mathcal{R}_{\beta+2,k}(\psi^{1/2}u)\|_{\beta+2}^2 + \|\mathcal{R}_{\beta,k+1}(u)\|_{\beta}^2, \label{R-rho} \\
    \text{and} \quad & \langle \mathcal{R}_{\beta,k}(\partial_{\rho}u), \mathcal{R}_{\beta+2,k}(u) \rangle_{\beta+1}= \frac{(\beta+2k+4-Q)}{2} \|\mathcal{R}_{\beta+2,k}(u)\|_{\beta+2}^2. \label{R-alpha-innprd}
\end{align}
\end{lemma}
\begin{proof}
The proof of \eqref{R-gamma} is obvious and follows from repeated applications of \eqref{T-alpha-gamma}. Using induction in $k$, we will prove \eqref{R-rho} and \eqref{R-alpha-innprd}. First, we will prove \eqref{R-alpha-innprd}. 

The equality \eqref{R-alpha-innprd} for $\beta=0$ is follows from \eqref{T-alpha-inner-product}. Let us assume that \eqref{R-alpha-innprd} is true for $k=l\in \mathbb{N}$ i.e. 
  \begin{align}\label{R-induction-l}
   \langle \mathcal{R}_{\beta,l}(\partial_{\rho}u), \mathcal{R}_{\beta+2,l}(u) \rangle_{\beta+1}= \frac{(\beta+2l+4-Q)}{2} \|\mathcal{R}_{\beta+2,l}(u)\|_{\beta+2}^2.   
  \end{align}
Now, using \eqref{T-alpha-deri} and \eqref{R-gamma}, we get
\begin{align}\label{R-inter}
  \mathcal{R}_{\beta, l}(\partial_{\varrho}T_{\beta+2l+4}(u))= \mathcal{R}_{\beta, l}(T_{\beta+2l+2}(\partial_{\varrho}u))-\mathcal{R}_{\beta, l}(\varrho^{-1}T_{\beta+2l+4}(u)) 
   =\mathcal{R}_{\beta, l+1}(\partial_{\varrho}u)-\frac{1}{\varrho}\mathcal{R}_{\beta+2, l+1}(u)
\end{align}

Hence, the equality \eqref{R-induction-l} and \eqref{R-inter} together imply
\begin{align}\label{R-innprdt-intermediate}
    &\langle \mathcal{R}_{\beta,l+1}(\partial_{\rho}u), \mathcal{R}_{\beta+2,l+1}(u) \rangle_{\beta+1}\\
    \nonumber= & \langle \mathcal{R}_{\beta, l}(\partial_{\varrho}T_{\beta+2l+4}(u), \mathcal{R}_{\beta+2, l+1}(u) \rangle_{\beta+1}  + \|\mathcal{R}_{\beta+2,l+1}(u)\|_{\beta+2}^2\\
    \nonumber= & \frac{(\beta+2l+4-Q)}{2} \|\mathcal{R}_{\beta+2,l+1}(u)\|_{\beta+2}^2 + \|\mathcal{R}_{\beta+2,l+1}(u)\|_{\beta+2}^2\\
    \nonumber=& \frac{(\beta+2(l+1)+4-Q)}{2} \|\mathcal{R}_{\beta+2,l+1}(u)\|_{\beta+2}^2.
\end{align}
This proves \eqref{R-alpha-innprd} for $k=l+1$. Hence, by induction \eqref{R-alpha-innprd} is true for all $k\in \mathbb{N}.$
Next, we prove \eqref{R-rho}. Using Lemma~\ref{lemma-T-alpha} and equality \eqref{T-alpha-deri}, \eqref{T-alpha-inner-product}, we get
\begin{align}\label{R-alpha-inter-1}
 \|\mathcal{R}_{\beta,1}(u)\|_{\beta}^2=& \|T_{\beta}\circ T_{\beta+2}(u)\|_{\beta}^2 \\ 
 \nonumber= & \|\psi^{1/2} \partial_{\varrho}T_{\beta}(u)\|_{\beta}^{2}-\frac{(Q-\beta-2)^2}{4} \|\psi^{1/2}T_{\beta+2}(u)\|_{\beta+2}^2\\
 \nonumber =& \|T_{\beta}(\psi^{1/2} \partial_{\varrho}u)-\frac{1}{\varrho}T_{\beta+2}(\psi^{1/2}u)\|_{\beta}^2- \frac{(Q-\beta-2)^2}{4} \|T_{\beta+2}(\psi^{1/2}u)\|_{\beta+2}^2\\
 \nonumber =& \|T_{\beta}(\psi^{1/2} \partial_{\varrho}u)\|_{\beta}^2+ \|T_{\beta+2}(\psi^{1/2} u)\|_{\beta+2}^2- 2 \langle T_{\beta}(\psi^{1/2} \partial_{\varrho}u,  T_{\beta+2}(\psi^{1/2} u) \rangle\\
 \nonumber & \quad - \frac{(Q-\beta-2)^2}{4} \|T_{\beta+2}(\psi^{1/2}u)\|_{\beta+2}^2\\
 \nonumber = & \|T_{\beta}(\psi^{1/2} \partial_{\varrho}u)\|_{\beta}^2 - \frac{(Q-\beta-4)^2}{4} \|T_{\beta+2}(\psi^{1/2}u)\|_{\beta+2}^2,
\end{align}
which proves \eqref{R-rho} for $k=0$. Let us assume that \eqref{R-rho} is true for $k=l$. Therefore, we write
\begin{align*}
    &\|R_{\beta, l+2}(u)\|_{\beta}^2=\|R_{\beta, l+1}(T_{\beta+2l+4}(u))\|_{\beta}^2\\
    \nonumber=&  \|R_{\beta, l}(\partial_{\varrho}(T_{\beta+2l+4}(\psi^{1/2}u)))\|_{\beta}^2-\frac{(Q-2l-\beta-4)^2}{4} \|R_{\beta+2, l}(T_{\beta+2l+4}(\psi^{1/2}u))\|_{\beta+2}^2.
\end{align*}
Using the equalities \eqref{R-alpha-inter-1}, \eqref{R-inter} and \eqref{R-innprdt-intermediate}, we get
\begin{align*}
 &\|R_{\beta, l+2}(u)\|_{\beta}^2\\
    \nonumber=&  \|R_{\beta, l+1}(\psi^{1/2}\partial_{\varrho}(u))-\varrho^{-1}R_{\beta+2, l+1}(\psi^{1/2}u)\|_{\beta}^2\\
    & \quad -\frac{(Q-2l-\beta-4)^2}{4} \|R_{\beta+2, l}(T_{\beta+2l+4}(\psi^{1/2}u))\|_{\beta+2}^2\\
    =& \|R_{\beta, l+1}(\psi^{1/2}\partial_{\varrho}(u))\|_{\beta}^2 + \|R_{\beta+2, l+1}(\psi^{1/2}u) \|_{\beta+2}^2 -2 \langle R_{\beta, l+1}(\psi^{1/2}\partial_{\varrho}(u)), R_{\beta+2, l+1}(\psi^{1/2}u)\rangle \\
    & \quad -\frac{(Q-2l-\beta-4)^2}{4} \|R_{\beta+2, l+1}(\psi^{1/2}u)\|_{\beta+2}^2\\
    = & \|R_{\beta, l+1}(\psi^{1/2}\partial_{\varrho}(u))\|_{\beta}^2 -\frac{(Q-2l-2-\beta-4)^2}{4} \|R_{\beta+2, l+1}(\psi^{1/2}u)\|_{\beta+2}^2.
\end{align*}
This implies \eqref{R-rho} is true for $k=l+1$. Hence, \eqref{R-rho} is true for all $k\in \mathbb{N}$. 

This completes the proof of the lemma. 
\end{proof}

Now, we will prove an iterative identity for $\|\mathcal{R}_{\beta, k}(\rlapg u)\|_{\beta}^2$ in the following lemma.

Let 
\begin{align}\label{Dm-constant}
    A_{m}= \frac{(Q+m)^2+ (Q-4-m)^2}{4}, \quad m \in \mathbb{R}.
\end{align}
\begin{lemma}\label{R-alpha-laplacian}
    Let $\beta \in \mathbb{R}$ and $k \in \mathbb{N}$. For $u \in C_0^\infty(\rno\setminus\{o\})$, there holds
    \begin{align}\label{eq:R-alpha-laplacian}
     &\|\mathcal{R}_{\beta, k}(\rlapg u)\|_{\beta}^2\\
    \nonumber =& \Lambda_{2}(Q,2,\beta+2k+2) \|\mathcal{R}_{\beta+4,k}(\psi u)\|_{\beta+4}^2 + A_{\beta+ 2k+2}\|\mathcal{R}_{\beta+2,k+1}(\psi^{1/2}u)\|_{\beta+2}^2 + \|\mathcal{R}_{\beta,k+2}(u)\|_{\beta}^2,  
    \end{align}
    where $\Lambda_{2}(Q,2,\beta+2k+2)$ and $A_{\beta+2k+2}$ are defined as in \eqref{big-lam}  and \eqref{Dm-constant} , respectively.
\end{lemma}
\begin{proof}
 Note that we can write as
 \begin{align*}
     \rlapg (u)= \psi \varrho^{1-Q}\partial_{\varrho}(\varrho^{Q-1}\partial_{\varrho} u).
 \end{align*}
 Therefore, using the above equality together with \eqref{R-gamma} and \eqref{R-rho}, we get
 \begin{align*}
     &\|\mathcal{R}_{\beta, k}(\rlapg u)\|_{\beta}^2\\
     =&\|\mathcal{R}_{\beta, k}(\psi \varrho^{1-Q}\partial_{\varrho}(\varrho^{Q-1}\partial_{\varrho} u) )\|_{\beta}^2\\
     =& \| \mathcal{R}_{\beta+2Q-2, k}(\psi^{1/2}\partial_{\varrho}(\psi^{1/2}\varrho^{Q-1}\partial_{\varrho} u) )\|_{\beta+2Q-2}^2\\
     =& \frac{(Q+2k+\beta+2)^2}{4} \|\mathcal{R}_{\beta+2Q,k}(\psi \varrho^{Q-1}\partial_{\varrho}u)\|_{\beta+2Q}^2 + \|\mathcal{R}_{\beta+2Q-2,k+1}(\psi^{1/2}\varrho^{Q-1}\partial_{\varrho}u)\|_{\beta+2Q-2}^2\\
     = & \frac{(Q+2k+\beta+2)^2}{4} \|\mathcal{R}_{\beta+2,k}(\psi\partial_{\varrho}u)\|_{\beta+2}^2 + \|\mathcal{R}_{\beta,k+1}(\psi^{1/2}\partial_{\varrho}u)\|_{\beta}^2\\
     =& \frac{(Q+2k+\beta+2)^2}{4} \left[ \frac{(Q-2k-\beta-6)^2}{4} \|\mathcal{R}_{\beta+4,k}(\psi u)\|_{\beta+4}^2 + \|\mathcal{R}_{\beta+2,k+1}(\psi^{1/2}u)\|_{\beta+2}^2 \right]\\
      & \quad + \frac{(Q-2k-\beta-6)^2}{4} \|\mathcal{R}_{\beta+2,k+1}(\psi^{1/2}u)\|_{\beta+2}^2 + \|\mathcal{R}_{\beta,k+2}(u)\|_{\beta}^2\\
     =& \Lambda_{2}(Q,2,\beta+2k+2) \|\mathcal{R}_{\beta+4,k}(\psi u)\|_{\beta+4}^2 + A_{\beta+ 2k+2}\|\mathcal{R}_{\beta+2,k+1}(\psi^{1/2}u)\|_{\beta+2}^2 + \|\mathcal{R}_{\beta,k+2}(u)\|_{\beta}^2.
 \end{align*}
 This completes the proof of the lemma.
\end{proof}
\begin{remark}
 For any $m\geq 1$, by iterating \eqref{eq:R-alpha-laplacian}, one can get the following:
 \begin{align}\label{rem:R-alpha-itera}
  \|\mathcal{R}_{\beta, k}(\rlapg^m u)\|_{\beta}^2
    =  \sum_{j=0}^{2m} B_{j,m,\beta,k} &\|\mathcal{R}_{\beta+2j, k+2m-j}( \psi^{j/2}u)\|_{\beta+2j}^2
 \end{align}
 where the coefficients $B_{j,m,\beta,k}$  are given by 
 \begin{align}\label{coeff-R-alpha-ite}
   B_{j,m+1,\beta,k} = \Lambda_{2}(Q,2,\beta+2k+2)B_{j-2,m,\beta+4,k} + A_{\beta+ 2k+2} B_{j-1,m,\beta+2,k+1} + B_{j,m,\beta,k+2}
 \end{align}
 for all $0\leq j \leq 2m+2$ with the convention $B_{0,0,\beta,k}=1$ for every $\beta, \, k$, and $B_{j,m,\beta,k}=0$ if $j<0$ or $j>2m$.
\end{remark}

We now conclude our discussion by presenting the proof of our final result.

\begin{proof}[Proof of Theorem~\ref{main-thm-4}] 
    First, we will prove \eqref{radial-higher-Rellich} by induction on $l$. Recall that $\psi=|\gradg \varrho|^2$. For $l=1$,  by \eqref{consequence} and \eqref{R-alpha-inter-1}, we get
    \begin{align}\label{radial-higher-Rellich-1}
        &\int_{\rno} \frac{\left|\rlapg u(x,t)\right|^2}{|\gradg \varrho(x,t)|^{2}\varrho^{\beta}(x,t)}\dx\dt\\
        \nonumber=&  \int_{\rno} \frac{\left|T_{-Q}(u)(x,t)\right|^2}{\varrho^{\beta}(x,t)}\dx\dt\\
        \nonumber=& \int_{\rno} \frac{\left|T_{\beta}( \partial_{\varrho} u)(x,t)\right|^2}{\varrho^{\beta}(x,t)}\dx\dt + \frac{(\beta+Q)^2}{4}\int_{\rno} \frac{\left| \partial_{\varrho}u(x,t)\right|^2\psi(x,t)}{\varrho^{\beta+2}(x,t)}\dx\dt\\
       \nonumber = & \int_{\rno} \frac{\left|\mathcal{R}_{\beta,1}(\psi^{-1/2} u)(x,t)\right|^2}{\varrho^{\beta}(x,t)}\dx\dt  + \frac{(Q-\beta-2)^2}{4}\int_{\rno} \frac{\left|T_{\beta+2}( u)(x,t)\right|^2}{\varrho^{\beta+2}(x,t)}\dx\dt \\
       \nonumber & \quad \quad \quad \quad + \frac{(\beta+Q)^2}{4}\int_{\rno} \frac{\left| \partial_{\varrho}u(x,t)\right|^2\psi(x,t)}{\varrho^{\beta+2}(x,t)}\dx\dt\\
        \nonumber=& \Lambda_2(Q,2,\beta) \int_{\rno}
        \frac{|u(x,t)|^2|\gradg \varrho(x,t)|^2}{\varrho^{\beta+4}(x,t)}\dx\dt + A_{\beta} \int_{\rno} \frac{\left|T_{\beta+2}(u)(x,t)\right|^2}{\varrho^{\beta+2}(x,t)}\dx\dt \\
        \nonumber & \quad \quad \quad \quad + \int_{\rno} \frac{\left|\mathcal{R}_{\beta,1}(u)(x,t)\right|^2}{|\gradg \varrho(x,t)|^{2}\varrho^{\beta}(x,t)}\dx\dt.
    \end{align}
    This proves the estimate \eqref{radial-higher-Rellich} for $l=1$. Let us assume that \eqref{radial-higher-Rellich} is true for $l=m\geq 1$. Therefore, we have 
    \begin{align}\label{main-m+1}
         &\int_{\rno}\frac{\left|\rlapg^{m+1} u(x,t)\right|^2}{|\gradg \varrho(x,t)|^{2(2m+1)}\varrho^{\beta}(x,t)}\dx\dt
        = \int_{\rno}\frac{\left|\rlapg^m (\psi^{-1}\rlapg u)(x,t)\right|^2}{|\gradg \varrho(x,t)|^{2(2m-1)}\varrho^{\beta}(x,t)}\dx\dt \\& =  C_{2m,2,\beta}\int_{\rno}\frac{|\rlapg u(x,t)|^2}{|\gradg \varrho(x,t)|^2 \varrho^{\beta+4m}(x,t)}{\rm d}x{\rm d}t +\sum_{j=0}^{2m-1} D_{j,m, \beta} \int_{\rno} \frac{|\mathcal{R}_{\beta+2j, 2m-1-j} (\rlapg u)(x,t)|^2}{|\gradg \varrho(x,t)|^{2(2m-j+1)}\varrho^{\beta+2j}(x,t)}  {\rm d}x{\rm d}t.\nonumber
    \end{align}
    Now, using \eqref{radial-higher-Rellich-1} with $\beta$ replaces with $\beta + 4m$, we get
    \begin{align}\label{ineter-laplacian}
        &\int_{\rno} \frac{\left|\rlapg u(x,t)\right|^2}{|\gradg \varrho(x,t)|^{2}\varrho^{\beta+4m}(x,t)}\dx\dt = \Lambda_2(Q,2,\beta+4m) \int_{\rno}
        \frac{|u(x,t)|^2|\gradg \varrho(x,t)|^2}{\varrho^{\beta+4m+4}(x,t)}\dx\dt \\
        \nonumber& \quad + A_{\beta+4m} \int_{\rno} \frac{\left|\mathcal{R}_{\beta+4m+2,0}(u)(x,t)\right|^2}{\varrho^{\beta+4m+2}(x,t)}\dx\dt 
         + \int_{\rno} \frac{\left|\mathcal{R}_{\beta+4m,1}(u)(x,t)\right|^2}{|\gradg \varrho(x,t)|^{2}\varrho^{\beta+4m}(x,t)}\dx\dt
    \end{align}
    Also, applying Lemma~\ref{R-alpha-laplacian}, for $0\leq j \leq 2m-1$ we obtain
    \begin{align*}
        &\int_{\rno}\frac{|\mathcal{R}_{\beta+2j, 2m-1-j}(\rlapg u)|^2}{|\gradg \varrho(x,t)|^{2(2m-j+1)}\varrho^{\beta+2j}(x,t)} \dx \dt\\
     \nonumber =& \Lambda_{2}(Q,2,\beta+4m) \int_{\rno} \frac{|\mathcal{R}_{\beta+2j+4,2m-1-j}(u)(x,t)|^2}{|\gradg \varrho(x,t)|^{2(2m-1-j)}\varrho^{\beta+2j+4}(x,t)} \dx\dt\\
     \nonumber & + A_{\beta+ 4m}\int_{\rno} \frac{|\mathcal{R}_{\beta+2j+2,2m-j}(u)(x,t)|^2}{|\gradg \varrho(x,t)|^{2(2m-j)}\varrho^{\beta+2j+2}(x,t)} \dx\dt \\
     \nonumber & \quad \quad \quad  + \int_{\rno} \frac{|\mathcal{R}_{\beta+2j,2m-j+1}(u)(x,t)|^2}{|\gradg \varrho(x,t)|^{2(2m-j-1)}\varrho^{\beta+2j}(x,t)} \dx\dt.
    \end{align*}
    Finally, using the equalities \eqref{ineter-laplacian} and \eqref{R-alpha-laplacian} in \eqref{main-m+1} and the fact $\Lambda_{2}(Q,2, \beta+4m)= \Lambda_{2m+2}(Q,2, \beta)$ and iterative relation \eqref{coeff-iteration} of coefficients, we get \eqref{radial-higher-Rellich} for $l=m+1$. Hence, by the principle of induction, we conclude the proof of part (i). 

   For part (ii), the case $k=1$, i.e., $l=0$, follows from \cite[Corollary~3.2]{GJR}. Then, in the rest of the proof, we will prove other cases of \eqref{Higher-order-odd}. Using Lemma~\ref{lemma-T-alpha} and \eqref{radial-higher-Rellich}, we write,
    \begin{align}\label{eq:higher-order-odd-inter}
        &\int_{\rno}\frac{\left|\rgradg \rlapg^l u(x,t)\right|^2}{|\gradg \varrho(x,t)|^{4l}\varrho^{\beta}(x,t)}\dx\dt\\
       \nonumber = &\frac{(Q-\beta-2)^2}{4} \int_{\rno}\frac{\left|\rlapg^l u(x,t)\right|^2}{|\gradg \varrho(x,t)|^{4l-2}\varrho^{\beta+2}(x,t)}\dx\dt + \int_{\rno}\frac{\left|T_{\beta+2} (\rlapg^l u)(x,t)\right|^2}{|\gradg \varrho(x,t)|^{4l}\varrho^{\beta}(x,t)}\dx\dt\\
       \nonumber =& \frac{(Q-\beta-2)^2}{4} C_{2l,2,\beta+2}\int_{\rno}\frac{|u(x,t)|^2|\gradg \varrho(x,t)|^2}{\varrho^{\beta+4l+2}(x,t)}\dx\dt\\
       \nonumber & + \frac{(Q-\beta-2)^2}{4} \sum_{j=0}^{2l-1} D_{j,l, \beta+2} \int_{\rno} \frac{|\mathcal{R}_{\beta+2j, 2l-1-j} (u(x,t))|^2}{|\gradg \varrho(x,t)|^{2(2l-1-j)}\varrho^{\beta+2j+2}(x,t)}  \dx\dt\\
       \nonumber  & \quad \quad  \quad \quad \quad + \int_{\rno}\frac{\left|\mathcal{R}_{\beta+2,0} (\rlapg^l u)(x,t)\right|^2}{|\gradg \varrho(x,t)|^{4l}\varrho^{\beta}(x,t)}\dx\dt.
    \end{align}
    Note that $\frac{(Q-\beta-2)^2}{4} C_{2l,2,\beta+2}= C_{2l+1,2,\beta+2}$. Therefore, using the estimate \eqref{rem:R-alpha-itera} in \eqref{eq:higher-order-odd-inter} and taking $\widetilde{D}_{j,l, \beta}=B_{j,l,\beta, 0}+ \frac{(Q-\beta-2)^2}{4} D_{j-1, l,\beta+2}$, we get the desired equality \eqref{Higher-order-odd}.
\end{proof}

\section*{Acknowledgments} 
The authors are grateful to the anonymous referees for their careful reading of the manuscript and for their valuable comments and suggestions, which have significantly improved the presentation of this work. A.~Banerjee is supported by the doctoral fellowship of the Indian Statistical Institute, Delhi Centre. R.~Basak is supported by the National Science and Technology Council of Taiwan under research grant numbers 113-2811-M-003-007/113-2811-M-003-039. P.~Roychowdhury is partially supported by the MUR-PRIN project No. 20227HX33Z ``Pattern formation in nonlinear phenomena'' granted by the European Union - Next Generation EU. P.~Roychowdhury is also grateful to IIT Hyderabad for its warm hospitality and excellent research environment during the completion of this work.

\section*{Data availability} No new data were collected or generated during this research.

\section*{Declarations} 
{\bf Conflict of interest.} The authors have no conflicts of interest to declare that are relevant to the content of this
article.

\end{document}